\documentclass[, preprint, 11pt]{imsart}

\RequirePackage[OT1]{fontenc}
\RequirePackage{amsthm,amsmath}
\RequirePackage[numbers]{natbib}
\RequirePackage[colorlinks,citecolor=blue,urlcolor=blue]{hyperref}
\usepackage{amsfonts}
\usepackage{color}
\usepackage{pifont}
\usepackage{amssymb}
\usepackage[english]{babel}
\usepackage[T1]{fontenc}
\usepackage{graphicx}
\usepackage{subfigure}
\usepackage{latexsym}
\usepackage{amstext}
\usepackage{mathrsfs}
\usepackage{hyperref}
\usepackage{dsfont}
\usepackage{graphics}
\usepackage{enumerate}
\usepackage{paralist}
\usepackage{color}
\usepackage{fullpage}

\numberwithin{equation}{section}

\newcommand{\beq}{\begin{eqnarray}}
\newcommand{\beqq}{\begin{eqnarray*}}
\newcommand{\eeq}{\end{eqnarray}}
\newcommand{\eeqq}{\end{eqnarray*}}

\newcommand{\rme}{\mathrm{e}}

\usepackage{mathtools}
\usepackage[english]{babel}
\usepackage[utf8]{inputenc}
\usepackage{amsmath}
\usepackage{graphicx}
\usepackage[colorinlistoftodos]{todonotes}
\usepackage{polynom}
\usepackage{amsfonts}
\usepackage{dsfont}
\usepackage{amsthm}
\usepackage{hyperref}
\usepackage{graphicx}
\newtheorem{theorem}{Theorem}[section]

\newtheorem{lemma}{Lemma}[section]
\newtheorem{proposition}[theorem]{Proposition}

\usepackage{fullpage}
\usepackage[T1]{fontenc}
\usepackage{hyperref}
\usepackage{xcolor}
\definecolor{link-color}{rgb}{0.15,0.4,0.15}
\hypersetup{
  colorlinks,
  linkcolor = {link-color}, citecolor = {link-color},
}

\usepackage{graphicx}
\usepackage{hyperref}
\usepackage{amsmath,amsthm,amssymb}
\usepackage{bbm}
\usepackage{tabularx}

\newtheorem{corollary}[theorem]{Corollary}
\theoremstyle{definition}
\newtheorem{remark}[theorem]{Remark}
\newtheorem{definition}[theorem]{Definition}
\newtheorem{example}{Example}[section]

\newcommand{\R}{\mathbb{R}}

\renewcommand{\P}{\mathbb{P}}

\newcommand{\1}{\mathbbm{1}}
\newcommand{\F}{\mathcal{F}}

\newcommand{\bP}{\pmb{\rm P}}

\newcommand{\h}{\mathcal{H}}

\newcommand{\s}{\mathcal{S}}

\renewcommand{\P}{\mathbb{P}}

    \def\d{{\textnormal d}}

\newenvironment{eqnarr}{\begin{IEEEeqnarray}{rCl}}{\end{IEEEeqnarray}\ignorespacesafterend}

\renewcommand{\eqref}[1]{\hyperref[#1]{(\ref*{#1})}}

    \def\beq{\begin{eqnarr}}
    \def\eeq{\end{eqnarr}}
    \def\beqq{\begin{eqnarray*}} 
    \def\eeqq{\end{eqnarray*}} 

        \def\d{{\rm d}}

    \def\d{{\textnormal d}}

\newcommand*{\pref}[1]{\hyperref[#1]{(\ref*{#1})}}
\newcommand*{\refpref}[2]{\hyperref[#2]{\ref*{#1}(\ref*{#2})}}

\definecolor{wco}{rgb}{0.5,0.2,0.3}

\startlocaldefs
\numberwithin{equation}{section}
\theoremstyle{plain}

\endlocaldefs

\begin{document}

\begin{frontmatter}
\title{Entrance laws at the origin of self-similar Markov processes in high dimensions}

\runtitle{Entrance laws at the origin of self-similar Markov processes in high dimension}

\begin{aug}

\author{\fnms{Andreas E. Kyprianou}\thanksref{t3}\thanksref{t2}\ead[label=e1]{a.kyprianou@bath.ac.uk, batisengul@gmail.com}},
\author{\fnms{Victor Rivero}\thanksref{t2}\ead[label=e2]{rivero@cimat.mx}},
\author{\fnms{Bat\i{} \c{S}eng\"ul}\thanksref{t3}\ead[label=e1]{a.kyprianou@bath.ac.uk, batisengul@gmail.com}}
\and
\author{\fnms{Ting Yang}\thanksref{t3}\thanksref{t4}\ead[label=e4]{yangt@bit.edu.cn}}

\thankstext{t3}{Supported by EPSRC grants EP/L002442/1}

\thankstext{t2}{Supported by EPSRC grants EP/M001784/1}

\thankstext{t4}{Supported by NSFC (Grant No. 11501029 and 11731009)}

\affiliation{University of Bath, CIMAT and  University of Bath, Beijing Institute of Technology}

\address{
Andreas Kyprianou
and Bat\i{} \c{S}eng\"ul\\
University of Bath\\
Department of Mathematical Sciences \\
Bath, BA2 7AY\\
 UK\\
\printead{e1}
}

\address{
Victor Rivero\\
CIMAT A. C.\\
Calle Jalisco s/n\\
Col. Valenciana\\
A. P. 402, C.P. 36000\\
Guanajuato, Gto.\\
Mexico\\
\printead{e2}
}
\address{
Ting Yang\\
School of Mathematics and Statistics\\
Beijing Key Laboratory on MCAACI\\
Beijing Institute of Technology\\
Beijing, 100081\\
P. R. China\\
\printead{e4}
}
\end{aug}

\begin{abstract}\hspace{0.1cm}
In this paper we consider the problem of finding entrance laws at the origin for self-similar Markov processes in
$\mathbb{R}^d$,
killed upon hitting the origin. Under suitable assumptions, we show the existence of an entrance law and the convergence to this law when the process is started close to the origin. We obtain an explicit description of the process started from the origin as the time reversal of the original self-similar Markov process conditioned to hit the origin.\end{abstract}



\end{frontmatter}
\setcounter{tocdepth}{1}
\tableofcontents

\part{Entrance laws of self-similar Markov processes}

\section{Introduction}\label{sec:introduction}

Suppose $\h$ is a locally compact
{subset}
of $\R^{d}\setminus\{0\}$ ($d\ge 1$). An $\h$-valued self-similar Markov process (ssMp for short) $(X,\P)=((X_t)_{t \geq 0},\{\P_{z}:z\in\h\})$ is an $\h$-valued c\`adl\`ag Markov process killed at $0$ with $\P_{z}\left(X_{0}=z\right)=1$, which fulfils the scaling property, namely, there exists an $\alpha>0$ such that for any $c>0$,
$$((cX_{c^{-\alpha}t})_{t \geq 0},\P_z) \text{ has the same law as }((X_{t})_{t\ge 0},\P_{cz}) \quad \forall z \in \h.$$
It follows from the scaling property that $\h=c\h$ for all $c>0$. Therefore $\h$ is necessarily a cone of $\R^{d}\setminus\{0\}$ which has the form
$$\h =\phi(\R\times \s)$$
where $\s$ is a locally compact
{subset}
of $\mathbb{S}^{d-1}$ and $\phi$ is the homeomorphism from $\R\times \mathbb{S}^{d-1}$ to $\R^{d}\setminus\{0\}$ defined by $\phi(y,\theta)=\theta\mathrm{e}^{y}$.

\smallskip

The crucial tool in the study of ssMp is the Lamperti-Kiu transform which we now describe. Suppose first that $(X,\P_z)$ is an $\h$-valued ssMp started at $z\in\h$ with index $\alpha>0$ and lifetime $\zeta$, then there exists a Markov additive process (MAP for short, see Section \ref{sec:preliminary} for a rigorous definition) $(\xi,\Theta)$ on $\R \times \s$ started at $(\log \|z\|,\arg(z))$ with lifetime $\zeta_{p}$ such that
  \begin{equation}\label{eq:lamperti_kiu}
   X_t  = \exp\{\xi_{\varphi(t)}\} \Theta_{\varphi(t)}\1_{\{t < \zeta\}} \qquad\forall t \geq 0,
  \end{equation}
   where $\varphi(t)$ is the time-change defined by
  \begin{equation}\label{eq:time_change}
   \varphi(t):=\inf\left\{s>0:\int_0^s \exp\{\alpha\xi_u\}\,{\rm d}u > t\right\},
  \end{equation}
  and $\zeta_{p}=\int_{0}^{\zeta}\|X_{s}\|^{-\alpha}{\rm d}s$.
 We denote the law of $(\xi,\Theta)$ started from $(y,\theta)\in \R\times \s$ by $\bP_{y,\theta}$.
  Conversely given a MAP $(\xi,\Theta)$ under $\bP_{y,\theta}$ with lifetime $\zeta_{p}$, the process $X$ defined by~\eqref{eq:lamperti_kiu} is a ssMp started from $z=\theta {\rm e}^y$ with lifetime $\zeta=\int_{0}^{\zeta_{p}}\mathrm{e}^{\alpha \xi_{s}}{\rm d}s$. Roughly speaking, a MAP is a natural extension of a L\'evy process in the sense that $\Theta$ is an arbitrary well behaved Markov process and $((\xi_{t},\Theta_{t})_{t\ge 0},\bP_{x,\theta})$ is equal in law to $((\xi_{t}+x,\Theta_{t})_{t\ge 0},\bP_{0,\theta})$ for all $x\in \R$ and $\theta\in\s $.
 Whilst MAPs have found a prominent role in e.g. classical applied probability models for queues and dams, c.f. \cite{AsmussenQueue} when $\Theta$ is a Markov chain, the case that $\Theta$ is a general Markov process has received somewhat less attention. Nonetheless a core base of literature exists in the general setting from the 1970s and 1980s thanks to e.g. \cite{Cinlar2,Cinlar1,KaspiWH,Kaspi83}.

  \smallskip

We denote $\h\cup\{0\}$ by $\h_{0}$. In this paper we look for entrance laws of ssMp at the origin,
  that is, the existence of a probability measure $\P_{0}$ such that the extension of $(X,\{\P_{z}:z\in \h_{0}\})$ is self-similar and in particular $\P_{0}=\mbox{w-}\lim_{\h\ni z\to 0}\P_{z}$ in the Skorokhod topology. In Theorem \ref{thm:main} we will prove a general result with as weak assumptions as our study of the underlying MAPs permits. However, the statement of this theorem comes relatively late in this paper on account of the large amount of fluctuation theory we must first develop for general MAPs,
  in order that
  the sufficient conditions make sense. It is quite natural to expect that
 conditions for the existence and stochastic continuity of an entrance law will be highly non-trivial as the process $\Theta$ could essentially take on any role as a regular Markov process.
   Nonetheless, we want to give a flavour of the main results. We give immediately below  the collection of conclusions we are aiming towards, i.e. (C1)-(C5), without addressing the technical assumptions.

\smallskip

The first two conclusions (C1) and (C2) seem rather specialist and pertain to analogues of classical fluctuation results for L\'evy processes, but now in the setting of MAPs. However they hold value in the sense that they provide key building blocks for some of the conclusions lower down.

  \smallskip

{\it {\bf (C1): Conditioning to remain negative:}  There exists a family of probability measures $\hat{\bP}^{\downarrow}=\{\hat{\bP}^{\downarrow}_{y,\theta}: y\le 0,\theta\in\s\}$ such that $((\xi,\Theta),\hat{\bP}^{\downarrow})$ is a right continuous Markov process taking values in $(-\infty,0]\times\s$. Moreover,
        For all $y<0$, $\theta\in\s $, $t\ge 0$ and $\Lambda\in \mathcal{F}_{t}$,
\begin{equation}\nonumber
\hat{\bP}^{\downarrow}_{y,\theta}\left(\Lambda\right)=\lim_{q\to 0+}\hat{\bP}_{y,\theta}\left(\Lambda,t<\mathbf{e}_{q}\,|\,\tau^{+}_{0}>\mathbf{e}_{q}\right),
\end{equation}
where $(\xi, \Theta)$ under $\hat{\bP}_{y,\theta}$ is equal in law to
$(-\xi, \Theta)$,
when $-\xi_0 = y\in\mathbb{R}$,  and $\Theta_0=\theta\in\mathcal{S}$, $\mathbf{e}_{q}$ is an independent and exponentially distributed random variable with parameter $q$ and $\tau^{+}_{0} = \inf\{t>0: \xi_t>0\}$.}

  \smallskip

{\it {\bf (C2): Stationary overshoots and undershoots:} For every $\theta\in \s$, the joint probability measures on $\s\times \R^{-}\times \s\times \R^{+}$
$$\bP_{0,\theta}\left(\Theta_{\tau^{+}_{x}-}\in {\rm d} v,\xi_{\tau^{+}_{x}-}-x\in {\rm d} y, \Theta_{\tau^{+}_{x}}\in {\rm d} \phi, \xi_{\tau^{+}_{x}}-x\in {\rm d} z\right)$$
converges weakly to a probability measure $\rho({\rm d} v, {\rm d} y, {\rm d}\phi, {\rm d} z)$ as $x\to +\infty$.

In particular,
$\bP_{0,\theta}\left(\xi_{\tau^{+}_{x}}-x\in {\rm d} z,\ \Theta_{\tau^{+}_{x}}\in {\rm d} \phi\right)$
converges weakly to a probability measure denoted by $\rho^{\ominus}({\rm d} z,{\rm d} \phi)$
and
$\bP_{0,\theta}\left(\xi_{\tau^{+}_{x}-}-x\in {\rm d} y,\ \Theta_{\tau^{+}_{x}-}\in {\rm d} v\right)$
converges weakly to a probability measure denoted by $\rho^{\oplus}({\rm d} y,{\rm d} v)$.
}

  \smallskip

  As alluded to above, we can use the former two main conclusions above to build a process which acts as an entrance law of the ssMp from the origin.

  \smallskip

{\it {\bf (C3): Candidate entrance law:}  Let $\P^\searrow_{z}$ denote the law of $X$ given by the Lamperti-Kiu transform \eqref{eq:lamperti_kiu} under $\hat\bP_{y,\theta}^\downarrow$ with $y=\log\|z\|$ and $\theta=\arg(z)$, and let $\varrho$ denote the image measure of $\rho^{\oplus}$ under the map $(y,\theta)\mapsto \theta\mathrm{e}^{y}$. Then the process $(X,\P^{\searrow}_{\varrho})$ has a finite lifetime $\bar{\zeta}$ with $X_{\bar{\zeta}-}=0$. Its time reversal process $((\widetilde{X}_{t}:=X_{(\bar{\zeta}-t)-})_{t<\bar{\zeta}},\P^{\searrow}_{\varrho})$ is a right continuous Markov process satisfying that $\widetilde{X}_{0}=0$ and $\widetilde{X}_{t}\not=0$ for all $t>0$. Moreover, $((\tilde{X}_{t})_{0<t<\bar{\zeta}},\P^{\searrow}_{\varrho})$ is a strong Markov process having the same transition rates as the ssMp $(X,\{\P_{z},z\in \mathcal{H}\})$ killed when exiting the unit ball.}

  \smallskip

  Moreover the stability of the overshoots and undershoots
  in the second main conclusion also helps
   with identifying the above candidate entrance law as unique in the sense of weak limits
  on the Skorokhod space.

\smallskip

{\it {\bf (C4): Uniqueness of the entrance law:} There exists a probability measure $\P_{0}$ such that
\begin{enumerate}
\item{$\mbox{w-}\lim_{z\to 0}\P_{z}=\P_{0}$
    in the weak sense of measures on the Skorokhod space.}

\item{$(X,\{\P_{z},z\in\h_{0}\})$ is a ssMp.}

\item{$(X,\{\P_{z},z\in \h_{0}\})$ is a Feller process.}\label{p:Feller}

\item{$((X_{t})_{t<\tau^{\ominus}_{r}},\P_0)$ is
equal in law with
$(r X_{(\bar{\zeta}-r^{-\alpha}t)-})_{t<r^{\alpha}\bar{\zeta}},\P^{\searrow}_{\varrho})$ for every $r>0$.}\label{p:distribution}

\item{Under $\P_{0}$ the process $X$ starts at $0$ and leaves $0$ instantaneously.}\label{p:instant}
\end{enumerate}
Here $\tau^{\ominus}_{r}=\inf\{t>0:\ \|X_{t}\|>r\}$.
 Moreover, $\P_{0}$ is the unique probability measure such that the extension $(X,\{\P_{z},z\in\h_{0}\})$ is a right continuous Markov process satisfying either \eqref{p:Feller} or \eqref{p:instant} listed above.}

  \smallskip

Finally we can reassert the stability of the underlying MAP over/undershoots  to generate the unique entrance law at the origin, but now in terms of the ssMp.

\smallskip

{\it {\bf (C5): Stability of the the process started at the origin:} For every $\delta>0$,
$((X_{\tau^{\ominus}_{\delta}-},X_{\tau^{\ominus}_{\delta}}),\P_{z})$ converges in distribution to $((X_{\tau^{\ominus}_{\delta}-},X_{\tau^{\ominus}_{\delta}}),\P_{0})$ as $z\to 0$,
and
\begin{align}
\mbox{w-}\lim_{\h\ni z\to 0}&\P_{z}\left(\arg(X_{\tau^{\ominus}_{1}-})\in {\rm d}v,\ \log\|X_{\tau^{\ominus}_{1}-}\|\in {\rm d}y,\ \arg(X_{\tau^{\ominus}_{1}})\in {\rm d}\phi,\ \log\|X_{\tau^{\ominus}_{1}}\|\in {\rm d}z\right)\nonumber\\
&=\P_{0}\left(\arg(X_{\tau^{\ominus}_{1}-})\in {\rm d}v,\ \log\|X_{\tau^{\ominus}_{1}-}\|\in {\rm d}y,\ \arg(X_{\tau^{\ominus}_{1}})\in {\rm d}\phi,\ \log\|X_{\tau^{\ominus}_{1}}\|\in {\rm d}z\right)\nonumber\\
&=\rho({\rm d}v,{\rm d}y,{\rm d}\phi,{\rm d}z).\nonumber
\end{align}}

In the case $d=1$ and the ssMp is positive, several works have established
   the limit $\P_{0}=\mbox{w-}\lim_{z\to 0}\P_{z}$
   using various techniques, see \cite{BC02,BS,BY02,CC,R07}. Recently, in the case when ssMp is allowed to take negative values as well, entrance laws were obtained in \cite{DDK}. Our contribution here is two-fold. Firstly we show, under suitable conditions, the existence of an entrance law at $0$ for a ssMp in any dimension. Secondly, our proof here uses a path reversal argument which follows the spirit of \cite{BS,DDK}, but works directly with the reversal of the ssMp rather than the underlying MAP. This appeals to the full strength  of Hunt-Nagasawa duality  as explored in e.g. \cite{Nagasawa, ChungWalsh}. We note that, in dimension $d=1$ or $d = 1/2$ (i.e. positive self-similar Markov processes), taking all fluctuation theory for granted in those settings (which means fluctuation theory of L\'evy processes for $d = 1/2$), our approach offers an alternative simple proof of the entrance laws.

  \smallskip

The rest of this paper is structured as follows. In Section \ref{sec:preliminary} we develop the fluctuation theory for general MAPs, which we believe is of independent interest and should be useful in studying ssMps. In Section \ref{sec:duality} we present the notions of duality as well as several time-reversal results about duality.
Among them, Lemma \ref{lem:walsh_killing_time} plays a key role in our path reversal argument.
 In Section \ref{sec:main results}, we present our working assumptions and the main result, Theorem \ref{thm:main}, which gives the existence of a weak limit of $\P_{z}$ as $z\to 0$, as well as the explicit law of the process started at the origin.
 The large number of assumptions given there, largely pertains to stability conditions that permit the aforesaid weak convergence.
  {In Section \ref{sec:application} we give two interesting examples to illustrate the main result. }Our main result is proved step by step through the arguments in Sections \ref{sec:condi-negative}-\ref{sec:convergence}: Firstly we define a family of probability measures $\{\hat{\bP}^{\downarrow}_{x,\theta},x\le 0,\theta\in\s\}$ under which the MAP $(\xi,\Theta)$ is conditioned to stay negative. Then we show both the overshoots and undershoots of the MAP $(\xi,\Theta)$ have stationary distributions, which we denote by $\rho^{\ominus}$ and $\rho^{\oplus}$ respectively. Starting from $((\xi,\Theta),\hat{\bP}^{\downarrow}_{\rho^{\oplus}})$ we construct by Lamperti-Kiu transform the process $(X,\P^{\searrow}_{\varrho})$ which is conditioned to stay inside the unit ball and hit the origin in a finite time. By time-reversing $(X,\P^{\searrow}_{\varrho})$ from its lifetime, we get the law of $(X,\P_{0})$ until first exit from a unit ball. Finally we prove $\P_{0}$ is the weak limit of $\P_{z}$ as $z\to 0$.

\smallskip

Notation:  Throughout this paper, we use
``$:=$'' as definition and ``$\stackrel{d}=$'' to mean ``equal in distribution''.
{Suppose $E$ is a locally compact separable metric space. Let $E_{\partial}=E\cup\{\partial\}$ (where $\partial\not\in E$) be the one-point compactification of $E$. Then $E_{\partial}$ is a compact separable metric space. For $T\in [0,+\infty]$, let $\mathbb D_E{[0,T)}$ denote the space of functions $\omega:{[0,T)}\rightarrow E_{\partial}$,}
such that there exists $\zeta=\zeta(\omega)\in [0,T]$, called the lifetime of $\omega$, with the property that  $t \mapsto \omega(t)$ is a c\`{a}dl\`{a}g function from ${[0,\zeta)}$ to $E$ and $\omega(t)=\partial$ for $t\geq \zeta$.
We endow the space $\mathbb D_E{[0,T)}$  with the Skorokhod topology which makes it into a Polish space.
We use the shorthand notation $\mathbb D_E=\mathbb D_E {[0,\infty)}$.
{Unless stipulated otherwise, every function $f$ on $E$ is automatically extended to $E_{\partial}$ by setting $f(\partial)=0$.}
For a point $x\in\R^{d}$, we use $\|x\|$ to denote its Euclidean norm.
For $q>0$, we use $\mathbf{e}_q$ to denote an independent exponential random variable with mean $1/q$.

\part{Fluctuation theory of Markov additive processes}

\section{Preliminaries}\label{sec:preliminary}

\subsection{Markov additive processes and L\'{e}vy systems}
Suppose $(\xi_{t},\Theta_{t})_{t\ge 0}$ is the coordinate process in $\mathbb{D}_{\R\times\s}$ and
$$\left((\xi,\Theta),\bP \right)=\left((\xi_{t},\Theta_{t})_{t\ge 0},\mathcal{F}_{\infty},(\mathcal{F}_{t})_{t\ge 0},\{\bP_{x,\theta}:(x,\theta)\in\R\times\s\}\right)$$ is a (possibly killed) Markov process with $\bP_{x,\theta}\left(\xi_{0}=x,\Theta_{0}
=\theta\right)=1$. Here $(\mathcal{F}_{t})_{t\ge 0}$ is the minimal augmented admissible filtration and $\mathcal{F}_{\infty}=\bigvee_{t=0}^{+\infty}\mathcal{F}_{t}$.

  \begin{definition}
    The process $\left((\xi,\Theta),\bP\right)$ is called a Markov additive process (MAP) on $\R\times\s$ if, for any $t \geq 0$, given $\{(\xi_s,\Theta_s), s\leq t\}$, the process $(\xi_{s+t}-\xi_t,\Theta_{s+t})_{s\geq 0}$ has the same law as $(\xi_{s},\Theta_{s})_{s\ge 0}$ under $\bP_{0,v}$ with $v=\Theta_t$. We call $\left((\xi,\Theta),\bP\right)$ a nondecreasing MAP if $\xi$ is a nondecreasing process on $\R$.
  \end{definition}

For a MAP process $\left((\xi,\Theta),\bP\right)$, we call $\xi$ the \textit{ordinate} and $\Theta$ the \textit{modulator}. By definition we can see that a MAP is translation invariant in $\xi$, i.e., $((\xi_{t},\Theta_{t})_{t\ge 0},\bP_{x,\theta})$ is equal in law to $((\xi_{t}+x,\Theta_{t})_{t\ge 0},\bP_{0,\theta})$ for all $x\in \R$ and $\theta\in\s $.

\smallskip

We assume throughout the paper that $(\Theta_{t})_{t\ge 0}$ is a Hunt process and $(\xi_{t})_{t\ge 0}$ is quasi-left continuous on $[0,\zeta)$. Then it is shown in \cite{Cinlar2} that there exist a continuous increasing additive functional $t\mapsto H_{t}$ of $\Theta$ and a transition kernel $\Pi$ from $\s $ to $\s \times \R$ satisfying
$$\Pi(\theta,\{(\theta,0)\})=0,\quad \quad \int_{\R}\left(1\wedge|y|^{2}\right)\Pi(\theta,\{\theta\}\times {\rm d}y)<+\infty\quad\forall \theta\in\s,$$
such that, for every nonnegative measurable function $f:\s \times \s \times \R\to \R^{+}$, every $\theta\in \s $ and $t\ge 0$,
\begin{align}
\bP_{0,\theta}&\left[\sum_{s\leq t}f(\Theta_{s-}, \Theta_{s}, \xi_{s}-\xi_{s-})\1_{\{\Theta_{s-}\neq \Theta_{s}\ \text{or}\ \xi_{s-}\neq \xi_{s}\}}\right]\nonumber\\
&=\bP_{0,\theta}\left[\int^{t}_{0}{\rm d}H_{s} \int_{\s \times \R}\Pi(\Theta_{s}, {\rm d}v, {\rm d}y)f(\Theta_{s}, v, y)\right].\nonumber
\end{align}
This pair $(H,\Pi)$ is said to be a \textit{L\'{e}vy system} for $\left((\xi,\Theta),\bP\right)$.
It can be shown that for every nonnegative predictable process $Z$ and nonnegative measurable function $g:\s \times \R \times \s \times \R\to \R^{+}$,
\begin{align}
\bP_{0,\theta}&\left[\sum_{s\le t}Z_{s}g(\Theta_{s-},\xi_{s-},\Theta_{s},\xi_{s})\1_{\{\Theta_{s-}\neq \Theta_{s}\ \text{or}\ \xi_{s-}\neq \xi_{s}\}}\right]\nonumber\\
&=\bP_{0,\theta}\left[\int^{t}_{0}{\rm d}H_{s}Z_{s}\int_{\s \times \R}\Pi(\Theta_{s},{\rm d}v,{\rm d}y)g(\Theta_{s},\xi_{s},v,\xi_{s}+y)\right]\label{levysys}
\end{align}
for all $\theta\in\s $ and $t\ge 0$.

\smallskip

The topic of MAPs are covered in various parts of the literature. We refer to \cite{AsmussenQueue,AA,CPR,Cinlar2,Cinlar1,KPal} to name but a few of the  texts and papers which give a general treatment.

\smallskip

For the remainder of the paper we will restrict ourselves to the setting that, up to killing of the MAP, $H_t = t$. On account of the bijection in \eqref{eq:time_change},  this naturally puts us in a restricted class of ssMps through the underlying driving MAP, however, as we will shortly see, it is on the MAP that we will impose additional assumptions.

 \subsection{Fluctuation theory for MAPs}\label{sec:fluctuation}

\begin{definition}
For any $y\in\R$, let $\tau^{+}_{y}:=\inf\{t>0:\xi_{t}>y\}$.
We say that $((\xi,\Theta),\bP)$ is
upwards regular
if
$$\bP_{0,\theta}\left(\tau^{+}_{0}=0\right)=1\quad\forall \theta\in\s.$$
  \end{definition}

 Suppose $(X,\P)=((X_{t})_{t\ge 0},\{\P_{z}:z\in\h\})$ is the ssMp associated to the MAP $((\xi,\Theta),\bP)$ via Lamperti-Kiu transform. We say $(X,\P)$ is \textit{sphere-exterior regular} if $((\xi,\Theta),\bP)$ is
 upwards regular.
 For $r>0$, let $\tau^{\ominus}_{r}:=\inf\{t>0:\ \|X_{t}\|>r\}$.
Immediately by the definition, $(X,\P)$ is sphere-exterior regular if and only if $\mathbb{P}_{z}\left(\tau^{\ominus}_{1}=0\right)=1$ for all $z\in \h$ with $\|z\|=1$.

\smallskip

In the remaining of this paper we assume that the MAP $((\xi,\Theta),\bP)$ is
upwards regular.
This assumption is not really necessary but nevertheless avoids a lot of unnecessary technicalities when we explore the fluctuation properties.

\subsubsection{Excursion from maximum/minimum.}\label{excmm}
Let $\bar{\xi}_{t}:=\sup_{s\le t}\xi_{s}$ and $U_{t}:=\bar{\xi}_{t}-\xi_{t}$. Then under $\bP_{0,\theta}$ the process $(\Theta_{t},\xi_{t},U_{t})_{t\ge 0}$ is an $\s \times \R\times \R^{+}$-valued right process started at $(\theta,0,0)$, whose transition semigroup on $(0,+\infty)$ is given by
$$P_{t}f(v,x,u):=\bP_{0,v}\left[f\left(\Theta_{t},\xi_{t}+x,u\vee \bar{\xi}_{t}-\xi_{t}\right)\right]$$
for every $t>0$ and every nonnegative measurable function $f:\s \times\R \times\R^{+}\to \R^{+}$.
We shall work with the canonical realization of $(\Theta_{t},\xi_{t},U_{t})_{t\ge 0}$ on the sample space
{$\mathbb{D}_{\s\times \R\times \R^{+}}$.}

\smallskip

We define
$\bar{M}:=\{t\ge 0 : U_{t}=0\}$ and $\bar{M}^{cl}$ its closure in $\R^{+}$.
Obviously the set $\R^{+}\setminus \bar{M}^{cl}$ is an open set and can be written as a union of intervals. We use $\bar{G}$ and $\bar{D}$, respectively, to denote the sets of left and right end points of such intervals.
Define $\bar{R}:=\inf\{t>0:\ t\in \bar{M}^{cl}\}$.
The upwards regularity implies that every point in $\s$ is regular for $\bar{M}$
{in the sense that $\bP_{0,\theta}\left(\bar{R}=0\right)=1$ for all $\theta\in\s$.}
Thus by
{\cite[Theorem (4.1)]{M}}
there exist a continuous additive functional $t\mapsto\bar{L}_{t}$ of
{$(\Theta_{t},\xi_{t},U_{t})_{t\ge 0}$}
which is carried by $\s \times \R\times\{0\}$ and a kernel $\mathfrak{P}$ from
{$\s \times \R\times \R^{+}$ into $\mathbb{D}_{\s\times\R\times \R^{+}}$}
satisfying
$\mathfrak{P}^{\theta,x,u}\left(\bar{R}=0\right)=0$ and $\mathfrak{P}^{\theta,x,u}\left(1-\mathrm{e}^{-\bar{R}}\right)\le 1$ such that
\begin{equation}\label{M1}
\bP_{0,\theta}\left[\sum_{s\in \bar{G}}Z_{s}\,f\circ\theta_{s}\right]=\bP_{0,\theta}\left[\int_{0}^{+\infty}Z_{s}\,\mathfrak{P}^{\Theta_{s},\xi_{s},0}(f){\rm d}\bar{L}_{s}\right]
\end{equation}
for any
{nonnegative predictable}
process $Z$ and
{any nonnegative function $f$ which is measurable with respect to $\sigma\left((\Theta_{t},\xi_{t},U_{t})_{t\ge 0}\right)$.}
{Moreover, by \cite[Theorem (5.1)]{M}, $\mathfrak{P}^{\theta,x,0}\left((\Theta_{0},\xi_{0},U_{0})\not=(\theta,x,0)\right)=0$, and under $\mathfrak{P}^{\theta,x,0}$, the process $(\Theta_{t},\xi_{t},U_{t})_{t>0}$ has the strong Markov property (as defined in \cite[(5.2)]{M}) with respect to $P_{t}$.}
In particular, if $f$ is measurable with respect to $\sigma((\Theta_{t},U_{t})_{t\ge 0})$, then the right-hand side of \eqref{M1} equals
$$\bP_{0,\theta}\left[\int_{0}^{+\infty}Z_{s}\,\mathfrak{P}^{\Theta_{s},0}(f){\rm d}\bar{L}_{s}\right],$$
{where $\mathfrak{P}$ denotes the kernel from $\s\times\R^{+}$ into $\mathbb{D}_{\s\times \R^{+}}$ for the process $(\Theta_{t},U_{t})_{t\ge 0}$ defined in the same way of \cite{M}.}
It is known (see, for example, \cite[Section 3]{Kaspi83}) that there is a nonnegative measurable function $\ell^+:\s \to \R^{+}$ such that
\begin{equation}\label{M2}
\int_{0}^{t}\1_{\{s\in \bar{M}\}}{\rm d}s=\int_{0}^{t}\1_{\{s\in \bar{M}^{cl}\}}{\rm d}s=\int_{0}^{t}\ell^+(\Theta_{s}){\rm d}\bar{L}_{s}\quad\forall t\ge 0\quad \bP_{0,\theta}\mbox{-a.s.}
\end{equation}
Let $\bar{L}^{-1}_{t}$ be the right inverse process of $\bar{L}_{t}$. Define $\xi^{+}_{t}:=\xi_{\bar{L}^{-1}_{t}}$ and $\Theta^{+}_{t}:=\Theta_{\bar{L}^{-1}_{t}}$ for all $t$ such that $\bar{L}^{-1}_{t}<+\infty$ and otherwise $\xi^{+}_{t}$ and $\Theta^{+}_{t}$ are both assigned to be the cemetery state $\partial$.
{One can verify by the strong Markov property that $(\bar{L}^{-1}_{t},\xi^{+}_{t},\Theta^{+}_{t})_{t\ge 0}$ defines a MAP, whose first two elements are ordinates. Similarly, both $(\xi^{+}_{t},\Theta^{+}_{t})_{t\ge 0} $ and $(\bar{L}^{-1}_{t},\Theta^{+}_{t})_{t\ge 0} $ are MAPs.}
These three processes are referred to as \textit{ascending ladder process}, \textit{ascending ladder height process} and \textit{ascending ladder time process}, respectively.

\smallskip

Suppose the set $\mathbb{R}^{+}\setminus \bar{M}^{cl}$ is written as a union of random intervals $(g,d)$. For such intervals, define
\begin{equation} \nonumber
(\epsilon^{(g)}_{s},\nu^{(g)}_{s}):=\begin{cases}
        (U_{g+s},\Theta_{g+s})
         \quad &\hbox{if } 0\le s<d-g, \smallskip \\
        (U_{d},\Theta_{d})
         \quad &\hbox{if } s\ge d-g.
        \end{cases}
\end{equation}
$(\epsilon^{(g)}_{s},\nu^{(g)}_{s})_{s\ge 0}$ is called an excursion from the maximum and $\zeta^{(g)}:=d-g$ is called its lifetime.
We use $\mathcal{E}$ to denote the collection
$\{(\epsilon^{(g)}_{s}(\omega),\nu^{(g)}_{s}(\omega))_{s\ge 0}:g\in \bar{G}(\omega),\ \omega\in
{\mathbb{D}_{\s\times\R\times \R^{+}}}
\}$,
and call it the space of excursions. Let ${\rm n}^{+}_{\theta}$ be the image measure of $\mathfrak{P}^{\theta,0}$ under the mapping that stops the path of $(\Theta_{t},U_{t})_{t\ge 0}$ at time $\bar{R}$.
A direct consequence of~\cite[equation (4.9)]{M} is that
for any
{nonnegative}
measurable functionals
{$F:\mathbb{D}_{\R^{+}\times \s }\to \R^{+}$}
and
{$G:\mathbb{R}^{+}\times \mathbb{D}_{\mathbb{R}\times \s }\to\R^{+}$,}
\begin{align}\label{eq:last exit}
   \bP_{y,\theta}&\left[ \sum_{g\in \bar{G}} G\big(g,(\xi_t,\Theta_t)_{t \leq g}\big) F\big(\epsilon^{(g)},{\nu}^{(g)}\big) \right]\nonumber\\
   &=\bP_{y,\theta} \left[ \int_0^\infty {\rm d}\bar{L}_s\, G\big(s,(\xi_t,\Theta_t)_{t\leq s}\big) \int_{\mathcal{E}} {\rm n}^{+}_{\Theta_s}({\rm d}\epsilon,{\rm d}{\nu}) F(\epsilon,{\nu})  \right].
  \end{align}
We call $\{{\rm n}^{+}_{\theta}:\theta\in \s\}$ the \textit{excursion measures at the maximum}.

\smallskip

The excursion measures at the minimum and
descending ladder process are defined analogously replacing $\xi$ by $-\xi$.

\subsubsection{Fluctuation identities.}

For $t>0$, define
$$\bar{g}_{t}:=\sup\{s\le t: s\in \bar{M}^{cl}\}\quad \mbox{ and }\quad\bar{\Theta}_{t}:=\Theta_{\bar{g}_{t}}\1_{\{\bar{\xi}_{t}=\xi_{\bar{g}_{t}}\}}+\Theta_{\bar{g}_{t}-}\1_{\{\bar{\xi}_{t}
>\xi_{\bar{g}_{t}}\}}.$$
By the right continuity of sample paths one can easily show that $\bar{g}_{t}$ is equal to $\sup\{s\le t:\ s\in\bar{M}\}$ with probability $1$. Since by quasi-left continuity, $\bP_{0,\theta}(\xi_{t}\not=\xi_{t-})=0$ for all $t>0$, we have $\bP_{0,\theta}\left(\bar{g}_{t}=\sup\{s<t:s\in \bar{M}\}\right)=1$ for all $t>0$.
We claim that $\bP_{0,\theta}$-almost surely $\bar{g}_{t}$ is not a jump time of the process $(\xi,\Theta)$ for every $\theta\in\s$. Otherwise, one can construct a stopping time $T$ such that
{$$\bP_{0,\theta}\left(\{T\in \bar{G}\}\cap \{\xi_{T-}\not=\xi_{T}\mbox{ or } \Theta_{T-}\not=\Theta_{T}\}\right)>0.$$}
Noting that by \eqref{M1}
 $$\bP_{0,\theta}\left(\sum_{g\in \bar{G}}\1_{\{U_{g}>0\mbox{ or }\Theta_{g-}\not=\Theta_{g}\}}\right)=\bP_{0,\theta}\left[\int_{0}^{+\infty}\mathfrak{P}^{\Theta_{s},\xi_{s},0}\left(U_{0}>0\mbox{ or }\Theta_{0}\not=0\right)\mathrm{d}\bar{L}_{s}\right]=0,$$
we get from the above inequality that
{$\bP_{0,\theta}\left(T\in \bar{G},\ \xi_{T-}<\xi_{T}\right)>0$.}
This brings a contradiction, since if we apply Markov property and upwards regularity at $T$, we get $\xi_{T+s}>\xi_{T}>\xi_{T-}$ for $s$ sufficiently small on the event
{$\{T\in \bar{G},\ \xi_{T-}<\xi_{T}\}$,} which is impossible.

\smallskip

The following identity is one of the key tools in extending identities from the fluctuation theory for L\'evy processes to MAPs, it is the base to establish a Wiener-Hopf type factorisation for MAPs.

  \begin{proposition}\label{prop:wiener_hopf}
    Suppose that $((\xi,\Theta),\bP)$ is a Markov additive process taking values in $\R\times \s $.
   Then for every bounded measurable functions $F,G:{[0,\infty)}\times \R\times \s\to \R$ and every $\theta\in\s $,
    \begin{align}
      \bP_{0,\theta}&\left[G\big(\bar g_{\mathbf{e}_q},\bar \xi_{\mathbf{e}_q},\bar{\Theta}_{\mathbf{e}_{q}}\big)  F\big(\mathbf{e}_q-\bar g_{\mathbf{e}_q},\bar \xi_{\mathbf{e}_q}-\xi_{\mathbf{e}_q},\Theta_{\mathbf{e}_q}\big) \right]\nonumber\\
      &=\int_{\R^{+}\times \s\times\R^{+} }\mathrm{e}^{-q r}G(r,z,v)\left[q \ell^+(v) F(0,0,v)+\emph{\rm n}^{+}_{v}\left(F(\mathbf{e}_{q},\epsilon_{\mathbf{e}_{q}},\nu_{\mathbf{e}_{q}})\1_{\{\mathbf{e}_{q}<\zeta\}}\right)\right]V^{+}_{\theta}({\rm d}r,{\rm d}v,{\rm d}z),\nonumber
    \end{align}
    where
    $$V^{+}_{\theta}({\rm d}r,{\rm d}v,{\rm d}z):=\bP_{0,\theta}\left[\int_{0}^{\bar{L}_{\infty}}\1_{\{\bar{L}^{-1}_{s}\in {\rm d}r,\ \Theta^{+}_{s}\in {\rm d}v,\ \xi^{+}_{s}\in {\rm d}z\}}{\rm d}s\right].$$
  \end{proposition}
  \proof
     It is known from the above argument that $\bP_{0,\theta}$-almost surely $\bar{g}_{\mathbf{e}_{q}}$ is not jump time of $(\xi,\Theta)$, and thus $(\bar{\xi}_{\mathbf{e}_{q}},\bar{\Theta}_{\mathbf{e}_{q}})=(\xi_{\bar{g}_{\mathbf{e}_{q}}},\Theta_{\bar{g}_{\mathbf{e}_{q}}})$ $\bP_{0,\theta}$-a.s.
    Then we have
    \begin{align}\label{eq:wh}
      \bP_{0,\theta}&\left[ F\big(\mathbf{e}_q-\bar g_{\mathbf{e}_q},\bar \xi_{\mathbf{e}_q}-\xi_{\mathbf{e}_q},\Theta_{\mathbf{e}_q}\big) G\big(\bar g_{\mathbf{e}_q},\bar \xi_{\mathbf{e}_q},\bar{\Theta}_{\mathbf{e}_{q}}\big) \right]\nonumber\\
      &=\bP_{0,\theta}\left[ F\big(\mathbf{e}_q-\bar g_{\mathbf{e}_q},\bar \xi_{\mathbf{e}_q}-\xi_{\mathbf{e}_q},\Theta_{\mathbf{e}_q}\big) G\big(\bar g_{\mathbf{e}_q},\bar \xi_{\mathbf{e}_q},\bar{\Theta}_{\mathbf{e}_{q}}\big) \1_{\{\xi_{\mathbf{e}_q}=\bar\xi_{\mathbf{e}_q}\}}\right]\nonumber\\
      &\qquad+\bP_{0,\theta}\left[ F\big(\mathbf{e}_q-\bar g_{\mathbf{e}_q},\bar \xi_{\mathbf{e}_q}-\xi_{\mathbf{e}_q},\Theta_{\mathbf{e}_q}\big) G\big(\bar g_{\mathbf{e}_q},\bar \xi_{\mathbf{e}_q},\bar{\Theta}_{\mathbf{e}_{q}}\big) \1_{\{\xi_{\mathbf{e}_q}<\bar\xi_{\mathbf{e}_q}\}}\right]\nonumber\\
      &=\bP_{0,\theta}\left[ F\big(0,0,\Theta_{\mathbf{e}_q}\big) G\big(\mathbf{e}_q,\xi_{\mathbf{e}_q},\Theta_{\mathbf{e}_q}\big) \1_{\{\xi_{\mathbf{e}_q}=\bar\xi_{\mathbf{e}_q}\}}\right]\nonumber\\
      &\qquad+\bP_{0,\theta}\left[\sum_{g\in \bar{G}}\1_{\{g<\mathbf{e}_q<g+\zeta^{(g)}\}} F\big(\mathbf{e}_q-g,\epsilon^{(g)}_{\mathbf{e}_q-g},{\nu}^{(g)}_{\mathbf{e}_q-g}\big) G\big(g,\bar \xi_{g},\Theta_{g}\big) \right].\nonumber
      \end{align}
     By \eqref{eq:last exit}, the above sum is equal to
      \begin{align}
      &\bP_{0,\theta}\left[ F\big(0,0,\Theta_{\mathbf{e}_q}\big) G\big(\mathbf{e}_q,\xi_{\mathbf{e}_q},\Theta_{\mathbf{e}_q}\big) \1_{\{\xi_{\mathbf{e}_q}=\bar\xi_{\mathbf{e}_q}\}}\right]\nonumber\\
      &\qquad+\bP_{0,\theta}\left[\int_0^\infty {\rm d}\bar L_s\1_{\{s<\mathbf{e}_q\}} G\big(s,\bar \xi_{s},\Theta_{s}\big) {\rm n}^{+}_{\Theta_s}\left(F\big(\mathbf{e}_q-s,\epsilon_{\mathbf{e}_q-s},{\nu}_{\mathbf{e}_q-s}\big)\1_{\{\mathbf{e}_{q}-s<\zeta\}}\right)\right],
    \end{align}
 For the second term we can use the memorylessness of the exponential distribution and a change of variable to yield
    \begin{align}
    \bP_{0,\theta}&\left[\int_0^\infty {\rm d}\bar L_s\1_{\{s<\mathbf{e}_q\}} G\big(s,\bar \xi_{s},\Theta_{s}\big) {\rm n}^{+}_{\Theta_s}\left(F\big(\mathbf{e}_q-s,\epsilon_{\mathbf{e}_q-s},{\nu}_{\mathbf{e}_q-s}\big)\1_{\{\mathbf{e}_{q}-s<\zeta\}}\right) \right] \nonumber\\
    &=\bP_{0,\theta}\left[\int_0^\infty {\rm d}\bar L_s\mathrm{e}^{-q s}G\big(s,\bar \xi_{s},\Theta_{s}\big) {\rm n}^{+}_{\Theta_s}\left(F\big(\mathbf{e}_q,\epsilon_{\mathbf{e}_q},{\nu}_{\mathbf{e}_q}\big)\1_{\{\mathbf{e}_{q}<\zeta\}}\right) \right]\nonumber\\
    &=\bP_{0,\theta}\left[\int_0^{\bar{L}_{\infty}} {\rm d}s\,\mathrm{e}^{-q\bar{L}^{-1}_{s}} G\big(\bar{L}^{-1}_{s},\xi^{+}_{s},\Theta^{+}_{s}\big) {\rm n}^{+}_{\Theta^{+}_s}\left(F\big(\mathbf{e}_q,\epsilon_{\mathbf{e}_q},{\nu}_{\mathbf{e}_q}\big)\1_{\{\mathbf{e}_{q}<\zeta\}}\right) \right].\label{prop3.1.1}
    \end{align}
For the first term in \eqref{eq:wh} we use (\ref{M2}) to get
\begin{align}
    &\bP_{0,\theta}\left[ F\big(0,0,\Theta_{\mathbf{e}_q}\big) G\big(\mathbf{e}_q,\xi_{\mathbf{e}_q},\Theta_{\mathbf{e}_q}\big) \1_{\{\xi_{\mathbf{e}_q}=\bar\xi_{\mathbf{e}_q}\}}\right]\nonumber\\
    &=q\bP_{0,\theta}\left[\int_{0}^{+\infty}\mathrm{e}^{-q t}G(t,\bar{\xi}_{t},\Theta_{t})F(0,0,\Theta_{t})\1_{\{t\in\bar{M}\}}{\rm d}t\right]\nonumber\\
    &=q\bP_{0,\theta}\left[\int_{0}^{+\infty}\mathrm{e}^{-q t}G(t,\bar{\xi}_{t},\Theta_{t})F(0,0,\Theta_{t})\ell^+(\Theta_{t}){\rm d}\bar{L}_{t}\right]\nonumber\\
    &=q\bP_{0,\theta}\left[\int_{0}^{\bar{L}_{\infty}}\mathrm{e}^{-q \bar{L}^{-1}_{s}}G(\bar{L}^{-1}_{s},\xi^{+}_{s},\Theta^{+}_{s})F(0,0,\Theta^{+}_{s})\ell^+(\Theta^{+}_{s}){\rm d}s\right].\label{prop3.1.2}
    \end{align}
By plugging \eqref{prop3.1.1} and \eqref{prop3.1.2} into \eqref{eq:wh} we get that
\begin{align}
&\bP_{0,\theta}\left[ F\big(\mathbf{e}_q-\bar g_{\mathbf{e}_q},\bar \xi_{\mathbf{e}_q}-\xi_{\mathbf{e}_q},\Theta_{\mathbf{e}_q}\big) G\big(\bar g_{\mathbf{e}_q},\bar \xi_{\mathbf{e}_q},\bar{\Theta}_{\mathbf{e}_{q}}\big) \right]\nonumber\\
&=\bP_{0,\theta}\big[\int_{0}^{\bar{L}_{\infty}}{\rm d}s\,\rme^{-q \bar{L}^{-1}_{s}}G\left(\bar{L}^{-1}_{s},\xi^{+}_{s},\Theta^{+}_{s}\right)\big(q \ell^+(\Theta^{+}_{s})F(0,0,\Theta^{+}_{s})
+{\rm n}^{+}_{\Theta^{+}_{s}}\left(F(\mathbf{e}_{q},\epsilon_{\mathbf{e}_{q}},\nu_{\mathbf{e}_{q}})\1_{\{\mathbf{e}_{q}<\zeta\}}\right)\big)\big].\nonumber
\end{align}
We have thus proved this proposition.\qed

\smallskip

\begin{corollary}\label{cor:dist WH}
For every $\theta\in\s $, we have
\begin{align}
\bP_{0,\theta}&\left(\bar{\xi}_{\mathbf{e}_{q}}\in {\rm d}z,\ \bar{\xi}_{\mathbf{e}_{q}}-\xi_{\mathbf{e}_{q}}\in {\rm d}w,\ \Theta_{\mathbf{e}_{q}}\in {\rm d}v\right)\nonumber\\
&=\delta_{0}({\rm d}w)\ell^+(v)\int_{0}^{+\infty}q\mathrm{e}^{-q r}V^{+}_{\theta}({\rm d}r,{\rm d}v,{\rm d}z)\nonumber\\
&+\int_{(r,u)\in \R^{+}\times \s }\mathrm{e}^{-q r}\emph{\rm n}^{+}_{u}\left(\epsilon_{\mathbf{e}_{q}}\in {\rm d}w,\ \nu_{\mathbf{e}_{q}}\in {\rm d}v,\ \mathbf{e}_{q}<\zeta\right)V^{+}_{\theta}({\rm d}r,{\rm d}u,{\rm d}z).\nonumber
\end{align}
\end{corollary}

\smallskip

The excursion measures allow us to gain some additional insight into the analytical form of the jumping measures of
the ascending ladder processes.

\begin{proposition}\label{prop:Pi+}
Suppose $((\xi,\Theta),\bP)$ is a MAP with L\'{e}vy system $(H,\Pi)$ where $H_{t}=t\wedge \zeta$. Then the ascending ladder process $(( \bar{L}^{-1},\xi^{+},\Theta^{+}),\bP)$ has a L\'{e}vy system
$(H^{+},\Gamma^{+})$
where $H^{+}_{t}=t\wedge \zeta^{+}$ and
$$
\Gamma^{+}(\theta,{\rm d}v, {\rm d}r, {\rm d}y)
=\delta_{0}({\rm d}r)\ell^+(\theta)\Pi(\theta,{\rm d}v,{\rm d}y)+{\rm n}^{+}_{\theta}\left(\Pi\left(\nu_{r},{\rm d}v,\epsilon_{r}+{\rm d}y\right),  r<\zeta\right){\rm d}r$$
for $\theta,v\in\s $, $r\ge 0$ and $y>0$.
Here $\zeta^{+}$ denotes the lifetime of $(\xi^{+},\Theta^{+})$.
{In particular, the ascending ladder height process $((\xi^{+},\Theta^{+}),\bP)$ has a L\'{e}vy system $(H^{+},\Pi^{+})$ where
$\Pi^{+}(\theta,{\rm d}v,{\rm d}y)=\Gamma^{+}(\theta,{\rm d}v,[0,+\infty),{\rm d}y)$ for $\theta,v\in\s $ and $y>0$.}
\end{proposition}

\proof To prove this proposition we apply the theory for L\'{e}vy systems and time-changed processes developed in \cite{F-G}.
We consider the strong Markov process $Y_{t}:=(\Theta_{t},\xi_{t},U_{t}, t)$ on the state space
{$\s \times \R\times \R^{+}\times \R^{+}$}
where $U_{t}=\bar{\xi}_{t}-\xi_{t}$.
Let $\bar{M}=\{t\ge 0:\ U_{t}=0\}$ and $\bar{R}=\inf\{t>0:\ t\in \bar{M}^{cl}\}$. It is known that the local time at the maximum $\bar{L}_{t}$ is a continuous additive functional carried by $\tilde{F}:=\s \times \R\times\{0\}\times\R^{+}$.
 The argument in the beginning of this subsection implies that almost surely the ``irregular part''
 (in the sense of \cite{F-G}) $G^{i}:=\{s\in \bar{G}:\ U_{s}\not=0\}$ is an empty set. Let $\tilde{Y}_{t}:=(\Theta^{+}_{t},\xi^{+}_{t},U^{+}_{t}, \bar{L}^{-1}_{t})$ be the time-changed process of $Y_{t}$ by the inverse local time $\bar{L}^{-1}_{t}$. It is a right process on the state space $\tilde{F}$. Then following the arguments and calculations in \cite[Section 5]{F-G}, one can get a L\'{e}vy system for this time-changed process. In fact, applying \cite[Theorem 5.2]{F-G} here, we have
 \begin{align}
\bP_{0,\theta}&\left[\sum_{s>0}F\left(\Theta^{+}_{s-},\xi^{+}_{s-}, \bar{L}^{-1}_{s-}, \Theta^{+}_{s},\xi^{+}_{s}, \bar{L}^{-1}_{s}\right)\1_{\{\xi^{+}_{s-}\not=\xi^{+}_{s}\}}\right]\nonumber\\
&=\bP_{0,\theta}\big[\int_{0}^{+\infty}{\rm d}s\int_{
{\s\times\R^{+}\times [s,+\infty)}
}F(\Theta^{+}_{s},\xi^{+}_{s}, s, v, y, u)\1_{\{\xi^{+}_{s}\not=y\}}
\big(\mathfrak{P}^{\Theta^{+}_{s},\xi^{+}_{s},0, s}\left(\Theta_{\bar{R}}\in {\rm d}v,\ \xi_{\bar{R}}\in {\rm d}y, \bar{R}\in {\rm d}u\right)\nonumber\\
&\hspace{3cm}+\ell^+(\Theta^{+}_{s})\Pi(\Theta^{+}_{s},{\rm d}v,{\rm d}y-\xi^{+}_{s})\delta_s({\rm d}u)\big)\big]\nonumber\\
&=\bP_{0,\theta}\big[\int_{0}^{+\infty}{\rm d}s\int_{\s\times\R^{+}\times\R^{+}}F(\Theta^{+}_{s},\xi^{+}_{s}, s, v, y, s+r)\1_{\{\xi^{+}_{s}\not=y\}}
\big(\mathfrak{P}^{\Theta^{+}_{s},\xi^{+}_{s},0, s}\left(\Theta_{\bar{R}}\in {\rm d}v,\ \xi_{\bar{R}}\in {\rm d}y, \bar{R}-s\in {\rm d}r\right)\nonumber\\
&\hspace{3cm}+\ell^+(\Theta^{+}_{s})\Pi(\Theta^{+}_{s},{\rm d}v,{\rm d}y-\xi^{+}_{s})\delta_0({\rm d}r)\big)\big]\label{prop3.4.0}
\end{align}
 for every nonnegative measurable function $F$.
Here $\mathfrak{P}^{\theta,x, 0, s}$ denotes the kernel $\mathfrak{P}^{\theta,x, 0}$ trivially extended to include the pure drift process issued from $s.$
{So, note that under $\mathfrak{P}^{\theta,x, 0, s}$ the process $(Y_{t})_{t>0}$ has the strong Markov property with respect to the same transition semigroup as $(Y,\bP_{x,\theta})$.}
Using this and the translation invariance, we have
\begin{equation}\label{prop3.4.1}
\mathfrak{P}^{\theta,x,0, s}\left[\1_{\{s+r<\bar{R},\xi_{\bar{R}}\not=x\}}f(\Theta_{\bar{R}},\xi_{\bar{R}},
{\bar{R}-s}
)\right]
=\mathfrak{P}^{\theta,0,0}\left[\1_{\{r<\bar{R}\}}\bP_{\xi_{r},\Theta_{r}}\left(f(\Theta_{\tau^{+}_{0}},\xi_{\tau^{+}_{0}}+x, r+\tau^{+}_{0})\1_{\{\xi_{\tau^{+}_{0}}>0\}}\right)\right]
\end{equation}
for any $r>0$ and nonnegative measurable function $f$. Since $(\xi,\Theta)$ has L\'{e}vy system $(H,\Pi)$ with $H_{t}=t\wedge \zeta$, we have
\begin{equation}\nonumber
\bP_{z,v}\left[f(\Theta_{\tau^{+}_{0}},\xi_{\tau^{+}_{0}}+x, r+\tau^{+}_{0})\1_{\{\xi_{\tau^{+}_{0}}>0\}}\right]
=\bP_{z,v}\left[\int_{0}^{\tau^{+}_{0}}{\rm d}t\int_{\s\times(-\xi_{t},+\infty)}f(w,\xi_{t}+x+y, r+t)\Pi(\Theta_{t},{\rm d}w,{\rm d}y)\right],
\end{equation}
where we used that, on the event $\{\xi_{\tau^{+}_{0}}>0\},$ $\tau^{+}_{0}$ is the first jump time of $\xi,$ that takes $\xi$ into the positive axis, and we apply (\ref{levysys}).
Plugging this in \eqref{prop3.4.1}, and using the Markov property under $\mathfrak{P}^{\theta,0,0}$, we have
\begin{align}
\mathfrak{P}^{\theta,x,0,s}&\left[\1_{\{s+r<\bar{R},\xi_{\bar{R}}\not=x\}}f(\Theta_{\bar{R}},\xi_{\bar{R}},
{\bar{R}-s}
\right]\nonumber\\
&=\mathfrak{P}^{\theta,0,0}\left[\1_{\{r<\bar{R}\}}\bP_{-U_{r},\Theta_{r}}\left(\int_{0}^{\tau^{+}_{0}}{\rm d}t\int_{\s\times(U_{t},+\infty)}f(w,-U_{t}+x+y, r+t)\Pi(\Theta_{t},{\rm d}w,{\rm d}y)\right)\right]\nonumber\\
&={\rm n}^{+}_{\theta}\left[\1_{\{r<\zeta\}}\int_{r}^{\zeta}{\rm d}t\int_{\s\times(\epsilon_{t},+\infty)}f(w,-\epsilon_{t}+x+y, t)\Pi(\nu_{t},{\rm d}w,{\rm d}y)\right].\nonumber
\end{align}
By letting $r\to 0+$, we get from above equation that
$$\1_{\{x\not=y\}}\mathfrak{P}^{\theta,x,0,s}\left(\Theta_{\bar{R}}\in {\rm d}v,\ \xi_{\bar{R}}\in {\rm d}y,
{\bar{R}-s}
\in {\rm d}t\right)={\rm d}t\,{\rm n}^{+}_{\theta}
\left[\int_{
{z\in}(\epsilon_{t},+\infty)}\1_{\{-\epsilon_{t}+z+x\in {\rm d}y\}}\Pi(\nu_{t},{\rm d}v,{\rm d}z)\right].$$
Plugging this in \eqref{prop3.4.0} yields that
\begin{align}
&\bP_{0,\theta}\big[\sum_{s> 0}F\left(\Theta^{+}_{s-},\xi^{+}_{s-},  \bar{L}^{-1}_{s-} ,\Theta^{+}_{s},\xi^{+}_{s},  \bar{L}^{-1}_{s}\right)\1_{\{\xi^{+}_{s-}\not=\xi^{+}_{s}\}}\big]\nonumber\\
&=\bP_{0,\theta}\big[\int_{0}^{+\infty}{\rm d}s\int_{
{\s\times\R^{+}\times\R^{+}}
} F(\Theta^{+}_{s},\xi^{+}_{s}, s, v,\xi^{+}_{s}+y, s+r)\big(\delta_{0}({\rm d}r)\ell^+(\Theta^{+}_{s})\Pi(\Theta^{+}_{s},{\rm d}v,{\rm d}y) \nonumber\\
&\quad+ {\rm n}^{+}_{\Theta^{+}_{s}}\left(\Pi\left(\nu_{r},{\rm d}v,\epsilon_{r}+{\rm d}y\right), r<\zeta\right) {\rm d}r\big)\big],\nonumber
\end{align}
which in turn yields the assertion of this proposition.\qed

\smallskip

\begin{remark}
\rm
Suppose $\xi$ is a non-killed $\R$-valued L\'{e}vy process with triplet $(a,\sigma^{2},\Pi)$ for which $0$ is regular for $(0,+\infty)$. This process can be viewed as the projection of a upwards regular MAP $(\xi,\Theta)$ where the modulator $\Theta$ is equal to a constant. Therefore all the above results we obtained for MAP can be applied to this L\'{e}vy process.
We use $\bP_{0}$ (resp. $\hat{\bP}_{0}$) to denote the law of $\xi$ (resp. $-\xi$) started from $0$.
It is a known fact that its ascending ladder process $(\bar{L}^{-1}_{t},\xi^{+}_{t})_{t\ge 0}$ is a (possibly killed) bivariate subordinator.
Let $\Pi^{+}$ be the L\'{e}vy measure of $\xi^{+}$. Proposition \ref{prop:Pi+} yields that for $y>0$,
\begin{equation}
\Pi^{+}(y,+\infty)=\ell^+\Pi(y,+\infty)+{\rm n}^{+}\left[\int_{0}^{+\infty}\Pi(\epsilon_{r}+y,+\infty){\rm d}r\right],\label{rm3.5.1}
\end{equation}
where $\ell^+$ is the drift coefficient of $\bar{L}^{-1}_{t}$ and ${\rm n}^{+}$ is the excursion measure at maximum.
It follows by Proposition \ref{prop:wiener_hopf} that for any nonnegative measurable function $F:\R\to \R^{+}$
\begin{eqnarray}
\bP_{0}\left[F(\bar{\xi}_{\mathbf{e}_{q}}-\xi_{\mathbf{e}_{q}})\right]&=&\frac{q\ell^+F(0)+{\rm n}^{+}\left[\int_{0}^{\zeta}q \mathrm{e}^{-q s}F(\epsilon_{s}){\rm d}s\right]}{\Phi(q)},\nonumber\\
\hat{\bP}_{0}\left[F(\bar{\xi}_{\mathbf{e}_{q}})\right]&=&\hat{\Phi}(q)\int_{\R^{+}\times\R^{+}}\mathrm{e}^{-q r}F(z)\hat{V}^{+}({\rm d}r,{\rm d}z)\label{rm3.5.2}
\end{eqnarray}
where $\hat{V}^{+}({\rm d}r,{\rm d}z):=\hat{\bP}_{0}\left[\int_{0}^{+\infty}\1_{\{\bar{L}^{-1}_{s}\in {\rm d}r,\xi^{+}_{s}\in {\rm d}z\}}{\rm d}s\right]$, and $\Phi(q)$ (resp. $\hat{\Phi}(q)$) is equal to the Laplace exponent of the (possibly killed) subordinator $(\bar{L}^{-1}_{t})_{t\ge 0}$ under $\bP_{0}$ (resp. $\hat{\bP}_{0}$).
The Wiener-Hopf factorization of L\'{e}vy process implies that
$\Phi(q)\hat{\Phi}(q)=\kappa q$ for some constant $\kappa>0$. We may and do assume $\kappa=1$. By this and \eqref{rm3.5.2}, we get
\begin{eqnarray}
\ell^+F(0)+{\rm n}^{+}\left[\int_{0}^{\zeta}F(\epsilon_{s}){\rm d}s\right]
&=&\lim_{q\to 0+}\frac{\bP_{0}\left[F(\bar{\xi}_{\mathbf{e}_{q}}-\xi_{\mathbf{e}_{q}})\right]\Phi(q)}{q}\nonumber\\
&=&\lim_{q\to 0+}\frac{\hat{\bP}_{0}\left[F(\bar{\xi}_{\mathbf{e}_{q}})\right]}{\hat{\Phi}(q)}\nonumber\\
&=&\int_{\R^{+}}F(z)\hat{U}^{+}({\rm d}z)\nonumber
\end{eqnarray}
where $\hat{U}^{+}({\rm d}z):=\hat{\bP}_{0}\left[\int_{0}^{+\infty}\1_{\{\xi^{+}_{t}\in {\rm d}z\}}{\rm d}t\right]$. In the second equality we use the fact that $(\bar{\xi}_{\mathbf{e}_{q}}-\xi_{\mathbf{e}_{q}},\bP_{0})\stackrel{d}{=}(\bar{\xi}_{\mathbf{e}_{q}},\hat{\bP}_{0})$.
Setting $F(\cdot)=\Pi(y+\cdot,+\infty)$ in above equation and plugging it in \eqref{rm3.5.1} we get $$\Pi^{+}(y,+\infty)=\int_{\R^{+}}\Pi(z+y,+\infty)U^{-}({\rm d}z)$$ for $y>0$. This is Vigon's identity for L\'{e}vy process.
\end{remark}

\smallskip

Define
$$U^{+}_{\theta}({\rm d}v,{\rm d}z):=\bP_{0,\theta}\left[\int_{0}^{\bar{L}_{\infty}}\1_{\{\Theta^{+}_{s}\in {\rm d}v,\ \xi^{+}_{s}\in {\rm d}z\}}{\rm d}s\right].$$

\begin{proposition}\label{prop:under&overshoot}
Suppose $((\xi,\Theta),\bP)$ is a MAP with L\'{e}vy system $(H,\Pi)$ where $H_{t}=t\wedge \zeta$.
Then for any $x>0$, $\theta\in\s $ and any nonnegative measurable functions $f,g:\s \times \R^{+}\to \R^{+}$,
\begin{align}
\bP_{0,\theta}&\left[f(\Theta_{\tau^{+}_{x}-},x-\xi_{\tau^{+}_{x}-})g(\Theta_{\tau^{+}_{x}},\xi_{\tau^{+}_{x}}-x)\1_{\{\xi_{\tau^{+}_{x}}>x\}}\right]\nonumber\\
&=\int_{\s \times [0,x]}U^{+}_{\theta}({\rm d}v,{\rm d}z)\big[\ell^+(v)f(v,x-z)G(v,x-z)\nonumber\\
&\hspace{2cm}+{\rm n}^{+}_{v}\big(\int_{0}^{\zeta}f(\nu_{s},x-z+\epsilon_{s})G(\nu_{s},x-z+\epsilon_{s}){\rm d}s\big)\big],
\label{eq:under&overshoot}
\end{align}
where $G(v,u):=\int_{\s \times (u,+\infty)}g(\phi,y-u)\Pi(v,{\rm d}\phi,{\rm d}y)$ for $v\in\s $ and $u\in\R$.
In particular,
\begin{align}
\bP_{0,\theta}&\left[g(\Theta_{\tau^{+}_{x}},\xi_{\tau^{+}_{x}}-x)\1_{\{\xi_{\tau^{+}_{x}}>x\}}\right]\nonumber\\
&=\int_{\s \times [0,x]}U^{+}_{\theta}({\rm d}v,{\rm d}z)\int_{\s \times (x-z,+\infty)}g(\phi,z+y-x)\Pi^{+}(v,{\rm d}\phi,{\rm d}y).
\label{eq:overshoot}
\end{align}
\end{proposition}

\proof Let $\Delta\xi_{s}:=\xi_{s}-\xi_{s-}$ for any $s>0$. By \eqref{levysys} we have
\begin{align}\label{prop3.7.1}
\bP_{0,\theta}&\left[f(\Theta_{\tau^{+}_{x}-},x-\xi_{\tau^{+}_{x}-})g(\Theta_{\tau^{+}_{x}},\xi_{\tau^{+}_{x}}-x)\1_{\{\xi_{\tau^{+}_{x}}>x\}}\right]\nonumber\\
&=\bP_{0,\theta}\left[\sum_{s\ge 0}f(\Theta_{s-},x-\xi_{s-})g(\Theta_{s},\xi_{s-}+\Delta\xi_{s}-x)\1_{\{\bar{\xi}_{s-}\le x,\xi_{s-}+\Delta\xi_{s}-x>0\}}\right]\nonumber\\
&=\bP_{0,\theta}\left[\int_{0}^{\zeta}\1_{\{\bar{\xi}_{s}\le x\}}f(\Theta_{s},x-\xi_{s}){\rm d}s\int_{\s \times \R^{+}}g(v,\xi_{s}+y-x)\1_{\{\xi_{s}+y-x>0\}}\Pi(\Theta_{s},{\rm d}v,{\rm d}y)\right]\nonumber\\
&=\bP_{0,\theta}\left[\int_{0}^{\zeta}\1_{\{\bar{\xi}_{s}\le x\}}f(\Theta_{s},x-\xi_{s})G(\Theta_{s},x-\xi_{s}){\rm d}s\right].
\end{align}
We set $F(y,v):=f(v,x-y)G(v,x-y)$, then the right-hand side of \eqref{prop3.7.1} equals
\begin{align}
\bP_{0,\theta}&\left[\int_{0}^{\zeta}\1_{\{\bar{\xi}_{s}\le x\}}F(\xi_{s},\Theta_{s}){\rm d}s\right]\nonumber\\
&=\bP_{0,\theta}\left[\int_{0}^{\zeta}\1_{\{\bar{\xi}_{s}\le x,s\in \bar{M}^{cl}\}}F(\xi_{s},\Theta_{s}){\rm d}s\right]
+\bP_{0,\theta}\left[\int_{0}^{\zeta}\1_{\{\bar{\xi}_{s}\le x,\ s\not\in \bar{M}^{cl}\}}F(\xi_{s},\Theta_{s}){\rm d}s\right]\nonumber\\
&=\bP_{0,\theta}\left[\int_{0}^{+\infty}\1_{\{\bar{\xi}_{s}\le x\}}F(\xi_{s},\Theta_{s})\ell^+(\Theta_{s}){\rm d}\bar{L}_{s}\right]+
\bP_{0,\theta}\left[\sum_{g\in \bar{G}}\1_{\{\bar{\xi}_{g}\le x\}}\int_{g}^{d}F(\xi_{s},\Theta_{s}){\rm d}s\right].\nonumber
\end{align}
By \eqref{eq:last exit} the second term equals
$$\bP_{0,\theta}\left[\int_{0}^{+\infty}\1_{\{\bar{\xi}_{s}\le x\}}{\rm n}^{+}_{\Theta_{s}}\left(\int_{0}^{\zeta}F\left(\bar{\xi}_{s}-\epsilon_{r},\nu_{r}\right)\right){\rm d}\bar{L}_{s}\right].$$
Hence we have
\begin{align}
\bP_{0,\theta}&\left[f(\Theta_{\tau^{+}_{x}-},x-\xi_{\tau^{+}_{x}-})g(\Theta_{\tau^{+}_{x}},\xi_{\tau^{+}_{x}}-x)\1_{\{\xi_{\tau^{+}_{x}}>x\}}\right]\nonumber\\
&=\bP_{0,\theta}\left[\int_{0}^{+\infty}\1_{\{\bar{\xi}_{s}\le x\}}\left(\ell^+(\Theta_{s})F(\xi_{s},\Theta_{s})+{\rm n}^{+}_{\Theta_{s}}\left(\int_{0}^{\zeta}F\left(\bar{\xi}_{s}-\epsilon_{r},\nu_{r}\right)\right)\right){\rm d}\bar{L}_{s}\right]\nonumber\\
&=\bP_{0,\theta}\left[\int_{0}^{\bar{L}_{\infty}}\1_{\{\xi^{+}_{s}\le x\}}\left(\ell^+(\Theta^{+}_{s})F(\xi^{+}_{s},\Theta^{+}_{s})+{\rm n}^{+}_{\Theta^{+}_{s}}\left(\int_{0}^{\zeta}F\left(\xi^{+}_{s}-\epsilon_{r},\nu_{r}\right)\right)\right)
{\rm d}s\right]\nonumber\\
&=\int_{\s\times [0,x]}U^{+}_{\theta}({\rm d}v,{\rm d}z)\left(\ell^+(v)F(v,z)+{\rm n}^{+}_{v}\left(\int_{0}^{\zeta}F(z-\epsilon_{r},\nu_{r}){\rm d}r\right)\right),\nonumber
\end{align}
which yields \eqref{eq:under&overshoot}. \eqref{eq:overshoot} follows directly from \eqref{eq:under&overshoot} and Proposition \ref{prop:Pi+}.\qed

\smallskip

We say a path of $\xi$ creeps across level $x$ if it enters $(x,+\infty)$ continuously, that is, the first passage time in $(x,+\infty)$ is not a jump time. The next lemma we present is about what happens on the event of creeping. It follows from \cite[Proposition (1.5) and Theorem (1.7)]{Cinlar1}.

  \begin{lemma}\label{lem:on creeping}
  Suppose the ascending ladder height process
  $((\xi^{+},\Theta^{+}),\bP)$
  has a L\'{e}vy system $(H^{+},\Pi^{+})$ where $H^{+}_{t}=t\wedge \zeta^{+}$. If the continuous part of $\xi^{+}$ can be represented by $\int_{0}^{t\wedge \zeta^{+}}a^{+}(\Theta^{+}_{s}){\rm d}s$ for some nonnegative measurable function $a^{+}$ on $\s $, then for every $\theta\in\s $,
  $\1_{\{a^{+}(v)>0\}} U^+_\theta({\rm d}v,{\rm d}x)$ has a
  {kernel}
  $u^+_\theta({\rm d}v,x)$ with respect to the Lebesgue measure ${\rm d}x$.
  Moreover, if we define $T^{+}_{x}:=\inf\{t>0:\xi^{+}_{t}>x\}$, then for any nonnegative measurable function $f:\s \times \s \times \R^{+}\times \R^{+}\to \R^{+}$ and almost every $x>0$,
\begin{equation}\label{lem2.1.2}
\bP_{0,\theta}\left(\xi^{+}_{T^{+}_{x}-}<x=\xi^{+}_{T^{+}_{x}}\right)=0,
\end{equation}
and
\begin{equation}\label{lem2.1.1}
\bP_{0,\theta}\left[f\left(\Theta^{+}_{T^{+}_{x}-},\Theta^{+}_{T^{+}_{x}},x-\xi^{+}_{T^{+}_{x}-},\xi^{+}_{T^{+}_{x}}-x\right)
\1_{\{\xi^{+}_{T^{+}_{x}}=x\}}\right]
=\int_{\s }a^{+}(v)f(v,v,0,0)u^{+}_{\theta}({\rm d}v,x).
\end{equation}
\end{lemma}

\smallskip

\begin{lemma}\label{lem:2.2}
Suppose the MAP $((\xi,\Theta),\bP)$ has a L\'{e}vy system $(H,\Pi)$ where $H_{t}=t\wedge \zeta$.
If $(x,\theta)\in (0,+\infty)\times\s$ satisfies that
\begin{equation}\label{lem2.2.1}
\bP_{0,\theta}\left(\xi_{\tau^{+}_{x}-}<x=\xi_{\tau^{+}_{x}}\right)=0,
\end{equation}
then
$$\bP_{0,\theta}\left(\Theta_{\tau^{+}_{x}-}\not=\Theta_{\tau^{+}_{x}},\xi_{\tau^{+}_{x}}=x\right)=0.$$
\end{lemma}
\proof For $x>0$, let $\tau_{[x,+\infty)}$ denote the first time when $\xi$ enters $[x,+\infty)$. The upwards regularity of $((\xi,\Theta),\bP)$ implies that $\tau_{[x,+\infty)}=\tau^{+}_{x}$ $\bP_{0,\theta}$-a.s. It follows by \eqref{lem2.2.1} and \eqref{levysys} that
 \begin{eqnarray}
 \bP_{0,\theta}\left(\Theta_{\tau^{+}_{x}-}\not=\Theta_{\tau^{+}_{x}},\ \xi_{\tau^{+}_{x}}=x\right)
 &=&\bP_{0,\theta}\left(\Theta_{\tau^{+}_{x}-}\not=\Theta_{\tau^{+}_{x}},\ \xi_{\tau^{+}_{x}-}=\xi_{\tau^{+}_{x}}=x\right)\nonumber\\
 &=&\bP_{0,\theta}\left[\sum_{s\ge 0}\1_{\{\xi_{r}<x,\forall r\in [0,s),\ \Theta_{s-}\not=\Theta_{s},\ \xi_{s-}=\xi_{s}=x\}}\right]\nonumber\\
 &=&\bP_{0,\theta}\left[\int_{0}^{+\infty}\1_{\{\xi_{r}<x,\forall r\in [0,s),\ \xi_{s}=x\}}\Pi(\Theta_{s},\s\setminus\{\Theta_{s}\},\{0\}){\rm d}s\right]\nonumber\\
 &=&\bP_{0,\theta}\left[\int_{0}^{+\infty}\1_{\{\tau_{[x,+\infty)}=s,\ \xi_{s}=x\}}\Pi(\Theta_{s},\s\setminus\{\Theta_{s}\},\{0\}){\rm d}s\right]\nonumber\\
 &=&0.\nonumber
 \end{eqnarray}
The last equality is because the integral inside $\bP_{0,\theta}$ equals $0$.\qed

\begin{proposition}\label{prop:on creeping}
Suppose the MAP $((\xi,\Theta),\bP)$ has a L\'{e}vy system $(H,\Pi)$ where $H_{t}=t\wedge \zeta$ and the continuous part of $\xi^{+}$ can be represented by $\int_{0}^{t\wedge \zeta^{+}}a^{+}(\Theta^{+}_{s}){\rm d}s$ for some nonnegative measurable function $a^{+}$ on $\s$.
Then for every $\theta\in\s$, every nonnegative measurable function $f:\s\times \s \times \R^{+}\times \R^{+}\to \R^{+}$ and almost every $x>0$,
\begin{equation}
\bP_{0,\theta}\left[f(\Theta_{\tau^{+}_{x}-},\Theta_{\tau^{+}_{x}},x-\xi_{\tau^{+}_{x}-},\xi_{\tau^{+}_{x}}-x)\1_{\{\xi_{\tau^{+}_{x}}=x\}}\right]
=\int_{\s}a^{+}(v)f(v,v,0,0)u^{+}_{\theta}({\rm d}v,x),\label{prop3.8.1}
\end{equation}
where $u^{+}_{\theta}({\rm d}v,x)$ is the
{kernel}
given in Lemma \ref{lem:on creeping}.
\end{proposition}

\proof
It is easy to see from Proposition \ref{prop:Pi+} that the conditions of Lemma \ref{lem:on creeping} holds under the assumptions of this proposition.
Fix an arbitrary $\theta\in\s$. Let $\mathcal{R}$ denote the set of points for which both identities in Lemma \ref{lem:on creeping} hold. Then 
{$\mathrm{Leb}\left(\R^{+}\setminus\mathcal{R}\right)=0$.}
We note that $(\xi_{\tau^{+}_{x}},\Theta_{\tau^{+}_{x}})=(\xi^{+}_{T^{+}_{x}},\Theta^{+}_{T^{+}_{x}})$.
If we can prove
$\bP_{0,\theta}\left(\xi_{\tau^{+}_{x}-}<x=\xi_{\tau^{+}_{x}}\right)=0$
for every $x\in\mathcal{R}$, then by Lemma \ref{lem:2.2} $\Theta_{\tau^{+}_{x}-}=\Theta_{\tau^{+}_{x}}$ $\bP_{0,\theta}$-a.s. on $\{\xi_{\tau^{+}_{x}}=x\}$, and \eqref{prop3.8.1} is a direct consequence of \eqref{lem2.1.1}. Now fix an arbitrary $x\in\mathcal{R}$.
Let $\tau_{[x,+\infty)}$ denote the first time when $\xi$ enters $[x,+\infty)$. \eqref{lem2.1.2} implies that
$\xi^{+}_{T^{+}_{x}-}=x$ $\bP_{0,\theta}$-a.s. on the event $\{\xi_{\tau^{+}_{x}-}<x=\xi_{\tau^{+}_{x}}\}$, which in turn implies that
$\tau_{[x,+\infty)}<\tau^{+}_{x}$ $\bP_{0,\theta}$-a.s. on $\{\xi_{\tau^{+}_{x}-}<x=\xi_{\tau^{+}_{x}}\}$.
 Hence $\bP_{0,\theta}\left(\xi_{\tau^{+}_{x}-}<x=\xi_{\tau^{+}_{x}}\right)=0$, otherwise $\bP_{0,\theta}\left(\tau_{[x,+\infty)}<\tau^{+}_{x}\right)>0$, which contradicts the upwards regularity of $(\xi,\Theta)$. Hence we complete the proof.\qed

\smallskip

We note that the result in Proposition \ref{prop:on creeping} holds only for almost every $x>0$. In the following we give  sufficient conditions under which it holds for every $x>0$.

{
\begin{lemma}\label{lem2.11}
Suppose $((\xi,\Theta),\bP)$ is a MAP in $\R\times \s$ and $(X,\P)$ is the ssMp underlying  $((\xi,\Theta),\bP)$ via Lamperti-Kiu transform. Then for any $\theta\in\s$ and $r_{n},r\in\R$ such that $\lim_{n\to+\infty}r_{n}=r$, the process $(X,\P_{\mathrm{e}^{r_{n}}\theta})$ converges to $(X,\P_{\mathrm{e}^{r}\theta})$ in distribution under the Skorokhod topology.
\end{lemma}

\proof We need to show that for an arbitrary Lipschitz continuous function $f:\mathbb{D}_{\R^{d}}\to \R$,
\begin{equation}\nonumber
\lim_{n\to+\infty}\P_{\mathrm{e}^{r_{n}}\theta}\left[f(X)\right]=\P_{\mathrm{e}^{r}\theta}\left[f(X)\right].
\end{equation}
Suppose the ssMp $(X,\P)$ has index $\alpha>0$.
By the scaling property of $X$, it suffices to show that
\begin{equation}
\lim_{n\to+\infty}\P_{\theta}\left[f\left(\left(\mathrm{e}^{r_{n}}X_{\mathrm{e}^{-\alpha r_{n}}t}\right)_{t\ge 0}\right)\right]=\P_{\theta}\left[f\left(\left(\mathrm{e}^{r}X_{\mathrm{e}^{-\alpha r}t}\right)_{t\ge 0}\right)\right].\label{lemn1.1}
\end{equation}
We use $\mathrm{d}(\cdot,\cdot)$ to denote the Prokhorov's metric in $\mathbb{D}_{\R^{d}}$ which is compatible with the Skorokhod convergence.
It follows from \cite[Proposition 3.5.3(c)]{EthKur} that for any $\omega_{n},\omega_{0}\in\mathbb{D}_{\R^{d}}$, one has $\lim_{n\to+\infty}\mathrm{d}(\omega_{n},\omega_{0})=0$ if and only if for every $T\in (0,+\infty)$, there exists a sequence of strictly increasing continuous functions $\{\lambda_{n}:[0,T]\to \R^{+},\ n\ge 1\}$ with $\lambda_{n}(0)=0$, such that
\begin{equation}\nonumber
\lim_{n\to+\infty}\sup_{t\in [0,T]}\left(\|\omega_{n}(t)-\omega_{0}\circ \lambda_{n}(t)\|\vee |\lambda_{n}(t)-t|\right)=0.
\end{equation}
For an arbitrary $\omega\in\mathbb{D}_{\R^{d}}$, by setting $\omega_{n}(t)=\mathrm{e}^{r_{n}}\omega\left(\mathrm{e}^{-\alpha r_{n}}t\right)$, $\omega_{0}(t)=\mathrm{e}^{r}\omega\left(\mathrm{e}^{-\alpha r}t\right)$ and $\lambda_{n}(t)=\mathrm{e}^{\alpha(r-r_{n})}t$ for all $t\ge 0$, one can easily show that
$$\sup_{t\in [0,T]}\left(\|\omega_{n}(t)-\omega_{0}\circ \lambda_{n}(t)\|\vee |\lambda_{n}(t)-t|\right)=\sup_{t\in [0,T]}\left(|\mathrm{e}^{r_{n}}-\mathrm{e}^{r}|\,\|\omega(\mathrm{e}^{-\alpha r_{n}}t)\|\vee |\mathrm{e}^{\alpha(r-r_{n})}-1|t\right)\to 0$$
as $n\to+\infty$, and hence $\mathrm{d}(\omega_{n},\omega_{0})\to 0$. This implies that the processes $(\mathrm{e}^{r_{n}}X_{\mathrm{e}^{-\alpha r_{n}}t})_{t\ge 0}$ converges to $(\mathrm{e}^{r}X_{\mathrm{e}^{-\alpha r}t})_{t\ge 0}$ $\P_{\theta}$-almost surely under the Skorokhod topology. Therefore \eqref{lemn1.1} follows from this and the bounded convergence theorem.\qed

\smallskip
}

\begin{proposition}\label{prop:sufficient}
Suppose the conditions in Proposition \ref{prop:on creeping} hold.
Then for every $\theta\in\s$ and every $x>0$,
$$\bP_{0,\theta}\left(\xi_{\tau^{+}_{x}-}<x\mbox{ or }\Theta_{\tau^{+}_{x}-}\not=\Theta_{\tau^{+}_{x}};\xi_{\tau^{+}_{x}}=x\right)=0,$$
and for every bounded continuous function $g:\R^{+}\times \s\times\s\times \R^{+}\times \R^{+}\to \R$, the function
$$x\mapsto \bP_{0,\theta}\left[g(\tau^{+}_{x},\Theta_{\tau^{+}_{x}-},\Theta_{\tau^{+}_{x}},x-\xi_{\tau^{+}_{x}-},\xi_{\tau^{+}_{x}}-x)\1_{\{\xi_{\tau^{+}_{x}}=x\}}\right]$$ is right continuous on $[0,+\infty)$.
{If, in addition, $a^{+}(v)>0$ for every $v\in\s$ or if $a^{+}(v)=0$ for every $v\in\s$, then}
{the kernel}
$u^{+}_{\theta}({\rm d}v,x)$ of
{$\1_{\{a^{+}(v)>0\}}U^{+}_{\theta}({\rm d}v,{\rm d}x)$} can take a unique version such that $x\mapsto a^{+}(v)u^{+}_{\theta}({\rm d}v,x)$ is right continuous on $(0,+\infty)$ in the sense of vague convergence. In this case, \eqref{prop3.8.1} holds for every $x>0$ and every nonnegative measurable function $f:\s\times\s\times \R^{+}\times \R^{+}\to \R^{+}$.
\end{proposition}

\proof For every $(x,\theta)\in\R^{+}\times \s$, let $p^{\theta}(x):=\bP_{0,\theta}\left(\xi_{\tau^{+}_{x}}=x\right)$,
$p^{\theta}_{1}(x):=\bP_{0,\theta}\left(\xi_{\tau^{+}_{x}-}=\xi_{\tau^{+}_{x}}=x\right)$ and $p^{\theta}_{2}(x):=p^{\theta}(x)-p^{\theta}_{1}(x)=\bP_{0,\theta}\left(\xi_{\tau^{+}_{x}-}<x=\xi_{\tau^{+}_{x}}\right)$. By Proposition \ref{prop:on creeping} we have $p^{\theta}_{2}(x)=0$ for almost every $x>0$.
{It follows by Proposition \ref{prop:under&overshoot} 
that}
\begin{eqnarray}
\bP_{0,\theta}\left(\xi_{\tau^{+}_{x}}>x\right)
&=&\int_{\mathcal{S}\times [0,x]}\bar{\Pi}^{+}_{v}(x-z)U^{+}_{\theta}({\rm d}v,{\rm d}z)
.\nonumber
\end{eqnarray}
Here
$\bar{\Pi}^{+}_{v}(u)=\Pi^{+}(v,\s,(u,+\infty))$.
Obviously from the above equation $x\mapsto p^{\theta}(x)=1-\bP_{0,\theta}\left(\xi_{\tau^{+}_{x}}>x\right)$ is right continuous on $[0,+\infty)$. Suppose $x_{n},x\in\R^{+}$ and $x_{n}\downarrow x$.
{Let $(X,\P)$ denote the ssMp underlying  $((\xi,\Theta),\bP)$ via Lamperti-Kiu transform. It follows by Lemma \ref{lem2.11} that}
$$\left(X,\P_{\theta\rme^{-x_{n}}}\right)\to \left(X,\P_{\theta\rme^{-x}}\right)$$
in distribution under the Skorokhod topology. For $n\ge 1$, let $(Y^{(n)},\P^{*})$ and $(Y ,\P^{*})$ be couplings of $(X,\P_{\theta\rme^{-x_{n}}})$ and $(X,\P_{\theta\rme^{-x}})$ respectively, such that $Y^{(n)}\rightarrow Y $ $\P^{*}$-a.s. in the Skorokhod topology.  Let $\varsigma_{0}:=\inf\{t\ge 0:\ \|Y_{t}\|>1\}$ and $\varsigma_{n}:=\inf\{t\ge 0:\ \|Y^{(n)}_{t}\|>1\}$ for $n\ge 0$.
Since $X$ is sphere-exterior regular, so is $Y$, which implies that $\|Y_{t}\|\not=1$ for any $t<\varsigma_{0}$ $\P^{*}$-a.s. In view of this, it follows by \cite[Theorem 13.6.4]{Whitt} that
$$(Y^{(n)}_{\varsigma_{n}-},Y^{(n)}_{\varsigma_{n}})\to (Y_{\varsigma_{0}-},Y_{\varsigma_{0}})\quad \P^{*}\mbox{-a.s.}$$
as $n\to+\infty$. Hence $\left((X_{\tau^{\ominus}_{1}-},X_{\tau^{\ominus}_{1}}),\P_{\theta\rme^{-x_{n}}}\right)$ converges in distribution to $\left((X_{\tau^{\ominus}_{1}-},X_{\tau^{\ominus}_{1}}),\P_{\theta\rme^{-x}}\right)$. This weak convergence yields that
\begin{eqnarray}
p^{\theta}_{1}(x)&=&\bP_{-x,\theta}\left(\xi_{\tau^{+}_{0}-}=\xi_{\tau^{+}_{0}}=0\right)\nonumber\\
&=&\P_{\theta\rme^{-x}}\left(X_{\tau^{\ominus}_{1}-}\in\mathbb{S}^{d-1},X_{\tau^{\ominus}_{1}}\in \mathbb{S}^{d-1}\right)\nonumber\\
&\ge&\limsup_{n\to+\infty}\P_{\theta\rme^{-x_{n}}}\left(X_{\tau^{\ominus}_{1}-}\in\mathbb{S}^{d-1},X_{\tau^{\ominus}_{1}}\in \mathbb{S}^{d-1}\right)\nonumber\\
&=&\limsup_{n\to+\infty}p^{\theta}_{1}(x_{n}).\nonumber
\end{eqnarray}
This and the right continuity of $p^{\theta}(\cdot)$ imply that
$\liminf_{n\to+\infty}p^{\theta}_{2}(x_{n})\ge p^{\theta}_{2}(x).$
Hence
\begin{equation}
p^{\theta}_{2}(x)=\bP_{0,\theta}\left(\xi_{\tau^{+}_{x}-}<x=\xi_{\tau^{+}_{x}}\right)=0\quad \forall x>0.\label{prop2.11.1}
\end{equation}
It then follows by Lemma \ref{lem:2.2} that
\begin{equation}\label{prop2.11.2}
\bP_{0,\theta}\left(\Theta_{\tau^{+}_{x}-}\not=\Theta_{\tau^{+}_{x}},\xi_{\tau^{+}_{x}}=x\right)=0,\quad\forall x>0.
\end{equation}
We need to show that
\begin{align}\label{prop2.11.3}
\lim_{n\to+\infty}&\bP_{0,\theta}\left[g(\tau^{+}_{x_{n}},\Theta_{\tau^{+}_{x_{n}}-},\Theta_{\tau^{+}_{x_{n}}},x_{n}-\xi_{\tau^{+}_{x_{n}}-},\xi_{\tau^{+}_{x_{n}}}-x_{n})\1_{\{\xi_{\tau^{+}_{x_{n}}}=x_{n}\}}\right]\nonumber\\
&=\bP_{0,\theta}\left[g(\tau^{+}_{x},\Theta_{\tau^{+}_{x}-},\Theta_{\tau^{+}_{x}},x-\xi_{\tau^{+}_{x}-},\xi_{\tau^{+}_{x}}-x)\1_{\{\xi_{\tau^{+}_{x}}=x\}}\right]
\end{align}
for any sequence $x_{n},x\in\R^{+}$, $x_{n}\downarrow x$ and any bounded continuous function $g:\R^{+}\times \s\times\s\times\R^{+}\times\R^{+}\to\R$.
Let $A_{n}:=\{\xi_{\tau^{+}_{x_{n}}}=x_{n}\}$ and $A:=\{\xi_{\tau^{+}_{x}}=x\}$. By the strong Markov property and the fact that $\lim_{y\to0+}p^{v}(y)=p^{v}(0)=1$ for every $v\in \s$, we have for every $\theta\in\s$
\begin{align}
\bP_{0,\theta}&\left(A\setminus A_{n}\right)=\bP_{0,\theta}\left(\xi_{\tau^{+}_{x}}=x,\xi_{\tau^{+}_{x_{n}}}>x_{n}\right)\nonumber\\
&=\bP_{0,\theta}\left(\bP_{0,\Theta_{\tau^{+}_{x}}}\left(\xi_{\tau^{+}_{x_{n}-x}}>x_{n}-x\right);\xi_{\tau^{+}_{x}}=x\right)\nonumber\\
&=\bP_{0,\theta}\left[\left(1-p^{\Theta_{\tau^{+}_{x}}}(x_{n}-x)\right)\1_{\{\xi_{\tau^{+}_{x}}=x\}}\right]\nonumber\\
&\to0,\quad \mbox{ as }n\to+\infty.\nonumber
\end{align}
Since $\bP_{0,\theta}\left(A_{n}\setminus A\right)-\bP_{0,\theta}\left(A\setminus A_{n}\right)=\bP_{0,\theta}(A_{n})-\bP_{0,\theta}(A)=p^{\theta}(x_{n})-p^{\theta}(x)\to 0$ as $n\to+\infty,$
we have
\begin{equation}
\bP_{0,\theta}\left(A\triangle A_{n}\right)=\bP_{0,\theta}\left(A_{n}\setminus A\right)+\bP_{0,\theta}\left(A\setminus A_{n}\right)\to 0\quad \mbox{ as }n\to+\infty.\label{prop2.11.4}
\end{equation}
Note that by \eqref{prop2.11.1} and \eqref{prop2.11.2}
\begin{align}
&\big|\bP_{0,\theta}\big[g(\tau^{+}_{x_{n}},\Theta_{\tau^{+}_{x_{n}}-},\Theta_{\tau^{+}_{x_{n}}},x_{n}-\xi_{\tau^{+}_{x_{n}}-},\xi_{\tau^{+}_{x_{n}}}-x_{n})\1_{\{\xi_{\tau^{+}_{x_{n}}}=x_{n}\}}\big] \nonumber\\
&\quad-\bP_{0,\theta}\big[g(\tau^{+}_{x},\Theta_{\tau^{+}_{x}-},\Theta_{\tau^{+}_{x}},x-\xi_{\tau^{+}_{x}-},\xi_{\tau^{+}_{x}}-x)\1_{\{\xi_{\tau^{+}_{x}}=x\}}\big]\big|\nonumber\\
&=\big|\bP_{0,\theta}\big[g({\tau^{+}_{x_{n}},}\Theta_{\tau^{+}_{x_{n}}},\Theta_{\tau^{+}_{x_{n}}},x_{n}-\xi_{\tau^{+}_{x_{n}}},\xi_{\tau^{+}_{x_{n}}}-x_{n})\1_{\{\xi_{\tau^{+}_{x_{n}}}=x_{n}\}}\big] \nonumber\\
&\quad-\bP_{0,\theta}\big[g({\tau^{+}_{x},}\Theta_{\tau^{+}_{x}},\Theta_{\tau^{+}_{x}},x-\xi_{\tau^{+}_{x}},\xi_{\tau^{+}_{x}}-x)\1_{\{\xi_{\tau^{+}_{x}}=x\}}\big]\big|\nonumber\\
&\le\big|\bP_{0,\theta}\big[g(\tau^{+}_{x_{n}},\Theta_{\tau^{+}_{x_{n}}},\Theta_{\tau^{+}_{x_{n}}},x_{n}-\xi_{\tau^{+}_{x_{n}}},\xi_{\tau^{+}_{x_{n}}}-x_{n})\big(\1_{\{\xi_{\tau^{+}_{x_{n}}}=x_{n}\}}
-\1_{\{\xi_{\tau^{+}_{x}}=x\}}\big)\big]\big|\nonumber\\
&\quad+\big|\bP_{0,\theta}\big[\big(g(\tau^{+}_{x_{n}},\Theta_{\tau^{+}_{x_{n}}},\Theta_{\tau^{+}_{x_{n}}},x_{n}-\xi_{\tau^{+}_{x_{n}}},\xi_{\tau^{+}_{x_{n}}}-x_{n})-g(\tau^{+}_{x},\Theta_{\tau^{+}_{x}},\Theta_{\tau^{+}_{x}},x-\xi_{\tau^{+}_{x}},\xi_{\tau^{+}_{x}}-x)\big)\1_{\{\xi_{\tau^{+}_{x}}=x\}}
\big]\big|\nonumber\\
&\le\|g\|_{\infty}\bP_{0,\theta}\big(A\triangle A_{n}\big)\nonumber\\
&\quad+\bP_{0,\theta}\big[\big|g(\tau^{+}_{x_{n}},\Theta_{\tau^{+}_{x_{n}}},\Theta_{\tau^{+}_{x_{n}}},x_{n}-\xi_{\tau^{+}_{x_{n}}},\xi_{\tau^{+}_{x_{n}}}-x_{n})
-g(\tau^{+}_{x},\Theta_{\tau^{+}_{x}},\Theta_{\tau^{+}_{x}},x-\xi_{\tau^{+}_{x}},\xi_{\tau^{+}_{x}}-x)\big|\big].\nonumber
\end{align}
We have $\tau^{+}_{x_{n}}\downarrow \tau^{+}_{x}$ $\bP_{0,\theta}$-a.s. by the upwards regularity of $(\xi,\Theta)$ and hence $\left(\Theta_{\tau^{+}_{x_{n}}},\xi_{\tau^{+}_{x_{n}}}\right)
\to \left(\Theta_{\tau^{+}_{x}},\xi_{\tau^{+}_{x}}\right)$ $\bP_{0,\theta}$-a.s. by the right continuity of $(\xi,\Theta)$. In view of this and \eqref{prop2.11.4}, \eqref{prop2.11.3} follows by letting $n\to+\infty$ in the above inequality.

\smallskip

By \eqref{prop2.11.1} and \eqref{prop2.11.2}, we have for every $x>0$ and every nonnegative measurable function $f:\s\times\s\times\R^{+}\times\R^{+}\to\R^{+}$,
\begin{eqnarray}
\bP_{0,\theta}\left[f(\Theta_{\tau^{+}_{x}-},\Theta_{\tau^{+}_{x}},x-\xi_{\tau^{+}_{x}-},\xi_{\tau^{+}_{x}}-x)\1_{\{\xi_{\tau^{+}_{x}}=x\}}\right]
&=&\bP_{0,\theta}\left[f(\Theta_{\tau^{+}_{x}},\Theta_{\tau^{+}_{x}},0,0)\1_{\{\xi_{\tau^{+}_{x}}=x\}}\right]\nonumber\\
&=&\int_{\s}f(v,v,0,0)\bP_{0,\theta}\left(\Theta_{\tau^{+}_{x}}\in {\rm d}v,\xi_{\tau^{+}_{x}}=x\right).\nonumber\\
\label{inviewof}
\end{eqnarray}
{ Let us momentarily assume that $a^+(v) >0$ for all $v\in \s$.} In view of \eqref{inviewof}, and Proposition \ref{prop:on creeping},
we can set the
{kernel}
$u^{+}_{\theta}({\rm d}v,x)$ of $U^{+}_{\theta}({\rm d}v,{\rm d}x)$ to be $\frac{1}{a^{+}(v)}\bP_{0,\theta}\left(\Theta_{\tau^{+}_{x}}\in {\rm d}v,\xi_{\tau^{+}_{x}}=x\right)$ for every $x>0$, in which case, the function
$$x\mapsto a^{+}(v)u^{+}_{\theta}({\rm d}v,x)=\bP_{0,\theta}\left(\Theta_{\tau^{+}_{x}}\in {\rm d}v,\xi_{\tau^{+}_{x}}=x\right)$$ is right continuous on $(0,+\infty)$ in the sense of vague convergence, because
$x\mapsto \bP_{0,\theta}\left[h(\Theta_{\tau^{+}_{x}});\xi_{\tau^{+}_{x}}=x\right]$ is right continuous on $(0,+\infty)$ for every bounded continuous function $h:\s\to\R$.
{For the other case that $a^+(v) = 0$ for all $v\in\s$, there is nothing to prove as, irrespective of our choice of $u^+_\theta$, the quantity  $a^{+}(v)u^{+}_{\theta}({\rm d}v,x)$ is right continuous as claimed (in fact it is identically equal to zero).} \qed

\subsection{Long time behavior of MAPs}\label{sec:slln}

It is well-known that for any $\R$-valued L\'{e}vy process $\chi$ one has $\chi_{t}/t\to \mathrm{E}\chi_{1}$ almost surely whenever $\mathrm{E}\chi_{1}$ is well-defined. Its proof relies on the classical strong law of large numbers. Following this, a L\'{e}vy process exhibits exactly one of the following behaviors: $\lim_{t\to+\infty}\chi_{t}=+\infty$ a.s., $\lim_{t\to+\infty}\chi_{t}=-\infty$ a.s. and  $\limsup_{t\to+\infty}\chi_{t}=-\liminf_{t\to+\infty}\chi_{t}=+\infty$ a.s. according as $\mathrm{E}\chi_{1}>,\ <,\ =0$.
This basic trichotomy is also true for the MAPs where $(\Theta_{t})_{t\ge 0}$ is a positive recurrent Markov process on a countable state space. We refer to \cite{Alsmeyer-Buckmann} and the references therein. In such case, let $\tau_{0}(i):=0$ and $\{\tau_{n}(i):n\ge 1\}$ denote the renewal sequence of successive return times to each state $i\in\mathcal{S}$. Then for each $i$, $\{\xi_{\tau_{n}(i)}:n\ge 0\}$ constitutes an ordinary random walk. In fact, a law of large numbers can be obtained by applying known results for these embedded random walks, but with considerable additional analysis. Regarding the more general situation when the modulator $\Theta$ has an uncountably infinite state space, we note that a natural substitute for $\{\tau_{n}(i):n\ge 1\}$ is a sequence of random times $\{R_{n}:n\ge 0\}$, in terms of which the process can be decomposed into independent and stationary blocks. In order to construct such random times, we assume the MAP satisfies the following Harris-type condition: There exist a constant $\delta>0$, a probability measure $\rho$ on $\mathcal{S}$ and a family of measures $\{\phi(\theta,\cdot):\theta\in\s \}$ on $\R$ with $\inf_{\theta\in\mathcal{S}}\phi(\theta,\R)>0$ such that
    \begin{equation}\label{condi:Harris type}\tag{HT}
    \bP_{0,\theta}\left(\xi_{\delta}\in \Gamma,\, \Theta_{\delta}\in A\right)\ge \phi(\theta,\Gamma)\rho(A)\quad\forall\theta\in\mathcal{S},\ A\in\mathcal{B}(\mathcal{S}),\ \Gamma\in\mathcal{B}(\R).
    \end{equation}
This section aims at providing the trichotomy regarding the almost sure behavior of $\xi_{t}$ as $t\to+\infty$ when condition \eqref{condi:Harris type} is satisfied.

 \smallskip

Define $M_{0}:=\Theta_{0}$, $S_{0}:=\xi_{0}$ and for any $n\ge 1$, define
    \begin{equation}\nonumber
    M_{n}:=\Theta_{n\delta},\ \Delta_{n}:=\xi_{n\delta}-\xi_{(n-1)\delta}\quad\mbox{and}\quad S_{n}:=S_{0}+\sum_{k=1}^{n}\Delta_{k}.
    \end{equation}
It is easy to verify that $((S_{n},M_{n})_{n\ge 0},\bP)$ is a discrete-time MAP satisfying
    \begin{equation}
    \bP_{0,\theta}\left(\Delta_{1}\in\Gamma,\ M_{1}\in A\right)\ge \phi(\theta,\Gamma)\rho(A)\label{condi:slln}
    \end{equation}
for all $\theta\in\mathcal{S},A\in\mathcal{B}(\s )$ and $\Gamma\in\mathcal{B}(\R)$.
    In particular we have
    \begin{equation}\nonumber
    \bP_{0,\theta}\left(M_{1}\in A\right)\ge \epsilon\rho(A)\quad\forall \theta\in\mathcal{S},\ A\in\mathcal{B}(\mathcal{S}),
    \end{equation}
    where $\epsilon:=\inf_{\theta\in\mathcal{S}}\phi(\theta,\R)>0$. This implies that $\{M_{n}:n\ge 0\}$ is an irreducible and strongly aperiodic Harris recurrent chain on $\mathcal{S}$. Given this and \eqref{condi:slln},
    it follows by \cite{Ney&NummelinI,Ney&NummelinII} that there exists a sequence of regeneration times $0\le R_{0}<R_{1}<\cdots<+\infty$ such that $\{R_{n+1}-R_{n}:n\ge 0\}$ is a sequence of independent and identically distributed nonnegative random variables, and that the random blocks $\{M_{R_{n}},\cdots,M_{R_{n+1}-1},\Delta_{R_{n}+1},\cdots,\Delta_{R_{n+1}}\}$ are independent with
    $$\bP_{0,\theta}\left[M_{R_{n}}\in A\,|\,\mathcal{G}_{R_{n-1}},\Delta_{R_{n}}\right]=\rho(A)\quad\forall A\in\mathcal{B}(\s),$$
    where $\mathcal{G}_{k}$ denotes the $\sigma$-field generated by $\{M_{0},\cdots,M_{k},\Delta_{1},\cdots,\Delta_{k}\}$.

   \smallskip

    We assume that $(\Theta_{t})_{t\ge 0}$ has an invariant distribution $\pi$.
    By \cite[Theorem 3.2]{AsmussenQueue} $\pi$ is uniquely determined by
    \begin{equation}
    \pi(A)=\frac{1}{\bP_{0,\rho}[R_{1}]}\bP_{0,\rho}\left[\sum_{j=0}^{R_{1}-1}\1_{\{M_{j}\in A\}}\right]\quad\forall A\in\mathcal{B}(\s )\label{stationarity}
    \end{equation}
    where $0<\bP_{0,\rho}[R_{1}]<+\infty$. It follows that
    \begin{eqnarray}
    \bP_{0,\pi}\left[S_{1}\right]
    &=&\frac{1}{\bP_{0,\rho}[R_{1}]}\sum_{j=0}^{+\infty}\int_{\s }\bP_{0,\rho}\left[S_{j+1}-S_{j}\,|\,M_{j}=\theta\right]\bP_{0,\rho}\left(M_{j}\in {\rm d}\theta,\,j\le R_{1}-1\right)\nonumber\\
    &=&\frac{1}{\bP_{0,\rho}[R_{1}]}\bP_{0,\rho}\left[\sum_{j=0}^{R_{1}-1}\left(S_{j+1}-S_{j}\right)\right]\nonumber\\
    &=&\frac{1}{\bP_{0,\rho}[R_{1}]}\bP_{0,\rho}\left[S_{R_{1}}\right],\label{eq:mean S1}
    \end{eqnarray}
    whenever $\bP_{0,\pi}\left[|S_{1}|\right]<+\infty$.
     The regeneration structure implies that $\left(S_{R_{n+1}}-S_{R_{n}}\right)$ is independent of $\{S_{k},k\le R_{n}\}$, and its distribution is independent of $n$.
     Let $N_{n}:=\sup\{k:\ R_{k}\le n\}$. We can write
    $$S_{n}=S_{R_{0}\wedge n}+\left[\left(S_{R_{1}}-S_{R_{0}}\right)+\cdots+\left(S_{R_{N_{n}}}-S_{R_{N_{n}-1}}\right)\right]+\left(S_{n}-S_{R_{N_{n}}}\right).$$
    It is easy to see that $S_{R_{0}\wedge n}/n\to 0$ almost surely, since $R_{0}$ is finite and $\lim_{n\to +\infty}S_{R_{0}\wedge n}=S_{R_{0}}<+\infty$ almost surely. Note that $\left(S_{R_{1}}-S_{R_{0}}\right)+\cdots+\left(S_{R_{N_{n}}}-S_{R_{N_{n}-1}}\right)$ is a random sum of i.i.d summands. In view of \eqref{eq:mean S1}, we have by the standard LLN and the elementary renewal theory that
    \begin{eqnarray}
    &&\lim_{n\to+\infty}\frac{\left(S_{R_{1}}-S_{R_{0}}\right)+\cdots+\left(S_{R_{N_{n}}}-S_{R_{N_{n}-1}}\right)}{n}\nonumber\\
    &=&\lim_{n\to+\infty}\frac{\left(S_{R_{1}}-S_{R_{0}}\right)+\cdots+\left(S_{R_{N_{n}}}-S_{R_{N_{n}-1}}\right)}{N_{n}}\cdot\frac{N_{n}}{n}\nonumber\\
    &=& \bP_{0,\rho}\left[S_{R_{1}}\right]\cdot \frac{1}{\bP_{0,\rho}[R_{1}]}=\bP_{0,\pi}\left[S_{1}\right]\quad \bP_{0,\theta}\mbox{-a.s.}\nonumber
    \end{eqnarray}
    Moreover, one can easily show by Borel-Cantelli lemma that $\left(S_{n}-S_{R_{N_{n}}}\right)/n\to 0$ $\bP_{0,\theta}$-a.s. if
    $$\bP_{0,\rho}\left[\max_{1\le k\le R_{1}}|S_{k}|\right]<+\infty.$$
    We have hence proved the following lemma.

    \begin{lemma}\label{lem:(A.1)}
    If $\bP_{0,\rho}\left[\max_{1\le k\le R_{1}}|S_{k}|\right]<+\infty$,
    then $\lim_{n\to+\infty}S_{n}/n=\bP_{0,\pi}[S_{1}]$ $\bP_{0,\theta}$-a.s. for every $\theta\in\mathcal{S}$.
    \end{lemma}

    \begin{lemma}\label{lem:A.4}
  If $\bP_{0,\pi}\left[\sup_{s\in[0,t]}|\xi_{s}|\right]$ is finite for some $t>0$, then it is finite for all $t>0$. Moreover,  $\bP_{0,\pi}\left[\sup_{s\in[0,\mathbf{e}_{q}]}|\xi_{s}|\right]$ is finite for all $q>0$.
  \end{lemma}
  \proof In this proof we use $\|\xi\|_{t}$ to denote $\sup_{s\in[0,t]}|\xi_{s}|$. Let $f(t):=\bP_{0,\pi}\left[\|\xi\|_{t}\right]$ for $t\ge 0$. We observe that for any $t,r>0$,
  \begin{equation}
  \|\xi\|_{t+r}\le \|\xi\|_{t}\vee\left(\sup_{s\in[t,t+r]}|\xi_{s}-\xi_{t}|+|\xi_{t}|\right)\le \|\xi\|_{t}+\sup_{s\in[t,t+r]}|\xi_{s}-\xi_{t}|.\label{lemA.4.1}
  \end{equation}
  By the Markov property and translation invariance in $\xi$, we have
  $$\bP_{0,\pi}\left[\sup_{s\in[t,t+r]}|\xi_{s}-\xi_{t}|\right]=\bP_{0,\pi}\left[\bP_{0,\Theta_{t}}\left[\|\xi\|_{r}\right]\right]
  =\bP_{0,\pi}\left[\|\xi\|_{r}\right]=f(r).$$
  The second equality is because $\pi$ is an invariant distribution of $(\Theta_{t})_{t\ge 0}$. Hence by \eqref{lemA.4.1} we get $f(t+r)\le f(t)+f(r)$.
   Given that $f(t)$ is finite for some $t>0$,
   $f$ is a nonnegative locally bounded subadditive function on $[0,+\infty)$. Hence there exist some constants $b,c>0$ such that $f(t)\le ct+b$ for all $t>0$. Consequently, $\bP_{0,\pi}\left[\|\xi\|_{\mathbf{e}_{q}}\right]=\int_{0}^{+\infty}q\mathrm{e}^{-q t}f(t){\rm d}t<+\infty$ for all $q>0$.\qed

\smallskip

\begin{proposition}\label{thm:slln}
    Suppose $((\xi,\Theta),\bP)$ is a MAP satisfying \eqref{condi:Harris type} and $\pi$ is an invariant distribution for $(\Theta_{t})_{t\ge 0}$. If
    \begin{equation}\label{condi:finite mean1}
    \bP_{0,\pi}\left[\sup_{s\in[0,1]}|\xi_{s}|\right]<+\infty
    \end{equation}
    then $\xi_{t}/t\to \bP_{0,\pi}[\xi_{1}]$ $\bP_{0,\theta}$-a.s. for  every $\theta\in\mathcal{S}$.
    \end{proposition}

 \proof Without loss of generality we assume that \eqref{condi:Harris type} holds for $\delta=1$. This proof works through for any $\delta>0$ with minor modifications. By Lemma \ref{lem:A.4}, condition \eqref{condi:finite mean1} implies that $\bP_{0,\pi}\left[\sup_{s\in[0,t]}|\xi_{s}|\right]<+\infty$ for all $t>0$ and $\bP_{0,\pi}[|\Delta_{1}|]=\bP_{0,\pi}[|\xi_{1}|]<+\infty$. We have
    \begin{eqnarray}
    \bP_{0,\rho}\left[\max_{1\le k\le R_{1}}|S_{k}|\right]
    &\le&\bP_{0,\rho}\left[\sum_{j=0}^{R_{1}-1}|\Delta_{j+1}|\right]\nonumber\\
    &=&\sum_{j=0}^{+\infty}\int_{\s }\bP_{0,\rho}\left[\left|\Delta_{j+1}\right|\,|\,M_{j}=\theta\right]\bP_{0,\rho}\left(M_{j}\in {\rm d}\theta,j\le R_{1}-1\right)\nonumber\\
    &=&\int_{\s }\bP_{0,\theta}\left[|\Delta_{1}|\right]\bP_{0,\rho}\left[\sum_{j=0}^{R_{1}-1}\1_{\{M_{j}\in {\rm d}\theta\}}\right]\nonumber\\
    &=&\bP_{0,\rho}[R_{1}]\bP_{0,\pi}[|\Delta_{1}|]<+\infty,\label{propA.2.0}
    \end{eqnarray}
    where in the last equality we use \eqref{stationarity}.
    It follows by Lemma \ref{lem:(A.1)} that $S_{n}/n\to \bP_{0,\pi}[S_{1}]=\bP_{0,\pi}[\xi_{1}]$ $\bP_{0,\theta}$-a.s. for  every $\theta\in\mathcal{S}$.
    Note that for any $t\in [ R_{k}, R_{k+1})$,
    \begin{equation}\nonumber
    \frac{S_{R_{k}}}{R_{k}}\,\frac{R_{k}}{ R_{k+1}}-\frac{\sup_{s\in [ R_{k}, R_{k+1}]}|\xi_{s}-S_{R_{k}}|}{ R_{k+1}}\le \frac{\xi_{t}}{t}\le
    \frac{S_{R_{k}}}{R_{k}}+\frac{\sup_{s\in [ R_{k}, R_{k+1}]}|\xi_{s}-S_{R_{k}}|}{ R_{k}}.
    \end{equation}
     It is known by the renewal theorem that $R_{k}/k\to \bP_{0,\rho}[R_{1}]$ $\bP_{0,\theta}$-a.s. Hence to prove $\xi_{t}/t\to \bP_{0,\pi}\left[\xi_{1}\right]$ $\bP_{0,\theta}$-a.s., it suffices to prove that
    \begin{equation}\label{propA.2.1}
    \frac{\sup_{s\in [R_{k},R_{k+1}]}|\xi_{s}-S_{R_{k}}|}{k}\to 0\quad \mbox{ as }k\to +\infty\quad\bP_{0,\theta}\mbox{-a.s.}
    \end{equation}
    for every $\theta\in\mathcal{S}$. The regeneration structure implies that
   $\{\sup_{s\in [R_{k},R_{k+1}]}|\xi_{s}-S_{R_{k}}|:k\ge 1\}$ under $\bP_{0,\theta}$
    is a family of i.i.d. random variables which have the same distribution as $\left(\sup_{s\in [0,R_{1}]}|\xi_{s}|,\bP_{0,\rho}\right)$. Hence by the second Borel-Cantelli lemma, \eqref{propA.2.1} holds if and only if
   \begin{equation}
   \bP_{0,\rho}\left[\sup_{s\in [0,R_{1}]}|\xi_{s}|\right]<+\infty.\label{propA.2.2}
   \end{equation}
   We note that
   \begin{equation}
   \sup_{s\in [0,R_{1}]}|\xi_{s}|\le\max_{0\le k\le R_{1}-1}|S_{k}|+\max_{0\le k\le R_{1}-1}\sup_{s\in [k,k+1]}|\xi_{s}-S_{k}|.\nonumber
   \end{equation}
   Applying similar calculation as in \eqref{propA.2.0} we can deduce that
   \begin{eqnarray}
   \bP_{0,\rho}\left[\max_{0\le k\le R_{1}-1}\sup_{s\in [k,k+1]}|\xi_{s}-S_{k}|\right]
   &\le& \bP_{0,\rho}\left[\sum_{k=0}^{R_{1}-1}\sup_{s\in [k,k+1]}|\xi_{s}-S_{k}|\right]\nonumber\\
   &=&\bP_{0,\rho}[R_{1}]\bP_{0,\pi}\left[\sup_{s\in [0,1]}|\xi_{s}|\right]<+\infty.\nonumber
   \end{eqnarray}
   Hence \eqref{propA.2.2} follows from this and \eqref{propA.2.0}, completing the proof.\qed

\smallskip

   \begin{proposition}\label{prop:tail}
   Suppose the conditions of Proposition \ref{thm:slln} hold. Then we have (a') $\xi_{t}\to +\infty$, (b') $\limsup_{t\to +\infty}\xi_{t}=+\infty$, $\liminf_{t\to+\infty}\xi_{t}=-\infty$ and (c') $\xi_{t}\to -\infty$ $\bP_{0,\theta}$-a.s. for  every $\theta\in\mathcal{S}$ according as (a) $\bP_{0,\pi}[\xi_{1}]>0$, (b) $\bP_{0,\pi}[\xi_{1}]=0$ and the increment distribution in each block is not concentrated at $0$ and (c) $\bP_{0,\pi}[\xi_{1}]<0$.
   \end{proposition}
   \proof It is immediate from Proposition \ref{thm:slln} that (a)$\Rightarrow$(a') and (c)$\Rightarrow$(c'). In case (b), we consider the sequence $\{S_{R_{k}}:k\ge 0\}$ which is a discrete-time random walk with mean $0$ and the increment distribution not concentrated at $0$. Hence
   $\limsup_{k\to+\infty}S_{R_{k}}=+\infty$ and $\liminf_{k\to+\infty}S_{R_{k}}=-\infty$ which implies (b').\qed

\begin{remark}\label{HTremark}
\rm
Let us make a brief remark on the condition \eqref{condi:Harris type}. This condition is of course not the most general condition under which the results of Propositions \ref{thm:slln} and \ref{prop:tail} hold.
We believe an extension is possible, at least to some extent. One direction is to assume Harris recurrence of $(M_{n})_{n\ge 0}$ alone. However, in this way, instead of having i.i.d increments, $\{S_{R_{n}}:n\ge 0\}$ has $1$-dependent and stationary increments. Therefore in all places where we apply results for ordinary random walks, extensions to the case of $1$-dependent and stationary increments are needed. Since this can not be done shortly, we have restricted this section to the case when condition \eqref{condi:Harris type} is satisfied.
\end{remark}

Hereafter we say that $\xi_{t}$ \textit{drifts to $+\infty$}, \textit{oscillates}, or \textit{drifts to $-\infty$ at $\theta$}, respectively, if $\lim_{t\to+\infty}\xi_{t}=+\infty$, $\limsup_{t\to +\infty}\xi_{t}=-\liminf_{t\to+\infty}\xi_{t}=+\infty$, or $\lim_{t\to +\infty}\xi_{t}=-\infty$ $\bP_{0,\theta}$-a.s.

  \begin{proposition}\label{prop:2:12}
For every $\theta\in\s$,
  \begin{equation} \nonumber
\int_{\s\times \R^{+} }\emph{\rm n}^{+}_{v}(\zeta=+\infty)U^{+}_{\theta}({\rm d}v,{\rm d}z)=\begin{cases}
         0
         \quad &\hbox{if $\xi_{t}$ oscillates or drifts to $+\infty$ at $\theta$,} \smallskip \\
         1
         \quad &\hbox{if $\xi_{t}$ drifts to $-\infty$ at $\theta$.}
        \end{cases}
\end{equation}

  \end{proposition}

 \proof Let $\bar{g}_{\infty}$ denote the last time when $\xi_{t}$ attains its running maximum. If $\xi_{t}$ oscillates or drifts to $+\infty$ at $\theta$, then $\bP_{0,\theta}\left(\bar{g}_{\infty}=+\infty\right)=1$.
 By Proposition \ref{prop:wiener_hopf} we have
 \begin{equation}\label{lemA.5.1}
 \bP_{0,\theta}\left[\mathrm{e}^{-\lambda \bar{g}_{\mathbf{e}_{q}}}\right]=\int_{\R^{+}\times \s\times \R^{+} }\mathrm{e}^{-\lambda r-qr}\left(q\ell^+(v)+{\rm n}^{+}_{v}\left(1-\mathrm{e}^{-q\zeta}\right)\right)V^{+}_{\theta}({\rm d}r,{\rm d}v,{\rm d}z)\quad\forall \lambda,q>0.
 \end{equation}
 Letting $q\to 0+$, we get by Fatou's lemma that
 \begin{equation}\nonumber
 0=\bP_{0,\theta}\left[\mathrm{e}^{-\lambda \bar{g}_{\infty}}\right]\ge \int_{\R^{+}\times \s\times \R^{+} }\mathrm{e}^{-\lambda r}{\rm n}^{+}_{v}(\zeta=+\infty)V^{+}_{\theta}({\rm d}r,{\rm d}v,{\rm d}z).
 \end{equation}
Then by letting $\lambda\to 0+$, we get by the monotone convergence theorem that
 $$\int_{\s\times \R^{+} }{\rm n}^{+}_{v}(\zeta=+\infty)U^{+}_{\theta}({\rm d}v,{\rm d}z)=0.$$
On the other hand, if $\xi_{t}$ drifts to $-\infty$ at $\theta$, then $\bP_{0,\theta}\left(\bar{g}_{\infty}<+\infty\right)=1$.
 Note that for any $0<q<\lambda/2$, the integrand in the right-hand side of \eqref{lemA.5.1} is bounded from above by $\mathrm{e}^{-\lambda r}\left(\frac{\lambda}{2}\ell^+(v)+{\rm n}^{+}_{v}\left(1-\mathrm{e}^{-\lambda \zeta/2}\right)\right)$ and
 \begin{equation}
 \int_{\R^{+}\times \s\times \R^{+} }\mathrm{e}^{-\lambda r}\left(\frac{\lambda}{2}\ell^+(v)+{\rm n}^{+}_{v}\left(1-\mathrm{e}^{-\lambda \zeta/2}\right)\right)V^{+}_{\theta}({\rm d}r,{\rm d}v,{\rm d}z)=\bP_{0,\theta}\left[\mathrm{e}^{-\frac{\lambda}{2} \bar{g}_{e_{\lambda/2}}}\right]<+\infty.\nonumber
 \end{equation}
 Hence by letting $q\to 0+$ in \eqref{lemA.5.1} and using the dominated convergence theorem in the right-hand side and the monotone convergence theorem in the left hand side we get
 $$\bP_{0,\theta}\left[\mathrm{e}^{-\lambda \bar{g}_{\infty}}\right]=\int_{\R^{+}\times \s\times \R^{+} }\mathrm{e}^{-\lambda r}{\rm n}^{+}_{v}(\zeta=+\infty)V^{+}_{\theta}({\rm d}r,{\rm d}v,{\rm d}z).$$
Letting $\lambda\to 0+$, we have
$$\int_{\s\times \R^{+} }{\rm n}^{+}_{v}(\zeta=+\infty)U^{+}_{\theta}({\rm d}v,{\rm d}z)=\bP_{0,\theta}\left(\bar{g}_{\infty}<+\infty\right)=1,$$
which completes the proof.\qed

\subsection{Invariant measures}

\begin{proposition}\label{prop:existence of invariant measure}
Suppose $((\xi,\Theta),\bP)$ is a MAP on $\R\times\s$ and $\nu$ is an invariant measure for the modulator $\Theta$. Then the measure
\begin{equation}\label{lem0.1.0}
\nu^{+}(\cdot):=\bP_{0,\nu}\left[\int_{0}^{1}\1_{\{\Theta_{s}\in\cdot\}}{\rm d}\bar{L}_{s}\right]
\end{equation}
is an invariant measure for the modulator $\Theta^{+}$ of the ascending ladder height process $((\xi^{+},\Theta^{+}),\bP)$. Moreover,
$\nu^{+}$ is finite if and only if $\bP_{0,\nu}\left[\bar{L}_{1}\right]<+\infty$.
\end{proposition}

\proof It suffices to show that
\begin{equation}
\int_{0}^{+\infty}\mathrm{e}^{-\alpha s}\bP_{0,\nu^{+}}\left[f(\Theta^{+}_{s})\right]{\rm d} s=\frac{1}{\alpha}\int_{\s}f(\theta)\nu^{+}(d\theta)\label{lem0.1.1}
\end{equation}
for any $\alpha>0$ and nonnegative measurable function $f:\s\to\R^{+}$. The left integral is equal to
\begin{eqnarray}
\bP_{0,\nu^{+}}\left[\int_{0}^{+\infty}\mathrm{e}^{-\alpha s}f(\Theta^{+}_{s}){\rm d} s\right]
&=&\bP_{0,\nu^{+}}\left[\int_{0}^{+\infty}\mathrm{e}^{-\alpha \bar{L}_{s}}f(\Theta_{s}){\rm d} \bar{L}_{s}\right]\nonumber\\
&=&\bP_{0,\nu}\left[\int_{0}^{1}\bP_{0,\Theta_{r}}\left[\int_{0}^{+\infty}\mathrm{e}^{-\alpha \bar{L}_{s}}f(\Theta_{s}){\rm d} \bar{L}_{s}\right]{\rm d} \bar{L}_{r} \right].\label{lem0.1.2}
\end{eqnarray}
Recall that $s\mapsto\bar{L}_{s}$ is an additive functional of $(\Theta_{t},\bar{\xi}_{t}-\xi_{t})_{t\ge 0}$. Hence the law of $(\bar{L}_{t},\Theta_{t})_{t\ge 0}$ under $\bP_{x,\theta}$ does not depend on $x$. The right hand side of \eqref{lem0.1.2} is equal to
\begin{eqnarray}
&&\bP_{0,\nu}\left[\int_{0}^{1}\bP_{\xi_{r},\Theta_{r}}\left[\int_{0}^{+\infty}\mathrm{e}^{-\alpha \bar{L}_{s}}f(\Theta_{s}){\rm d} \bar{L}_{s}\right]{\rm d} \bar{L}_{r} \right]\nonumber\\
&=&\bP_{0,\nu}\left[\int_{0}^{1}{\rm d} \bar{L}_{r}\int_{r}^{+\infty}\mathrm{e}^{-\alpha\left(\bar{L}_{s}-\bar{L}_{r}\right)}f(\Theta_{s}){\rm d} \bar{L}_{s}\right]\nonumber\\
&=&\bP_{0,\nu}\left[\int_{0}^{+\infty}{\rm d} \bar{L}_{s}\mathrm{e}^{-\alpha \bar{L}_{s}}f(\Theta_{s})\int_{0}^{1\wedge s}\mathrm{e}^{\alpha\bar{L}_{r}}{\rm d} \bar{L}_{r}\right]\nonumber\\
&=&\frac{1}{\alpha}\left[\bP_{0,\nu}\left[\int_{0}^{1}\mathrm{e}^{-\alpha\bar{L}_{s}}f(\Theta_{s})\left(\mathrm{e}^{\alpha\bar{L}_{s}}-1\right){\rm d} \bar{L}_{s}\right]+\bP_{0,\nu}\left[\int_{1}^{+\infty}\mathrm{e}^{-\alpha\bar{L}_{s}}f(\Theta_{s})\left(\mathrm{e}^{\alpha\bar{L}_{1}}-1\right){\rm d} \bar{L}_{s}\right]\right]\nonumber\\
&=&\frac{1}{\alpha}\left[\bP_{0,\nu}\left[\int_{0}^{1}f(\Theta_{s}){\rm d} \bar{L}_{s}\right]
-\bP_{0,\nu}\left[\int_{0}^{+\infty}\mathrm{e}^{-\alpha \bar{L}_{s}}f(\Theta_{s}){\rm d} \bar{L}_{s}\right]\right.\nonumber\\
&&\quad\quad \left.+\bP_{0,\nu}\left[\int_{1}^{+\infty}\mathrm{e}^{-\alpha\left(\bar{L}_{s}-\bar{L}_{1}\right)}f(\Theta_{s}){\rm d} \bar{L}_{s}\right]\right].\label{lem0.1.3}
\end{eqnarray}
In the first equality we use the Markov property and the additivity of $\bar{L}_{s}$. Using these facts again we have
\begin{eqnarray}
\bP_{0,\nu}\left[\int_{1}^{+\infty}\mathrm{e}^{-\alpha\left(\bar{L}_{s}-\bar{L}_{1}\right)}f(\Theta_{s}){\rm d} \bar{L}_{s}\right]
&=&\bP_{0,\nu}\left[\bP_{\xi_{1},\Theta_{1}}\left[\int_{0}^{+\infty}\mathrm{e}^{-\alpha\bar{L}_{r}}f(\Theta_{r}){\rm d} \bar{L}_{r}\right]\right]\nonumber\\
&=&\bP_{0,\nu}\left[\bP_{0,\Theta_{1}}\left[\int_{0}^{+\infty}\mathrm{e}^{-\alpha\bar{L}_{r}}f(\Theta_{r}){\rm d} \bar{L}_{r}\right]\right]\nonumber\\
&=&\bP_{0,\nu}\left[\int_{0}^{+\infty}\mathrm{e}^{-\alpha\bar{L}_{r}}f(\Theta_{r}){\rm d} \bar{L}_{r}\right].\label{lem0.1.4}
\end{eqnarray}
In the last equality we use the fact that $\bP_{0,\nu}\left(\Theta_{1}\in\cdot\right)=\nu(\cdot)$. In view of \eqref{lem0.1.4}, the right hand side of \eqref{lem0.1.3} equals
$$\frac{1}{\alpha}\bP_{0,\nu}\left[\int_{0}^{1}f(\Theta_{s}){\rm d} \bar{L}_{s}\right]=\frac{1}{\alpha}\int_{\s}f(\theta)\nu^{+}(d\theta).$$
Hence we get \eqref{lem0.1.1}.\qed

\smallskip

\begin{corollary}\label{cor:existence of invariant distribution}
Suppose the modulator $\Theta$ of $((\xi,\Theta),\bP)$ has an invariant distribution $\pi$. If  $\bP_{0,\pi}\left[\bar{L}_{1}\right]>0$ and $\inf_{\theta\in\s}\left[\ell^{+}(\theta)+{\rm n}^{+}_{\theta}\left(1-\mathrm{e}^{-\zeta}\right)\right]>0$, then the measure $\pi^{+}$ defined by
$$\pi^{+}(\cdot):=\frac{1}{\bP_{0,\pi}\left[\bar{L}_{1}\right]}\bP_{0,\pi}\left[\int_{0}^{1}\1_{\{\Theta_{s}\in\cdot\}}{\rm d}\bar{L}_s\right]$$
is an invariant distribution for the modulator $\Theta^{+}$ of $((\xi^{+},\Theta^{+}),\bP)$.
\end{corollary}

\proof By Proposition \ref{prop:existence of invariant measure}, it suffices to show that $\bP_{0,\pi}\left[\bar{L}_{1}\right]<+\infty$.
Let $$c:=\inf_{\theta\in\s}\left[\ell^{+}(\theta)+{\rm n}^{+}_{\theta}\left(1-\mathrm{e}^{-\zeta}\right)\right]\in (0,+\infty).$$ By \eqref{M2} and \eqref{eq:last exit} we have
\begin{eqnarray}
\bP_{0,\pi}\left[\bar{L}_{1}\right]&\le&\frac{1}{c}\bP_{0,\pi}\left[\int_{0}^{1}\ell^{+}(\Theta_{s})+{\rm n}^{+}_{\Theta_{s}}\left(1-\mathrm{e}^{-\zeta}\right){\rm d}\bar{L}_{s}\right]\nonumber\\
&=&\frac{1}{c}\left[\bP_{0,\pi}\left[\int_{0}^{1}\1_{\{s\in\bar{M}\}}{\rm d}s\right]+\bP_{0,\pi}\left[\sum_{g_{i}\in\bar{G},g_{i}\le 1}\left(1-\mathrm{e}^{-\zeta^{(g_{i})}}\right)\right]\right]\nonumber\\
&\le&\frac{1}{c}\left[1+\bP_{0,\pi}\left[\sum_{g_{i}\in\bar{G},g_{i}\le 1}\left(1\wedge \zeta^{(g_{i})}\right)\right]\right].\nonumber
\end{eqnarray}
 We note that among all the excursions that start in the time interval $[0,1]$, there is, at most, one excursion having a lifetime longer than $1$, and the sum of lifetimes of other excursions does not exceed $1$. Hence $\bP_{0,\pi}\left[\sum_{g_{i}\in\bar{G},g_{i}\le 1}\left(1\wedge \zeta^{(g_{i})}\right)\right]\le 2$ and $\bP_{0,\pi}\left[\bar{L}_{1}\right]<+\infty$.\qed

\section{Duality}\label{sec:duality}

    In this section we present the notion of duality as well as several results about duality. Here we suppose that $E$ is a Polish space and $\mu$ is a $\sigma$-finite Radon measure on $E$. Suppose that $(X,\P)$ and $(\hat{X},\mathbb Q)$ are two, possibly killed, right continuous strong Markov processes having left limits in $E$ except perhaps at their lifetime. We use $\zeta$ and $\hat{\zeta}$ respectively to denote their lifetimes. We take the convention that
    {$X_{0-}=X_{0}$ and $\hat{X}_{0-}=\hat{X}_{0}$.}

    \begin{definition}
      {Two processes}
      $(X,\P)$ and $(\hat{X},\mathbb Q)$ are dual with respect to $\mu$ if for every
      {nonnegative measurable functions $f,g:E\rightarrow \R^{+}$} and every $t \geq 0$,
      \[
        \int_E \mu({\rm d}x) g(x)\P_x[f(X_t),t<\zeta]=\int_E \mu({\rm d}x) f(x)\mathbb Q_x[g(\hat{X}_t),t<\hat{\zeta}].
      \]
    \end{definition}

    Note that there is no requirement that $\mu$ is a finite measure. The notion of duality is closely linked with reversibility. The following result is from~\cite[Theorem 2.1]{Walsh}.

    \begin{lemma}\label{lem:walsh_fixed_time}
      Suppose that $(X,\P)$ and $(\hat{X},\mathbb Q)$ are dual with respect to $\mu$, then,
      \[
        \int_E\mu({\rm d}x) \P_x\left[F\left((X_s)_{s \leq t}\right)\1_{\{t < \zeta\}}\right] = \int_E\mu({\rm d}x) \mathbb Q_x\left[F\left((\hat{X}_{(t-s)-})_{s \leq t}\right)\1_{\{t < \zeta\}}\right]
      \]
    for every $t\ge 0$ and nonnegative functional $F:\mathbb{D}_{E}[0,t]\to \R^{+}$.
    \end{lemma}

    {We present a result on the time reversal from the lifetime, which is an application of \cite[Theorem 3.5]{Nagasawa}. We also refer to \cite[Theorem 13.34]{ChungWalsh} for a simple proof in the special case where the resolvent kernels of $(X,\P)$ and $(\hat{X},\mathbb{Q})$ are absolutely continuous with respect to $\mu$.}

    \begin{lemma}\label{lem:walsh_killing_time}
      Suppose that $(X,\P)$ and $(\hat{X},\mathbb Q)$ are dual with respect to $\mu$. If the process $X$ has initial distribution
      $\eta$ and a finite lifetime $\zeta$ such that
      \begin{equation}
      \int_E \mu({\rm d}x)f(x)=\int_E {\eta}({\rm d}x)\P_x\left[\int_0^\zeta f(X_t){\rm d}t\right]\label{condi:Nagasawa}
      \end{equation}
      for every nonnegative measurable function $f:E\rightarrow\R^{+}$, then
      $\left((X_{(\zeta-t)-})_{0<t<\zeta},\P_{\eta}\right)$ is a right continuous strong Markov process having the same transition rates as $(\hat{X},\mathbb{Q})$.
    \end{lemma}

 We remark here that in general the measure $\eta$ appearing in \eqref{condi:Nagasawa} may not exist.
  If exists, it is uniquely determined by the reference measure $\mu$, see, for example, \cite[Theorem 2.12 and Section 6]{Getoor}.

  \smallskip

  Throughout the remainder of this paper, we assume that the process $((\xi,\Theta),\widetilde{\bP})
  $ is a MAP satisfying that
  $\widetilde{\bP}_{y,v}\left(\xi_{0}=y,\Theta_{0}=v\right)=1$, and is linked to $((\xi,\Theta),\bP)$ through the following weak reversability property: There exists a probability measure $\pi$ on $\s$ with full support such that
  \begin{equation}\label{condi:weak reversability}\tag{WR}
    \bP_{0,\theta}(\xi_t \in {\rm d}z; \Theta_t \in {\rm d}v)\pi({\rm d}\theta) = \widetilde{\bP}_{0,v}(\xi_t \in {\rm d}z; \Theta_t \in {\rm d}\theta)\pi({\rm d}v) \quad \forall t \geq 0.
  \end{equation}
 By integrating \eqref{condi:weak reversability} over variable $z$, it follows that the Markov processes $((\Theta_{t})_{t\ge 0},\{\bP_{0,\theta},\theta\in\s\})$ and  $((\Theta_{t})_{t\ge 0},\{\widetilde{\bP}_{0,\theta},\theta\in\s\})$ are dual with respect to the measure $\pi$.
Hereafter we denote by $\hat{\bP}_{x,\theta}$ the law of $(-\xi,\Theta)$ under $\widetilde{\bP}_{-x,\theta}$.
We will use the notation ~$\hat{}$~ to specify the mathematical quantities related to the process $((\xi,\Theta),\hat{\bP})$.
In the following we give some examples for a MAP to be weakly reversible. Each example corresponds to a well-known class of ssMps via Lamperti-Kiu transform.

 \begin{example}\label{eg1}
 Suppose $\s=\{s_{1},\cdots,s_{n}\}$ is a finite set.
 It is known that the process $((\xi,
 \Theta),\bP)$ is a MAP on $\R\times \s$ if and only if $((\Theta_{t})_{t\ge 0},\{\bP_{x,\theta}:\theta\in\s\})$ is a (possibly killed) Markov chain on $\s$ whose law does not depend on $x$, and for each $s_{j},s_{k}\in \s$ there exist a (non-killed) L\'{e}vy process $\xi^{j}$ and an $\R$-valued random variable $\Xi_{j,k}$ such that when $\Theta$ is in state $s_{j}$, $\xi$ evolves according to an independent copy of $\xi^{j}$, and when $\Theta$ changes from $s_{j}$ to another state $s_{k}$, $\xi$ has an additional jump which is an independent copy of $\Xi_{j,k}$ and until the next jump of $\Theta$, $\xi$ evolves according to an independent copy of $\xi^{k}$, and so on, until the lifetime of $\Theta$.
  For such a MAP condition \eqref{condi:weak reversability} is equivalent to require that there is a MAP $((\xi,\Theta),\widetilde{\bP})$ on $\R\times\s$ and a probability measure $\pi$ on $\s$ such that $\pi_{j}=\pi(\{s_{j}\})>0$ for $1\le j\le n$ and
  \begin{equation}\label{eg1.1}
  \pi_{j}\widetilde{\bP}_{0,s_{j}}\left[\rme^{i\lambda \xi_{t}}\1_{\{\Theta_{t}=s_{k}\}}\right]=\pi_{k}\bP_{0,s_{k}}\left[\rme^{i\lambda \xi_{t}}\1_{\{\Theta_{t}=s_{j}\}}\right]\quad \forall t\ge 0,\ \lambda\in\R,\ 1\le j,k\le n.
  \end{equation}
We let $(q_{j,k})_{1\le j,k\le n}$ denote the intensity matrix of the Markov chain $\Theta$, $\psi_{j}(\lambda)$ denote the characteristic exponent of the L\'{e}vy process $\xi^{j}$ and $J_{j,k}(\lambda)$ denote the characteristic function of the random variable $\Xi_{j,k}$. The matrix
 $${\boldsymbol F}(\lambda):=\mathrm{diag}(-\psi_{1}(\lambda),\cdots,-\psi_{n}(\lambda))+\left(q_{j,k}J_{j,k}(\lambda)\right)_{1\le j,k\le n}\quad \forall \lambda\in\R$$
 is called the characteristic matrix exponent of the MAP $((\xi,\Theta),\bP)$ because
 $$\bP_{0,s_{j}}\left[\rme^{i\lambda  \xi_{t}}\1_{\{\Theta_{t}=s_{k}\}}\right]=\left(\rme^{\boldsymbol F(\lambda)t}\right)_{j,k}\quad\forall t\ge 0,\ 1\le j,k\le n.$$
 Equation \eqref{eg1.1}, in terms of the characteristic matrix exponent, is equivalent to
 $$\widetilde{\boldsymbol F}(\lambda)=\boldsymbol {\triangle}^{-1}_{\pi}{\boldsymbol F}(\lambda)^{\mathrm{T}}\boldsymbol{\triangle}_{\pi}\quad\forall \lambda\in\R,$$
 where $\boldsymbol {\triangle}_{\pi}=\mathrm{diag}(\pi_{1},\cdots,\pi_{n})$.
 Condition \eqref{condi:weak reversability} is satisfied, in particular, if the process $\Theta$ is dual with itself with respect to a probability measure $\pi$ and
 $\Xi_{j,k}\stackrel{d}{=}\Xi_{k,j}$ for all $1\le j,k\le n$, in which case we can take $\widetilde{\bP}=\bP$.
 \end{example}

\begin{example}\label{eg2}
Suppose $\partial$ is an isolated extra state and the transition probabilities of $\left((\xi,\Theta),\bP\right)$ have the following form:
\begin{equation} \nonumber
\left\{ \begin{aligned}
         &\bP_{x,\theta}\left(\xi_{t}\in {\rm d}y,\ \Theta_{t}\in {\rm d}v\right)=\mathrm{e}^{-\lambda t}\mathrm{P}^{\xi'}_{x}\left(\xi'_{t}\in {\rm d}y\right)\mathrm{P}^{\Theta'}_{\theta}\left(\Theta'_{t}\in {\rm d}v\right),  \\
  &\bP_{x,\theta}\left((\xi_{t},\Theta_{t})=\partial\right)=1-\mathrm{e}^{-\lambda t}
                          \end{aligned} \right.
\end{equation}
for all $t\ge 0$ and $(x,\theta)\in\R\times \s$, where $\lambda\ge 0$ is a constant, $(\xi',\mathrm{P}^{\xi'}_{x})$ is a non-killed $\R$-valued L\'{e}vy process started from $x$ and $(\Theta',\mathrm{P}^{\Theta'}_{\theta})$ is a non-killed $\s$-valued Markov process started from $\theta$. Then condition \eqref{condi:weak reversability} is satisfied if and only if there exists an $\s$-valued Markov process $((\Theta'_{t})_{t\ge 0},\{\widetilde{\mathrm{P}}^{\Theta'}_{\theta},\theta\in\s\})$, which is dual to $((\Theta'_{t})_{t\ge 0},\{\mathrm{P}^{\Theta'}_{\theta},\theta\in\s\})$ with respect to a probability measure $\pi$ on $\s$. In this case, we can take the MAP $((\xi,\Theta),\widetilde{\bP})$ to be such that its transition probabilities have the following form:
\begin{equation} \nonumber
\left\{ \begin{aligned}
         &\widetilde{\bP}_{x,\theta}\left(\xi_{t}\in {\rm d}y,\ \Theta_{t}\in {\rm d}v\right)=\mathrm{e}^{-\lambda t}\mathrm{P}^{\xi'}_{x}\left(\xi'_{t}\in {\rm d}y\right)\widetilde{\mathrm{P}}^{\Theta'}_{\theta}\left(\Theta'_{t}\in {\rm d}v\right),  \\
  &\widetilde{\bP}_{x,\theta}\left((\xi_{t},\Theta_{t})=\partial\right)=1-\mathrm{e}^{-\lambda t}
                          \end{aligned} \right.
\end{equation}
for all $t\ge 0$ and $(x,\theta)\in \R\times\s$.
\end{example}

\begin{example}\label{eg3}
Suppose $\s=\mathbb{S}^{d-1}$ and for any orthogonal transformation $\mathcal{O}$ of $\mathbb{S}^{d-1}$ and $(x,\theta)\in \R\times\mathbb{S}^{d-1}$, $\left((\xi,\Theta),\bP_{x,\theta}\right)$ is equal in law with $\left((\xi,\mathcal{O}(\Theta)),\bP_{x,\mathcal{O}^{-1}(\theta)}\right)$. In view of this property, if $X$ is the ssMp associated with $(\xi,\Theta)$ by Lamperti-Kiu transform, then $X$ is a rotationally invariant Markov process on $\R^{d}$. Hence its norm $(\|X_{t}\|)_{t\ge 0}$ is a positive ssMp, which in turn implies that $\xi$ alone is a L\'{e}vy process. In this case condition \eqref{condi:weak reversability} is satisfied with $\widetilde{\bP}=\bP$ and $\pi$ being the uniform measure on the sphere $\mathbb{S}^{d-1}$. We refer to \cite[Proposition 3.2]{ACGZ} for a proof.
\end{example}

\begin{proposition}\label{lem:easy switch}
     The processes $((\xi,\Theta),\bP)$ and $((\xi,\Theta),\hat\bP)$ are dual with respect to the measure ${\rm Leb} \otimes \pi$, where ${\rm Leb}$ is the Lebesgue measure on $\R$.
  \end{proposition}

\proof Suppose $f,g:\R\times \s\to \R^{+}$ are nonnegative measurable functions. By an application of Fubini's theorem, a change of variable and condition \eqref{condi:weak reversability} we get
\begin{align}
\int_{\R\times\s}&\mathrm{d}x\pi({\rm d} \theta)f(x,\theta)\bP_{x,\theta}\left[g(\xi_{t},\Theta_{t})\right]\nonumber\\
&=\int_{\R\times\s}\mathrm{d}x\pi({\rm d} \theta)f(x,\theta)\bP_{0,\theta}\left[g(x+\xi_{t},\Theta_{t})\right]\nonumber\\
&=\int_{\R\times\s}\mathrm{d}y\pi({\rm d} \theta)\bP_{0,\theta}\left[f(y-\xi_{t},\theta)g(y,\Theta_{t})\right]\nonumber\\
&=\int_{\R\times\s}\mathrm{d}y\pi({\rm d} \theta)\int_{\R\times\s}\bP_{0,\theta}\left(\xi_{t}\in {\rm d} z,\Theta_{t}\in{\rm d} \nu\right)f(y-z,\theta)g(y,\nu)\nonumber\\
&=\int_{\R\times \s}\mathrm{d}y\pi({\rm d} \nu)\int_{\R\times\s}\widetilde{\bP}_{0,\nu}\left(\xi_{t}\in {\rm d} z,\Theta_{t}\in{\rm d} \theta\right)f(y-z,\theta)g(y,\nu)\nonumber\\
&=\int_{\R\times \s}{\rm d} y \pi({\rm d}\nu)g(y,\nu)\widetilde{\bP}_{0,\nu}\left[f(y-\xi_{t},\Theta_{t})\right]\nonumber\\
&=\int_{\R\times \s}{\rm d} y \pi({\rm d}\nu)g(y,\nu)\hat{\bP}_{0,\nu}\left[f(y+\xi_{t},\Theta_{t},)\right]\nonumber\\
&=\int_{\R\times \s}{\rm d} y \pi({\rm d}\nu)g(y,\nu)\hat{\bP}_{y,\nu}\left[f(\xi_{t},\Theta_{t})\right]\nonumber
\end{align}
for all $t\ge 0$. Hence we complete the proof.\qed

\smallskip

\begin{lemma}\label{lem:time reversal}
Suppose $t>0$. For every $x\in \R,$ the process $(\xi_{(t-s)-}-\xi_{t},\Theta_{(t-s)-})_{ 0\leq s\leq t }$ under $\bP_{x,\pi}$ has the same law as $(\xi_{s},\Theta_{s})_{0\leq s\leq t}$ under $\hat{\bP}_{0,\pi}$.
\end{lemma}
\proof
In order to prove this lemma it suffices to consider the finite dimensional distributions. Let $n\ge 1$ be a fixed integer. For $0\leq k\leq n$ we take nonnegative measurable functions $f_{k}:\s \times \R\to \R^{+}$ and $0=t_{0}<t_{1}<t_{2}<\cdots<t_{n}<t_{n+1}=t.$ Let $g:\R\to\R^{+}$ be a nonnegative measurable function. We know by Proposition \ref{lem:easy switch} and Lemma \ref{lem:walsh_fixed_time} that the process $\left((\xi_{(t-s)-},\Theta_{(t-s)-})_{0\le s\le t }, \bP\right)$ has the same law as $\left((\xi_{s},\Theta_{s})_{0\le s\le t }, \hat{\bP}\right)$ both started according to the measure $\mathrm{Leb}\otimes\pi$. Using this and the quasi-left continuity of $\xi$, we have

\begin{align}
\int_{\R\times\s}&{\rm d}x\pi({\rm d}\theta) g(x)\bP_{x,\theta}\left[f_{0}(\Theta_{(t-t_{0})-},\xi_{(t-t_{0})-}-\xi_{t})
\cdots f_{n}(\Theta_{(t-t_{n})-},\xi_{(t-t_{n})-}-\xi_{t})\right]\nonumber\\
&=\int_{\R\times\s}{\rm d}x\pi({\rm d}\theta) \bP_{x,\theta}\left[f_{0}(\Theta_{(t-t_{0})-},\xi_{(t-t_{0})-}-\xi_{t-})
\cdots f_{n}(\Theta_{(t-t_{n})-},\xi_{(t-t_{n})-}-\xi_{t-})g(\xi_{(t-t_{n+1})-})\right]\nonumber\\
&=\int_{\R\times\s}{\rm d}x\pi({\rm d}\theta) \hat{\bP}_{x,\theta}\left[f_{0}(\Theta_{t_{0}},\xi_{t_{0}}-\xi_{0})
\cdots f_{n}(\Theta_{t_{n}},\xi_{t_{n}}-\xi_{0})g(\xi_{t_{n+1}})\right]\nonumber\\
&=\int_{\R\times\s}{\rm d}x\pi({\rm d}\theta) \hat{\bP}_{0,\theta}\left[f_{0}(\Theta_{t_{0}},\xi_{t_{0}})
\cdots f_{n}(\Theta_{t_{n}},\xi_{t_{n}})g(\xi_{t_{n+1}}+x)\right].\nonumber
\end{align}
By Fubini's theorem and a change of variable, the integral on the right hand side is equal to
\begin{align}
&\hat{\bP}_{0,\pi}\left[f_{0}(\Theta_{t_{0}},\xi_{t_{0}})
\cdots f_{n}(\Theta_{t_{n}},\xi_{t_{n}})\int_{\R}g(\xi_{t_{n+1}}+x){\rm d}x\right]\nonumber\\
&=\int_{\R\times\s}{\rm d}x\pi({\rm d}\theta) g(x)\hat{\bP}_{0,\theta}\left[f_{0}(\Theta_{t_{0}},\xi_{t_{0}})
\cdots f_{n}(\Theta_{t_{n}},\xi_{t_{n}})\right],\nonumber
\end{align}
Since $g$ is arbitrary, it follows by above equations that $\{(\xi_{s},\Theta_{s}), 0\leq s\leq t \}$ under $\hat{\bP}_{0,\pi}$ has the same law as $\{(\xi_{(t-s)-}-\xi_{t},\Theta_{(t-s)-}), 0\leq s\leq t \}$ under ${\bP}_{x,\pi}$ for almost every $x\in\R$. We observe that the law of the latter does not depend on $x$, thus the claim holds for every $x\in\R$.\qed

\smallskip

The upwards regularity of $((\xi,\Theta),\bP)$ implies that almost surely the local maxima of $\xi$ during a finite time interval are all distinct. In view of this and Lemma \ref{lem:time reversal}, we have the following result.

\begin{proposition}\label{prop:equal dist}
For every $t>0$, $\left(\Theta_{0},t-\bar{g}_{t},\Theta_{t},\bar{\xi}_{t}-\xi_{t},\bar{g}_{t},\bar{\Theta}_{t},\bar{\xi}_{t}\right)$ under $\hat{\bP}_{0,\pi}$
is equal in distribution to $\left(\Theta_{t},\bar{g}_{t},\Theta_{0},\bar{\xi}_{t},t-\bar{g}_{t},\bar{\Theta}_{t},\bar{\xi}_{t}-\xi_{t}\right)$ under $\bP_{0,\pi}$.
\end{proposition}

\section{MAPs conditioned to stay negative}\label{sec:condi-negative}

In this section we assume that
$((\xi,\Theta),\hat{\bP})$ is an upwards regular MAP.
Define
$$\hat{H}^{+}_{\theta}(y):=\hat{\bP}_{y,\theta}\left(\tau^{+}_{0}=+\infty\right),\quad\forall (y,\theta)\in\R\times \s .$$
Obviously, $\hat{H}^{+}_{\theta}(y)=0$ for all $y\ge 0$.

\begin{proposition}\label{prop:martingale}
Assume that
\begin{equation}
\label{prea7}
\hat{{\rm n}}^{+}_{v}(\zeta=+\infty)>0\text{  for every }v\in\s,
\end{equation}
then
\begin{itemize}
\item[(i)] $\hat{H}^{+}_{\theta}(y)>0$ for all $\theta\in\s $ and $y<0$;

\item[(ii)] $\hat{H}^{+}_{\Theta_{t}}(\xi_{t})\1_{\{t<\tau^{+}_{0}\}}$ is a $\hat{\bP}_{y,\theta}$-martingale for every $y<0$ and $\theta\in\s $.
\end{itemize}
\end{proposition}

\proof (i) For $y<0$ and $\theta\in\s$,
$$\hat{H}^{+}_{\theta}(y)=\hat{\bP}_{0,\theta}\left(\tau^{+}_{-y}=+\infty\right)=\lim_{q\to 0+}\hat{\bP}_{0,\theta}\left(\tau^{+}_{-y}>\mathbf{e}_{q}\right).$$
It follows by Proposition \ref{prop:wiener_hopf} that
\begin{eqnarray}
\hat{\bP}_{0,\theta}\left(\tau^{+}_{-y}>\mathbf{e}_{q}\right)
&=&\hat{\bP}_{0,\theta}\left(\bar{\xi}_{\mathbf{e}_{q}}\le -y\right)\nonumber\\
&=&\int_{\R^{+}\times\s\times \R^{+}}\rme^{-q r}\1_{\{z\le -y\}}\left(q \hat\ell^+(v)+\hat{{\rm n}}^{+}_{v}\left(1-\rme^{-q \zeta}\right)\right)\hat{V}^{+}_{\theta}({\rm d} r,{\rm d} v,{\rm d} z).\nonumber
\end{eqnarray}
Hence by condition \eqref{prea7} and Fatou's lemma,
$$\hat{H}^{+}_{\theta}(y)\ge \int_{\s\times [0,-y]}\hat{{\rm n}}^{+}_{v}(\zeta=+\infty)\hat{U}^{+}_{\theta}({\rm d} v,{\rm d} z)>0.$$

(ii) By the Markov property of $((\xi,\Theta),\hat{\bP})$, we have for any $y<0$ and $\theta\in\s $,
\begin{eqnarray}
\hat{\bP}_{y,\theta}\left[\hat{H}^{+}_{\Theta_{t}}(\xi_{t})\1_{\{t<\tau^{+}_{0}\}}\right]
&=&\hat{\bP}_{y,\theta}\left[\hat{\bP}_{\xi_{t},\Theta_{t}}\left(\tau^{+}_{0}=+\infty\right)\1_{\{t<\tau^{+}_{0}\}}\right]\nonumber\\
&=&\hat{\bP}_{y,\theta}\left(\tau^{+}_{0}=+\infty\right)=\hat{H}^{+}_{\theta}(y).\nonumber
\end{eqnarray}
Using this and the Markov property of $((\xi,\Theta),\hat{\bP})$ we prove that $\hat{H}^{+}_{\Theta_{t}}(\xi_{t})1_{\{t<\tau^{+}_{0}\}}$ is a $\hat{\bP}$-martingale.\qed

\smallskip

Under the conditions of Proposition \ref{prop:martingale} we can define probability measures $\hat{\bP}^{\downarrow}_{y,\theta}$ on the Skorokhod space $\mathbb{D}_{\R\times\s}$ by
\begin{equation}\nonumber
\left.\frac{{\rm d}\hat{\bP}^{\downarrow}_{y,\theta}}{{\rm d}\hat{\bP}_{y,\theta}}\right|_{\mathcal{F}_{t}}:=\frac{\hat{H}^{+}_{\Theta_{t}}(\xi_{t})}{\hat{H}^{+}_{\theta}(y)}1_{\{t<\tau^{+}_{0}\}}\quad
\forall y<0,\ \theta\in\s ,\ t\ge 0.
\end{equation}
It follows by the theory of Doob's $h$-transform that for every $y<0$ and $\theta\in\s$ the process $((\xi,\Theta),\hat{\bP}^{\downarrow}_{y,\theta})$ is a strong Markov process on the state space $(0,+\infty)\times \s $ with semigroup $(\hat{P}^{\downarrow}_{t})_{t\ge 0}$ given by
$$\hat{P}^{\downarrow}_{t}f(z,\nu)=\frac{1}{\hat{H}^{+}_{\theta}(z)}\hat{\bP}_{z,\nu}\left[\hat{H}^{+}_{\Theta_{t}}(\xi_{t})f(\xi_{t},\Theta_{t})1_{\{t<\tau^{+}_{0}\}}\right]\quad\forall z<0,\nu\in\s,t\ge 0.$$
Since $\hat{H}^{+}_{\Theta_{t}}(\xi_{t})1_{\{t<\tau^{+}_{0}\}}$ is a $\hat{\bP}$-martingale, the semigroup $(\hat{P}^{\downarrow}_{t})_{t\ge 0}$ is Markovian and accordingly the process $((\xi,\Theta),\hat{\bP}^{\downarrow})$ has an infinite lifetime .

\begin{proposition}\label{prop:cond <0}
Suppose that \eqref{prea7} holds.
For all $y<0$, $\theta\in\s $, $t\ge 0$ and $\Lambda\in \mathcal{F}_{t}$,
\begin{equation}\nonumber
\hat{\bP}^{\downarrow}_{y,\theta}\left(\Lambda\right)=\lim_{q\to 0+}\hat{\bP}_{y,\theta}\left(\Lambda,t<\mathbf{e}_{q}\,|\,\tau^{+}_{0}>\mathbf{e}_{q}\right).
\end{equation}
\end{proposition}

\proof Note that by the Markov property of $((\xi,\Theta),\hat{\bP})$,
\begin{eqnarray}
\hat{\bP}_{y,\theta}\left(\Lambda;\ t<\mathbf{e}_{q}<\tau^{+}_{0}\right)
&=&\int_{t}^{+\infty}q \mathrm{e}^{-qs}\hat{\bP}_{y,\theta}(\Lambda;\ s<\tau^{+}_{0}){\rm d}s\nonumber\\
&=&\int_{0}^{+\infty}q \mathrm{e}^{-q(s+t)}\hat{\bP}_{y,\theta}(\Lambda;\ s+t<\tau^{+}_{0}){\rm d}s\nonumber\\
&=&\mathrm{e}^{-qt}\hat{\bP}_{y,\theta}\left(1_{\{\Lambda,t<\tau^{+}_{0}\}}\hat{\bP}_{\xi_{t},\Theta_{t}}(\tau^{+}_{0}>\mathbf{e}_{q})\right).\nonumber
\end{eqnarray}
Thus by the bounded convergence theorem,
\begin{eqnarray}
\lim_{q\to 0+}\hat{\bP}_{y,\theta}\left(\Lambda,\, t<\mathbf{e}_{q}\,|\,\tau^{+}_{0}>\mathbf{e}_{q}\right)
&=&\lim_{q\to 0+}\mathrm{e}^{-qt}\hat{\bP}_{y,\theta}\left(1_{\{\Lambda,t<\tau^{+}_{0}\}}\frac{\hat{\bP}_{\xi_{t},\Theta_{t}}(\tau^{+}_{0}>\mathbf{e}_{q})}{\hat{\bP}_{y,\theta}(\tau^{+}_{0}>\mathbf{e}_{q})}\right)\nonumber\\
&=&\hat{\bP}_{y,\theta}\left(\frac{\hat{H}^{+}_{\Theta_{t}}(\xi_{t})}{\hat{H}^{+}_{\theta}(y)}1_{\{\Lambda,t<\tau^{+}_{0}\}}\right)=\hat{\bP}^{\downarrow}_{y,\theta}\left(\Lambda\right).\nonumber
\end{eqnarray}\qed

\smallskip
The process $((\xi,\Theta),\hat{\bP}^{\downarrow})$ is referred to as the \textit{MAP conditioned to stay negative}.

\begin{proposition}\label{prop:cond at 0}
{Suppose that \eqref{prea7} holds.}
For every $\theta\in\s $, there exists a probability measure $\hat{\bP}^{\downarrow}_{0,\theta}$ on $\mathbb{D}_{\R\times\s}$ satisfying that
$\xi_{0}=0$ and $\xi_{t}\not=0$ for all $t>0,$ $\hat{\bP}^{\downarrow}_{0,\theta}$-a.s., and that
the process $(\xi_{t},\Theta_{t})_{t>0}$ under $\hat{\bP}^{\downarrow}_{0,\theta}$ is a strong Markov process with the same transition rates as $((\xi,\Theta), \{\hat{\bP}^{\downarrow}_{y,\theta}:y<0,\theta\in\s\})$. Moreover, we have
\begin{equation}\label{def:condi at 0}
\hat{\bP}^{\downarrow}_{0,\theta}\left[f(\xi_{t},\Theta_{t})\1_{\{t<\zeta\}}\right]=\frac{\hat{{\rm n}}^{+}_{\theta}
\left[\hat{H}^{+}_{\nu_{t}}(-\epsilon_{t})f(-\epsilon_{t},\nu_{t})\1_{\{t<\zeta\}}\right]}{\hat{{\rm n}}^{+}_{\theta}(\zeta=+\infty)}
\end{equation}
for any $t>0$ and nonnegative measurable function $f:\R\times \s\to\R^{+}$.
\end{proposition}
\proof
{Suppose $\hat{\mathfrak{P}}$ is the kernel from $\s\times \R\times\R^{+}$ to $\mathbb{D}_{\s\times\R\times \R^{+}}$ with respect to the process $((\Theta_{t},\xi_{t},U_{t})_{t\ge 0},\hat{\bP})$ defined in the same way of \cite{M} (see also the arguments in Section \ref{sec:fluctuation}).}
{Then under $\hat{\mathfrak{P}}^{\theta,0,0}$ the process $(\Theta_{t},\xi_{t},U_{t})_{t\ge 0}$ starts from $(\theta,0,0)$ and $(\Theta_{t},\xi_{t},U_{t})_{t>0}$ has the strong Markov property with respect to the same transition semigroup as $((\Theta_{t},\xi_{t},U_{t})_{t\ge 0},\hat{\bP}_{0,\theta})$. Let $\bar{R}=\inf\{t>0:t\in\bar{M}^{cl}\}$.}
Note that $\hat{H}^{+}_{\theta}(y)=\lim_{q\to 0+}\hat{\bP}_{y,\theta}\left(\tau^{+}_{0}>\mathbf{e}_{q}\right)$ for $y<0$ and $\theta\in\s $. It follows from the Markov property and the bouned convergence theorem that
\begin{eqnarray}
\hat{\mathfrak{P}}^{\theta,0,0}\left[\hat{H}^{+}_{\Theta_{t}}(\xi_{t})\1_{\{t<\bar{R}\}}\right]
&=&\lim_{q\to 0+}\hat{\mathfrak{P}}^{\theta,0,0}\left[\hat{\bP}_{\xi_{t},\Theta_{t}}\left(\tau^{+}_{0}>\mathbf{e}_{q}\right)\1_{\{t<\bar{R}\}}\right]\nonumber\\
&=&\lim_{q\to 0+}\mathrm{e}^{q t}\hat{\mathfrak{P}}^{\theta,0,0}\left(t<\mathbf{e}_{q}<\bar{R}\right)\nonumber\\
&=&\lim_{q\to 0+}\hat{{\rm n}}^{+}_{\theta}\left(t<\mathbf{e}_{q}<\zeta\right)\nonumber\\
&=&\hat{{\rm n}}^{+}_{\theta}(\zeta=+\infty).\nonumber
\end{eqnarray}
Thus we can define a probability measure $\hat{\bP}^{\downarrow}_{0,\theta}$ on $\mathbb{D}_{\R\times \s}$ by
\begin{equation}
\hat{\bP}^{\downarrow}_{0,\theta}(A):=\frac{1}{\hat{{\rm n}}^{+}_{\theta}(\zeta=+\infty)}\hat{\mathfrak{P}}^{\theta,0,0}\left[\hat{H}^{+}_{\Theta_{t}}(\xi_{t})\1_{\{t<\bar{R}\}}\1_{A}\right]\quad\forall A\in\mathcal{F}_{t},\ t>0.\label{prop3.15.1}
\end{equation}
One can easily show from the properties of $\hat{\mathfrak{P}}^{\theta,0,0}$ that under $\hat{\bP}^{\downarrow}_{0,\theta}$ the process $\xi_{t}$ leaves $0$ instantaneously and never hits $0$ again, and that the process $(\xi_{t},\Theta_{t})_{t\ge 0}$ is a Markov process whose transition rates  satisfy
$$\hat{\bP}^{\downarrow}_{0,\theta}\left[\xi_{t+s}\in A,\Theta_{t+s}\in B\,|\, \xi_{s},\Theta_{s}\right]
=\hat{\bP}^{\downarrow}_{\xi_{s},\Theta_{s}}\left[\xi_{t}\in A,\Theta_{t}\in B\right]$$
for all $t,s\ge 0$, $A\in \mathcal{B}(\R)$ and $B\in\mathcal{B}(\s )$.
Note that, by definition, under $\hat{\mathfrak{P}}^{\theta,0,0}$, $U_{t}$ equals $-\xi_{t}$ for $t<\bar{R}$. Hence by \eqref{prop3.15.1} for every $t>0$ and nonnegative measurable function $f:\R\times \s \to \R^{+}$, we have
\begin{eqnarray}
\hat{\bP}^{\downarrow}_{0,\theta}\left[f(\xi_{t},\Theta_{t})\1_{\{t<\zeta\}}\right]
&=&\frac{1}{\hat{{\rm n}}^{+}_{\theta}(\zeta=+\infty)}\hat{\mathfrak{P}}^{\theta,0,0}\left[\hat{H}^{+}_{\Theta_{t}}(-U_{t})f(-U_{t},\Theta_{t})\1_{\{t<\bar{R}\}}\right]\nonumber\\
&=&\frac{1}{\hat{{\rm n}}^{+}_{\theta}(\zeta=+\infty)}\hat{{\rm n}}^{+}_{\theta}\left[\hat{H}^{+}_{\nu_{t}}(-\epsilon_{t})f(-\epsilon_{t},\nu_{t})\1_{\{t<\zeta\}}\right].\nonumber
\end{eqnarray}
In the second equality we use the fact that $\hat{{\rm n}}^{+}_{\theta}$ is the image measure of $(U_{t},\Theta_{t})_{t<\bar{R}}$ under $\hat{\mathfrak{P}}^{\theta,0,0}$.\qed

\smallskip

\begin{remark}
\rm
 Suppose $\s=\{s_{1},s_{2},\cdots,s_{n}\}$ is a finite set and $((\xi,\Theta),\bP)$ is a MAP taking values in $\R\times\s$. For simplicity we assume the random variables $\Xi_{j,k}$ introduced in Example \ref{eg1} are such that $\Xi_{j,k}\stackrel{d}{=}\Xi_{k,j}$ for all $1\le j,k\le n$. Suppose the process $(\Theta,\{\bP_{x,\theta},\theta\in\s\})$ is irreducible and hence ergodic. Its invariant distribution is denoted by $\pi=(\pi_{1},\pi_{2},\cdots,\pi_{n})$. In this case condition \eqref{condi:weak reversability} is satisfied by taking $\widetilde{\bP}_{0,v}$ to be $\bP_{0,v}$. Hence $\hat{\bP}_{0,v}$ is the law of $(-\xi,\Theta)$ under $\bP_{0,v}$.
    Let $\hat{\phi}_{j}(q):=\hat{{\rm n}}^{+}_{j}(1-\mathrm{e}^{-q\zeta})$ for $1\le j\le n$ and $q>0$. It is proved in \cite{DDK} that
    $$\lim_{q\to 0+}\frac{\hat{\phi}_{j}(q)}{\hat{\phi}_{k}(q)}=\lim_{q\to 0+}\frac{\hat{{\rm n}}^{+}_{j}(\zeta=+\infty)+\hat{{\rm n}}^{+}_{j}(1-\mathrm{e}^{-q\zeta},\zeta<+\infty)}
    {\hat{{\rm n}}^{+}_{k}(\zeta=+\infty)+\hat{{\rm n}}^{+}_{k}(1-\mathrm{e}^{-q\zeta},\zeta<+\infty)}=\frac{\pi_{j}}{\pi_{k}}.$$
    It follows that if $\hat{{\rm n}}^{+}_{j}(\zeta=+\infty)>0$ for some (then for all) $1\le j\le n$, then there is a constant $c>0$ independent of $j$ such that $\hat{{\rm n}}^{+}_{j}(\zeta=+\infty)=c\pi_{j}$. Since $\hat{\bP}_{y,s_{i}}\left(\tau^{+}_{0}=+\infty\right)=\lim_{q\to 0+}\hat{\bP}_{0,s_{i}}\left(\bar{\xi}_{\mathbf{e}_{q}}\le -y\right)$, we get by Proposition \ref{prop:wiener_hopf} and the bounded convergence theorem that $$\hat{\bP}_{y,s_{i}}\left(\tau^{+}_{0}=+\infty\right)=c\sum_{j=1}^{n}\hat{U}^{+}_{ij}(-y)\pi_{j}$$
    where $\hat{U}^{+}_{ij}(-y)=\hat{\bP}_{0,s_{i}}\left[\int_{0}^{\bar{L}_{\infty}}1_{\{\xi^{+}_{t}\le -y, \Theta^{+}_{t}=s_{j}\}}{\rm d}t\right]$. In \cite{DDK}, $\sum_{j=1}^{n}\hat{U}^{+}_{ij}(-y)\pi_{j}$ is used as the harmonic function to define a martingale change of measure under which the MAP is conditioned to stay negative.

\end{remark}
\begin{remark}
\rm
 Suppose $((\xi,\Theta),\bP)$ is a MAP where $\xi$ is a (possibly killed) L\'{e}vy process on $\mathbb{R}$ whose law is independent of $\Theta$ and $\Theta$ has an invariant distribution.
In this case condition \eqref{condi:weak reversability} is satisfied by taking $\widetilde{\bP}_{0,v}=\bP_{0,v}$ and hence $\hat{\bP}_{0,v}$ is the law of $(-\xi,\Theta)$ under $\bP_{0,v}$. We assume that for $\xi$, $0$ is regular for both $(-\infty,0)$ and $(0,+\infty)$, in which case, both $((\xi,\Theta),\bP)$ and $((-\xi,\Theta),\bP)$ are upwards regular.
We claim that \eqref{prea7} is
 satisfied if and only if the L\'{e}vy process $\xi_{t}$ drifts to $+\infty$.
 To see this, we first recall some known facts about L\'{e}vy processes. Let $\underline{L}_{t}$ be the local time of $\xi$ at the running minima and ${\rm n}^{-}$ be the excursion measures at the minimum. In fact, ${\rm n}^{-}$ equals $\hat{{\rm n}}^{+}$ which is the excursion measure at the maximum of the dual process $-\xi$. Since  $0$ is regular for $(-\infty,0)$, there is a continuous version of $\underline{L}_{t}$ and a
 {nonnegative}
 constant ${l}^{-}$ such that almost surely
 $\int_{0}^{t}\1_{\{\xi_{s}=\inf_{r\in [0,s]}\xi_{r}\}}{\rm d}s={l}^{-}\,\underline{L}_{t}$
 for all $t\ge 0$. In this case, the inverse local time $\underline{L}^{-1}_{t}$ is a
 {(possibly killed) subordinator}
 with Laplace exponent given by $\hat{\Phi}(q)={l}^{-}q+{\rm n}^{-}(1-\mathrm{e}^{-q\zeta})$.
 It follows that $\underline{L}_{\infty}$ is exponentially distributed with parameter ${\rm n}^{-}(\zeta=+\infty)$. Hence
 ${\rm n}^{-}(\zeta=+\infty)>0$ if and only if $\xi_{t}$ drifts to $+\infty$, in which case \cite{CD2005} showed further that ${\rm n}^{+}(\zeta)=l^++{\rm n}^{+}\left(1-\mathrm{e}^{-\zeta}\right) <+\infty$ where ${\rm n}^{+}$ denotes the excursion measure at the maximum of $\xi$ and $l^+$ is the drift coefficient for the inverse local time at the maximum.
 \end{remark}

\section{Stationary overshoots and undershoots of MAPs}\label{sec:stationary overshoots}

Throughout this section we will assume that the modulator of $((\xi,\Theta),\bP)$
\begin{equation}
 \label{prea2}
\Theta\text{ is positive recurrent with invariant distribution }\pi\text{ which is fully supported on }\s.
\end{equation}

\begin{definition}
For $q>0$, let $\{T^{(q)}_{n}:n\ge 0\}$ be a sequence of random variables such that $T^{(q)}_{0}=0$ and $\{T^{(q)}_{n+1}-T^{(q)}_{n}:n\ge 0\}$ are independent and exponentially distributed random variables with mean $1/q$. Define
$$M^{(q),+}_{n}:=\Theta^{+}_{T^{(q)}_{n}}\quad\forall n\ge 0.$$
We call $M^{(q),+}: = \{M^{(q),+}_{n}:n\ge 0\}$ the
$q$-embedded chain
of the process $(\Theta^{+}_{t})_{t\ge 0}$. Moreover, in the spirit of \cite{MT}, we say that $\Theta^+$ is a (nonarithmetic aperiodic) Harris recurrent process if $\Theta^{+}$ has a (nonarithmetic aperiodic) Harris recurrent $q$-embedded chain for some $q>0$.
\end{definition}

Under the preceding assumption \eqref{prea2},
together with the assumption that
\begin{equation}
\label{prea8}
\inf_{v\in\s}\left[\ell^{+}(v)+{\rm n}^{+}_{v}\left(1-\mathrm{e}^{-\zeta}\right)\right]>0\text{  and  }{\rm n}^{+}_{v}(\zeta)<+\infty\text{ for every }v\in\s,
\end{equation}
it follows by Corollary \ref{cor:existence of invariant distribution} that
\begin{equation}\label{def:pi+}
\pi^{+}(\cdot)=\frac{1}{\bP_{0,\pi}\left[\bar{L}_{1}\right]}\bP_{0,\pi}\left[\int_{0}^{1}\1_{\{\Theta_{s}\in\cdot\}}{\rm d}\bar{L}_{s}\right]
\end{equation}
is an invariant distribution for $\Theta^{+}$ and hence for $M^{(q),+}$. It  follows by \cite[Theorem (5.1)]{Kaspi83} that
\begin{equation}\nonumber
\pi({\rm d}v)=\frac{1}{c_{\pi^{+}}}\left[\ell^+(v)\pi^{+}({\rm d}v)+\int_{\s }{\rm n}^{+}_{\theta}\left(\int_{0}^{\zeta}\1_{\{\nu_{t}\in {\rm d}v\}}{\rm d}t\right)\pi^{+}({\rm d}\theta)\right]
\end{equation}
where $c_{\pi^{+}}:=\int_{\s }\left[\ell^+(\theta)+{\rm n}^{+}_{\theta}(\zeta)\right]\pi^{+}({\rm d}\theta)$ is a positive constant.

\begin{lemma}\label{fromKaspi}
Assume that \eqref{prea2} and \eqref{prea8} hold and, further, that
$\bP_{0,\pi^{+}}\left[\xi^{+}_{1}\right]<+\infty$ where $\pi^{+}$ given in \eqref{def:pi+} is fully supported on $\s$.
Suppose that the continuous part of $\xi^{+}$ can be represented by $\int_{0}^{t}a^{+}(\Theta^{+}_{s}){\rm d}s$ for some {nonnegative} measurable function $a^{+}$ on $\s$.
Then for all $q>0$, we have
\begin{equation}\nonumber
\mu^{+}:=\int_{\s }a^{+}(\phi)\pi^{+}({\rm d}\phi)+\int_{\s \times \R^{+}}\bar{\Pi}^{+}_{\phi}(y)\pi^{+}({\rm d}\phi){\rm d}y=q\bP_{0,\pi^{+}}\left[\xi^{+}_{\mathbf{e}_{q}}\right]<+\infty,
\end{equation}
where $\bar{\Pi}^{+}_{\phi}(y):=
\Pi^{+}(\phi,\s, (y,+\infty))
$.
\end{lemma}

\proof
Using that $\bP_{0,\pi^{+}}\left[\xi^{+}_{1}\right]<+\infty$
and the subadditivity of $t\mapsto \bP_{0,\pi^{+}}\left[\xi^{+}_{t}\right]$, one can show in the same way as in the proof of Lemma \ref{lem:A.4} that
$\bP_{0,\pi^{+}}\left[\xi^{+}_{t}\right]<+\infty$ for all $t>0$ and $\bP_{0,\pi^{+}}\left[\xi^{+}_{\mathbf{e}_{q}}\right]<+\infty$ for all $q>0$. We note that for every $t>0$,
$$\xi^{+}_{t}=\int_{0}^{t}a^{+}(\Theta^{+}_{s}){\rm d}s+\sum_{0\le s\le t}\Delta \xi^{+}_{s}\1_{\{\Delta\xi^{+}_{s}>0\}},$$
where $\Delta\xi^{+}_{s}=\xi^{+}_{s}-\xi^{+}_{s-}$.
By Proposition \ref{prop:Pi+}
and Fubini's theorem, we have
\begin{eqnarray}
\bP_{0,\theta}\left[\sum_{0\le s\le t}\Delta \xi^{+}_{s}\1_{\{\Delta\xi^{+}_{s}>0\}}\right]
&=&\bP_{0,\theta}\left[\int_{0}^{t}{\rm d}s\int_{\s\times \R^{+}}y\Pi^{+}(\Theta^{+}_{s},{\rm d}\phi,{\rm d}y)\right]\nonumber\\
&=&\bP_{0,\theta}\left[\int_{0}^{t}{\rm d}s\int_{0}^{+\infty}\bar{\Pi}^{+}_{\Theta^{+}_{s}}(z){\rm d}z\right]
\end{eqnarray}
for every $\theta\in\s$. Hence
\begin{eqnarray}
\bP_{0,\pi^{+}}[\xi^{+}_{t}]&=&\bP_{0,\pi^{+}}\left[\int_{0}^{t}{\rm d}s\left(a^{+}(\Theta^{+}_{s})+\int_{0}^{+\infty}\bar{\Pi}^{+}_{\Theta^{+}_{s}}(z){\rm d}z\right)\right]\nonumber\\
&=&t\left(\int_{\s }a^{+}(\phi)\pi^{+}({\rm d}\phi)+\int_{\s \times \R^{+}}\bar{\Pi}^{+}_{\phi}(z)\pi^{+}({\rm d}\phi){\rm d}z\right)\nonumber\\
&=&t\mu^{+}.\nonumber
\end{eqnarray}
In the second equality we use the fact that $\pi^{+}$ is an invariant distribution for $\Theta^{+}$.
Consequently we have
$$\bP_{0,\pi^{+}}\left[\xi^{+}_{\mathbf{e}_{q}}\right]
=q\int_{0}^{+\infty}\mathrm{e}^{-q t}\bP_{0,\pi^{+}}\left[\xi^{+}_{t}\right]{\rm d}t=\frac{\mu^{+}}{q}.$$\qed

\smallskip

Under the assumptions of  Lemma \ref{fromKaspi}, the measure $\rho^{\ominus}$ given below is a probability measure on $\R^{+}\times \s$,
\begin{equation}\label{def:rho-}
\rho^{\ominus}({\rm d}z,{\rm d}v):=\frac{1}{\mu^{+}}\left[a^{+}(v)\pi^{+}({\rm d}v)\delta_{0}({\rm d}z)+\1_{\{z>0\}}\int_{\s \times \R^{+}}\pi^{+}({\rm d}\phi){\rm d}y\,\Pi^{+}(\phi,{\rm d}v,{\rm d}z+y)\right].
\end{equation}
We will show in the following that $\rho^{\ominus}$ is the stationary distribution for the overshoots of the MAP, assuming additionally that,
the modulator
\begin{equation}\label{prea6}
\Theta^{+}\mbox{ of }((\xi^{+},\Theta^{+}),\bP) \mbox{ is a nonarithmetic aperiodic Harris recurrent process.}
\end{equation}
The key of the proof is an application of Markov renewal theory developed in \cite{Alsmeyer1994}.
Suppose that $\{M^{(q),+}_{n}=\Theta^{+}_{T^{(q)}_{n}}:n\ge 0\}$ is a nonarithmetic aperiodic Harris recurrent $q$-embedded chain of $(\Theta^{+},\bP)$.
Define
$$S^{(q),+}_{n}:=\xi^{+}_{T^{(q)}_{n}},\quad N^{(q),+}_{n}:=\bar{L}^{-1}_{T^{(q)}_{n}}\quad \forall n\ge 0.$$
One can show that both $(M^{(q),+}_{n},S^{(q),+}_{n})_{n\ge 0}$ and $(M^{(q),+}_{n},N^{(q),+}_{n})_{n\ge 0}$ are Markov renewal processes in the sense of \cite{Alsmeyer1994}. We shall first consider the process $(M^{(q),+}_{n},S^{(q),+}_{n})_{n\ge 0}$.
For every $\theta\in \s $, let
\begin{equation}
F_{\theta}({\rm d}v,{\rm d}z):=\bP_{0,\theta}(M^{(q),+}_{1}\in {\rm d}v,S^{(q),+}_{1}\in {\rm d}z)=\int_{0}^{+\infty}q\mathrm{e}^{-q t}\bP_{0,\theta}\left(\Theta^{+}_{t}\in {\rm d}v,\ \xi^{+}_{t}\in {\rm d}z\right){\rm d}t.\label{def:renewaldist}
\end{equation}
Let $F^{*0}_{\theta}({\rm d}v,{\rm d}z):=\delta_{\theta}({\rm d}v)\delta_{0}({\rm d}z)$ and $F^{*n}_{\theta}$ be the $n$-th convolution of $F_{\theta}$ for $n\ge 1$.
Then $\sum_{n=0}^{+\infty}F^{*n}_{\theta}({\rm d}v,{\rm d}z)$ is the renewal measure of $(M^{(q),+}_{n},S^{(q),+}_{n})_{n\ge 0}$.
Note that $\bP_{0,\pi^{+}}\left[S^{(q),+}_{1}\right]=\bP_{0,\pi^{+}}\left[\xi^{+}_{\mathbf{e}_{q}}\right]=\mu^{+}/q$.
Given \eqref{prea6},
it follows by \cite[Theorem 2.1]{Alsmeyer1994} that
\begin{equation}
\lim_{y\to +\infty}\int_{\s \times [0,y]}g(v,y-z)\sum_{n=0}^{+\infty}F^{*n}_{\theta}({\rm d}v,{\rm d}z)=\frac{q}{\mu^{+}}\int_{\s \times \R^{+}}g(v,z)\pi^{+}({\rm d}v){\rm d}z\label{eq:limit of Fn}
\end{equation}
for every $\theta\in\s$ and every measurable function $g:\s \times \R^{+}\to \R$ satisfying the following two conditions:
\begin{itemize}
\item[(i)] for each $v\in\s$, the set of discontinuous points of $z\mapsto g(v,z)$ has zero Lebesgue measure;
\item[(ii)] $\int_{\s}\sum_{n=0}^{+\infty}\sup_{z\in[np,(n+1)p)}|g(v,z)|\pi^{+}({\rm d}v)<+\infty$ for some $p>0$.
\end{itemize}
We use $\mathcal{M}$ to denote the space of measurable functions satisfying both of the above conditions.
In view of the fact that $\bP_{0,\theta}\left(T^{(q)}_{n}\in {\rm d}t\right)=\frac{q^{n}t^{n-1}}{(n-1)!}\rme^{- q t}{\rm d}t$ for $n\ge 1$, we have
\begin{eqnarray}
U^{+}_{\theta}({\rm d}v,{\rm d}z)&=&\int_{0}^{+\infty}\bP_{0,\theta}\left(\Theta^{+}_{t}\in {\rm d}v,\ \xi^{+}_{t}\in {\rm d}z\right){\rm d}t\nonumber\\
&=&\sum_{n=1}^{+\infty}\int_{0}^{+\infty}\mathrm{e}^{-q t}\frac{(q t)^{n-1}}{(n-1)!}\bP_{0,\theta}\left(\Theta^{+}_{t}\in {\rm d}v,\ \xi^{+}_{t}\in {\rm d}z\right){\rm d}t\nonumber\\
&=&\frac{1}{q}\sum_{n=1}^{+\infty}\bP_{0,\theta}\left(M^{(q),+}_{n}\in {\rm d}v,\ S^{(q),+}_{n}\in {\rm d}z\right)\nonumber\\
&=&\frac{1}{q}\left[\sum_{n=0}^{+\infty}F^{*n}_{\theta}({\rm d}v,{\rm d}z)-\delta_{\theta}({\rm d}v)\delta_{0}({\rm d}z)\right].\nonumber
\end{eqnarray}
This and \eqref{eq:limit of Fn} imply that for every $\theta\in\s$ and every $g\in\mathcal{M}$,
\begin{equation}
\lim_{y\to+\infty}\int_{\s \times [0,y]}g(v,y-z)U^{+}_{\theta}({\rm d}v,{\rm d}z)=\frac{1}{\mu^{+}}\int_{\s \times \R^{+}}g(v,z)\pi^{+}({\rm d}v){\rm d}z.\label{eq:limit of U+}
\end{equation}

\smallskip

\begin{remark}\label{re:sufficient for M}
\rm
It is easy to see that $g\in\mathcal{M}$ if, in particular, $z\mapsto g(v,z)$ is right continuous on $[0,+\infty)$ and there is a measurable function $f:\s\times\R^{+}\to\R^{+}$ such that $|g(v,z)|\le f(v,z)$ for all $(v,z)\in\s\times\R^{+}$, $z\mapsto f(v,z)$ is a monotone function on $\R^{+}$ and $\int_{\s\times\R^{+}}f(v,z)\pi^{+}({\rm d}v){\rm d}z<+\infty$. In fact this sufficient condition for $g\in\mathcal{M}$ is easy to be verified and will be used later in our proofs where the Markov renewal theory is applied.
\end{remark}

\begin{proposition}\label{prop:stationary overshoots}
Suppose
{\eqref{prea6} and the conditions in Proposition \ref{prop:sufficient} \and Lemma \ref{fromKaspi} hold.}
For every $\theta\in \s$, the joint probability measures on $\s\times \R^{-}\times \s\times \R^{+}$
$$\bP_{0,\theta}\left(\Theta_{\tau^{+}_{x}-}\in {\rm d} v,\xi_{\tau^{+}_{x}-}-x\in {\rm d} y, \Theta_{\tau^{+}_{x}}\in {\rm d} \phi, \xi_{\tau^{+}_{x}}-x\in {\rm d} z\right)$$
converges weakly to a probability measure $\rho$ given by
\begin{eqnarray}
\rho({\rm d} v,{\rm d} y,{\rm d} \phi, {\rm d} z)
&:=&\frac{1}{\mu^{+}}\big[\1_{\{y<0,z>0\}}\ell^+(v)\Pi(v,{\rm d} \phi,{\rm d} z-y)\pi^{+}({\rm d} v){\rm d} y\nonumber\\
&&\left.+\1_{\{y<0,z>0\}}{\rm d} y\int_{\s}\pi^{+}({\rm d} \varphi){\rm n}^{+}_{\varphi}\big(\int_{0}^{\zeta}\1_{\{\epsilon_{s}\le -y,\nu_{s}\in{\rm d} v\}}\Pi(v,{\rm d}\phi,{\rm d} z-y){\rm d} s\big)\right.\nonumber\\
&&+a^{+}(v)\pi^{+}({\rm d} v)\delta_{0}({\rm d} y)\delta_{0}({\rm d} z)\delta_{v}({\rm d} \phi)\big]\nonumber
\end{eqnarray}
as $x\to +\infty$.
In particular,
$\bP_{0,\theta}\left(\xi_{\tau^{+}_{x}}-x\in {\rm d} z,\ \Theta_{\tau^{+}_{x}}\in {\rm d} \phi\right)$
converges weakly to $\rho^{\ominus}({\rm d} z,{\rm d} \phi)$ given by \eqref{def:rho-},
and $\bP_{0,\theta}\left(\xi_{\tau^{+}_{x}-}-x\in {\rm d} y,\ \Theta_{\tau^{+}_{x}-}\in {\rm d} v\right)$ converges weakly to a probability measure $\rho^{\oplus}({\rm d} y,{\rm d} v)$ given by
\begin{eqnarray}
\rho^{\oplus}({\rm d} y,{\rm d} v)&:=&\frac{1}{\mu^{+}}\big[a^{+}(v)\pi^{+}({\rm d} v)\delta_{0}({\rm d} y)+\1_{\{y<0\}}
\ell^+(v)\bar{\Pi}_{v}(-y)\pi^{+}({\rm d} v){\rm d} y\nonumber\\
&&+\1_{\{y<0\}}{\rm d} y\int_{\s}\pi^{+}({\rm d} \phi){\rm n}^{+}_{\phi}\big(\int_{0}^{\zeta}\bar{\Pi}_{v}(-y)\1_{\{\epsilon_{r}\le -y,\nu_{r}\in {\rm d} v\}}{\rm d} r\big)\big].\nonumber
\end{eqnarray}
Here $\bar{\Pi}_{v}(-y):=\Pi(v,\s,(-y,+\infty))$.
\end{proposition}

\proof First we claim that $\rho$ given above is a probability measure. Integrating $\rho({\rm d} v,{\rm d} y,{\rm d} \phi, {\rm d} z)$ over the variables $v$ and $y$, we get that
\begin{align}
\frac{1}{\mu^{+}}&\big[a^{+}(\phi)\pi^{+}({\rm d}\phi)\delta_{0}({\rm d}z)+\1_{\{z>0\}}\int_{\s\times \R^{+}}\pi^{+}({\rm d}v){\rm d}y\,\ell^+(v)\Pi(v,{\rm d}\phi,{\rm d}z+y)\nonumber\\
&\quad+\1_{\{z>0\}}\int_{\s}\pi^{+}(d\varphi){\rm n}^{+}_{\varphi}\big(\int_{0}^{\zeta}{\rm d}s\int_{\epsilon_{s}}^{+\infty}{\rm d}y\,\1_{\{\nu_{s}\in {\rm d}v\}}\Pi(\nu_{s},{\rm d}\phi,{\rm d}z+y)\big)\big]\nonumber\\
&=\frac{1}{\mu^{+}}\big[a^{+}(\phi)\pi^{+}({\rm d}\phi)\delta_{0}({\rm d}z)+\1_{\{z>0\}}\int_{\s\times \R^{+}}\pi^{+}({\rm d}v){\rm d}y\,\ell^+(v)\Pi(v,{\rm d}\phi,{\rm d}z+y)\nonumber\\
&\quad+\1_{\{z>0\}}\int_{\s\times \R^{+}}\pi^{+}(d\varphi){\rm d}u\,{\rm n}^{+}_{\varphi}\big(\int_{0}^{\zeta}\1_{\{\nu_{s}\in {\rm d}v\}}\Pi(\nu_{s},{\rm d}\phi,{\rm d}z+\epsilon_{s}+u)\big)\big]\nonumber\\
&=\frac{1}{\mu^{+}}\big[a^{+}(\phi)\pi^{+}({\rm d}\phi)\delta_{0}({\rm d}z)+\1_{\{z>0\}}\int_{\s\times \R^{+}}\pi^{+}({\rm d}v){\rm d}y\Pi^{+}(v,{\rm d}\phi,{\rm d}z+y)\big]\nonumber\\
&=\rho^{\ominus}({\rm d}z,{\rm d}\phi).\nonumber
\end{align}
The first equality follows from a change of variable and Fubini's theorem, and the second equality follows from Proposition \ref{prop:Pi+}. This implies that $\rho$ is a probability measure and $\rho^{\ominus}$ is its marginal law. Similarly, by integrating $\rho({\rm d} v,{\rm d} y,{\rm d} \phi, {\rm d} z)$ over the variables $\phi$ and $z$, we can show that $\rho^{\oplus}$ is also a marginal law of $\rho$. Next we prove the weak convergence.
Suppose $f,g:\s\times\R\to\R$ are bounded continuous functions. It follows by Proposition \ref{prop:under&overshoot} that for any $x>0$,
\begin{align}
\bP_{0,\theta}&\left[f(\Theta_{\tau^{+}_{x}-},\xi_{\tau^{+}_{x}-}-x)g(\Theta_{\tau^{+}_{x}}, \xi_{\tau^{+}_{x}}-x)\1_{\{\xi_{\tau^{+}_{x}>x}\}}\right]\nonumber\\
&=\int_{\s\times [0,x]}U^{+}_{\theta}({\rm d} v,{\rm d} z)\big[\ell^+(v)f(v,z-x)G(v,x-z)\nonumber\\
&\hspace{1cm}+{\rm n}^{+}_{v}\big(\int_{0}^{\zeta}f(\nu_{s},z-x-\epsilon_{s})G(\nu_{s},x-z+\epsilon_{s}){\rm d} s\big)\big].
\end{align}
where $G(v,u)=\int_{\s\times (u,+\infty)}g(\phi,y-u)\Pi(v,{\rm d} \phi,{\rm d} y)$.
One can easily show that the condition given in Remark \ref{re:sufficient for M} is satisfied by the function
$$(v,z)\mapsto \ell^+(v)f(v,-z)G(v,z)+{\rm n}^{+}_{v}\left(\int_{0}^{\zeta}f(\nu_{s},-z-\epsilon_{s})G(\nu_{s},z+\epsilon_{s}){\rm d} s\right).$$
Hence by \eqref{eq:limit of U+}, the integral in the right-hand side converges to
\begin{equation}\label{prop6.1.1}
\frac{1}{\mu^{+}}\int_{\s\times \R^{+}}\pi^{+}({\rm d} v){\rm d} z\left[\ell^+(v)f(v,-z)G(v,z)+{\rm n}^{+}_{v}\left(\int_{0}^{\zeta}f(\nu_{s},-z-\epsilon_{s})G(\nu_{s},z+\epsilon_{s}){\rm d} s\right)\right].
\end{equation}
By Fubini's theorem, we have
\begin{align}
\int_{\s\times \R^{+}}&\pi^{+}({\rm d} v){\rm d} z {\rm n}^{+}_{v}\left(\int_{0}^{\zeta}f(\nu_{s},-z-\epsilon_{s})G(\nu_{s},z+\epsilon_{s}){\rm d} s\right)\nonumber\\
&=\int_{\s}\pi^{+}({\rm d} v){\rm n}^{+}_{v}\left(\int_{0}^{+\infty}\int_{0}^{\zeta}f(\nu_{s},-z-\epsilon_{s})G(\nu_{s},z+\epsilon_{s}){\rm d} s{\rm d} z\right)\nonumber\\
&=\int_{\s}\pi^{+}({\rm d} v){\rm n}^{+}_{v}\left(\int_{0}^{\zeta}{\rm d} s \int_{\epsilon_{s}}^{+\infty}f(\nu_{s},-y)G(\nu_{s},y){\rm d} y\right)\nonumber\\
&=\int_{\s\times \R^{+}}\pi^{+}({\rm d} v){\rm d} y\, {\rm n}^{+}_{v}\left(\int_{0}^{\zeta}\1_{\{\epsilon_{s}\le y\}}f(\nu_{s},-y)G(\nu_{s},y){\rm d} s\right).\label{prop6.1.2}
\end{align}

Next we deal with the creeping event $\{\xi_{\tau^{+}_{x}}=x\}$.
Note that
\begin{align*}F_{\theta}({\rm d}v,{\rm d}z)&=\int_{0}^{+\infty}q\rme^{-q t}\bP_{0,\theta}\left(\Theta^{+}_{t}\in {\rm d}v,\xi^{+}_{t}\in {\rm d}z\right){\rm d}t\\
&=q\int_{0}^{+\infty}\bP_{0,\theta}\left(\Theta^{+}_{t}\in {\rm d}v,\xi^{+}_{t}\in {\rm d}z,t<\mathbf{e}_{q}\right){\rm d}t.
\end{align*}
This equation implies that $F_{\theta}({\rm d}v,{\rm d}z)/q$ can be viewed as the potential measure of the nondecreasing MAP $(\xi^{+},\Theta^{+})$ killed by an independent exponential time $\mathbf{e}_{q}$. In fact, we can verify that this killed process is a nondecreasing MAP and satisfies all the conditions in Lemma \ref{lem:on creeping}. Hence by Lemma \ref{lem:on creeping}
{$\1_{\{a^{+}(v)>0\}}F_{\theta}({\rm d}v,{\rm d}z)$} has a
{kernel}
${f}_{\theta}({\rm d}v,z)$ with respect to the Lebesgue measure ${\rm d}z$ such that
\begin{equation}\label{prop6.1.3}
\bP_{0,\theta}\left(h(\Theta^{+}_{T^{+}_{x}});\xi^{+}_{T^{+}_{x}}=x,T^{+}_{x}<\mathbf{e}_{q}\right)=\frac{1}{q}\int_{\s}a^{+}(v)h(v){f}_{\theta}({\rm d}v,x)
\end{equation}
for every nonnegative measurable function $h:\s\to\R^{+}$ and almost every $x>0$. We claim that $x\mapsto \bP_{0,\theta}\left(h(\Theta^{+}_{T^{+}_{x}});\xi^{+}_{T^{+}_{x}}=x,T^{+}_{x}<\mathbf{e}_{q}\right)=\bP_{0,\theta}\left(\rme^{-q T^{+}_{x}}h(\Theta^{+}_{T^{+}_{x}});\xi^{+}_{T^{+}_{x}}=x\right)$ is right continuous on $[0,+\infty)$ if in particular $h$ is a bounded continuous function. To see this, we take arbitrary $\{x_{n}:n\ge 1\}$ and $x\in\R^{+}$ satisfying $x_{n}\downarrow x$.
Since $\xi^{+}_{T^{+}_{x}}=\xi_{\tau^{+}_{x}}$ we have
\begin{align}
&\big|\bP_{0,\theta}\big(\rme^{-q T^{+}_{x_{n}}}h(\Theta^{+}_{T^{+}_{x_{n}}});\xi^{+}_{T^{+}_{x_{n}}}=x_{n}\big)-\bP_{0,\theta}\big(\rme^{-q T^{+}_{x}}h(\Theta^{+}_{T^{+}_{x}});\xi^{+}_{T^{+}_{x}}=x\big)\big|\nonumber\\
&\le\big|\bP_{0,\theta}\big[\rme^{-q T^{+}_{x_{n}}}h(\Theta^{+}_{T^{+}_{x_{n}}})\big(\1_{\{\xi^{+}_{T^{+}_{x_{n}}}=x_{n}\}}-\1_{\{\xi^{+}_{T^{+}_{x}}=x\}}\big)\big]\big|\nonumber\\
&\quad+\big|\bP_{0,\theta}\big[\rme^{-q T^{+}_{x_{n}}}h(\Theta^{+}_{T^{+}_{x_{n}}})-\rme^{-q T^{+}_{x}}h(\Theta^{+}_{T^{+}_{x}});\xi^{+}_{T^{+}_{x}}=x\big]\big|\nonumber\\
&\le\|h\|_{\infty}\bP_{0,\theta}\big(\{\xi_{\tau^{+}_{x_{n}}}=x_{n}\}\triangle \{\xi_{\tau^{+}_{x}}=x\}\big)+\bP_{0,\theta}\big[\big|\rme^{-q T^{+}_{x_{n}}}h(\Theta^{+}_{T^{+}_{x_{n}}})-\rme^{-q T^{+}_{x}}h(\Theta^{+}_{T^{+}_{x}})\big|\big].\nonumber
\end{align}
In view of \eqref{prop2.11.4} and the fact that $T^{+}_{x_{n}}\downarrow T^{+}_{x}$ and $\Theta^{+}_{T^{+}_{x_{n}}}\to \Theta^{+}_{T^{+}_{x}}$ $\bP_{0,\theta}$-a.s., we get by the above inequality and the bounded convergence theorem that
$$\lim_{n\to+\infty}\bP_{0,\theta}\left(\rme^{-q T^{+}_{x_{n}}}h(\Theta^{+}_{T^{+}_{x_{n}}});\xi^{+}_{T^{+}_{x_{n}}}=x_{n}\right)=\bP_{0,\theta}\left(\rme^{-q T^{+}_{x}}h(\Theta^{+}_{T^{+}_{x}});\xi^{+}_{T^{+}_{x}}=x\right).$$
Hence we prove the claim. Now we set  ${f}_{\theta}({\rm d}v,x)=\frac{q}{a^{+}(v)}\bP_{0,\theta}\left(\Theta^{+}_{T^{+}_{x}}\in {\rm d}v,\xi^{+}_{T^{+}_{x}}=x,T^{+}_{x}<\mathbf{e}_{q}\right)$ for every $x>0$ {in the setting that $a^+(v)>0$ for all $v\in \s$ and otherwise, we set $a^{+}(v){f}_{\theta}({\rm d}v,x)$ as identically equal to zero for all $v\in\s$}. The above arguments shows that $x\mapsto a^{+}(v){f}_{\theta}({\rm d}v,x)$ is right continuous on $(0,+\infty)$ in the sense of vague convergence and \eqref{prop6.1.3} holds for every $x>0$ and every nonnegative measurable function $h:\s\to\R$.
{ Reverting back to the setting $a^+(v) >0$ for all $v\in \s$}, since
\begin{eqnarray}
U^{+}_{\theta}({\rm d}v,{\rm d}z)&=&\frac{1}{q}\sum_{n=0}^{+\infty}F^{*(n+1)}_{\theta}({\rm d}v,{\rm d}z)\nonumber\\
&=&\frac{1}{q}\int_{\s\times [0,z]}F_{\phi}({\rm d}v,{\rm d}z-y)\sum_{n=0}^{+\infty}F^{*n}_{\theta}({\rm d}\phi,{\rm d}y),\nonumber
\end{eqnarray}
we can take the
{kernel $u^{+}_{\theta}({\rm d}v,z)$ of $U^{+}_{\theta}({\rm d}v,{\rm d}z)$}
to be such that
\begin{equation}\label{prop6.1.4}
u^{+}_{\theta}({\rm d}v,z)=\frac{1}{q}\int_{\s \times [0,z]}{f}_{\phi}({\rm d}v,z-y)\sum_{n=0}^{+\infty}F^{*n}_{\theta}({\rm d}\phi,{\rm d}y)\quad \forall z>0.
\end{equation}
For $n\ge 1$,
\begin{eqnarray}
F^{*n}_{\theta}({\rm d}v,{\rm d}z)&=&\int_{0}^{+\infty}\bP_{0,\theta}\left(\Theta^{+}_{T^{(q)}_{n}}\in {\rm d}v,\xi^{+}_{T^{(q)}_{n}}\in {\rm d}z\right){\rm d}t\nonumber\\
&=&\int_{0}^{+\infty}\frac{q^{n}t^{n-1}}{(n-1)!}\rme^{-q t}\bP_{\theta}(\Theta^{+}_{t}\in {\rm d}v,\xi^{+}_{t}\in {\rm d}z){\rm d}t,\nonumber
\end{eqnarray}
Obviously $F^{*n}_{\theta}({\rm d}v,{\rm d}z)$ is absolutely continuous with respect to $U^{+}_{\theta}({\rm d}v,{\rm d}z)$, and hence $F^{*n}_{\theta}({\rm d}v,{\rm d}z)$ has a
{kernel}
with respect to the Lebesgue measure ${\rm d}z$ which is denoted by $f^{*n}_{\theta}({\rm d}v,z)$. In view of this, $u^{+}_{\theta}({\rm d}v,z)$ given in \eqref{prop6.1.4} can be represented by
$$u^{+}_{\theta}({\rm d}v,z)=\frac{1}{q}{f}_{\theta}({\rm d}v,z)+\frac{1}{q}\int_{0}^{z}{\rm d}y\int_{\s}{f}_{\phi}({\rm d}v,z-y)\sum_{n=1}^{+\infty}f^{*n}_{\theta}({\rm d}\phi,y).$$
Using this expression and the fact that $z\mapsto a^{+}(v){f}_{\phi}({\rm d}v,z)$ is right continuous on $(0,+\infty)$, we can show that $x\mapsto a^{+}(v)u^{+}_{\theta}({\rm d}v,x)$ is right continuous on $(0,+\infty)$ in the sense of vague convergence. Hence $u^{+}_{\theta}({\rm d}v,z)$ given in \eqref{prop6.1.4} is the
{kernel}
taken in Proposition \ref{prop:sufficient}, and we have
\begin{align}
\bP_{0,\theta}&\left[f(\Theta_{\tau^{+}_{x}-},\xi_{\tau^{+}_{x}-}-x)g(\Theta_{\tau^{+}_{x}},\xi_{\tau^{+}_{x}}-x)\1_{\{\xi_{\tau^{+}_{x}}=x\}}\right]\nonumber\\
&=\int_{\s}a^{+}(v)f(v,0)g(v,0)u^{+}_{\theta}({\rm d}v,x)\nonumber\\
&=\frac{1}{q}\int_{\s\times [0,x]}\sum_{n=0}^{+\infty}F^{*n}_{\theta}({\rm d}\phi,{\rm d}y)\int_{\s}f(v,0)g(v,0)a^{+}(v){f}_{\phi}({\rm d}v,x-y)
\label{useMAPR}
\end{align}
for every $x>0$. Again by Remark \ref{re:sufficient for M} we can show that $(\phi,z)\mapsto \int_{\s}f(v,0)g(v,0)a^{+}(v){f}_{\phi}({\rm d}v,z)=q\bP_{0,\phi}\left[f(\Theta^{+}_{T^{+}_{z}},0)g(\Theta^{+}_{T^{+}_{z}},0);T^{+}_{z}<\mathbf{e}_{q}\right]\in\mathcal{M}$.
Hence by \eqref{eq:limit of Fn} the integral on the right-hand side of \eqref{useMAPR} converges, as $x\to+\infty$, towards
\begin{align}
\frac{1}{\mu^{+}}&\int_{\s \times \R^{+}}\pi^{+}({\rm d}\phi){\rm d}y\int_{\s }a^{+}(v)f(v,0)g(v,0){f}_{\phi}({\rm d}v,y)\nonumber\\
&=\frac{1}{\mu^{+}}\int_{\s }\pi^{+}({\rm d}\phi)\int_{\s \times \R^{+}}a^{+}(v)f(v,0)g(v,0)F_{\phi}({\rm d}v,{\rm d}y)\nonumber\\
&=\frac{1}{\mu^{+}}\int_{\s }\pi^{+}({\rm d}\phi)\bP_{0,\phi}\left[a^{+}(M^{(q),+}_{1})f(M^{(q),+}_{1},0)g(M^{(q),+}_{1},0)\right]\nonumber\\
&=\frac{1}{\mu^{+}}\int_{\s }a^{+}(v)f(v,0)g(v,0)\pi^{+}({\rm d}v).\label{prop6.1.5}
\end{align}
In the final equality we use the fact that $\pi^{+}$ is an invariant distribution for $(M^{(q),+}_{n})_{n\ge 0}$.
{ In the setting that $a^+(v) = 0$ for all $v\in \s$, the limit in \eqref{prop6.1.5} is trivial.}

Combining \eqref{prop6.1.1}, \eqref{prop6.1.2} and \eqref{prop6.1.5} we get
\begin{align}
\bP_{0,\theta}&\left[f(\Theta_{\tau^{+}_{x}-},\xi_{\tau^{+}_{x}-}-x)g(\Theta_{\tau^{+}_{x}},\xi_{\tau^{+}_{x}}-x)\right]\nonumber\\
&\to\frac{1}{\mu^{+}}\big[\int_{\s\times \R^{+}}\pi^{+}({\rm d}v){\rm d}z \ell^+(v)f(v,-z)G(v,z)+\int_{\s }\pi^{+}({\rm d}v)a^{+}(v)f(v,0)g(v,0)\nonumber\\
&+\int_{\s\times \R^{+}}\pi^{+}({\rm d}v){\rm d}y\, {\rm n}^{+}_{v}\big(\int_{0}^{\zeta}\1_{\{\epsilon_{s}\le y\}}f(\nu_{s},-y)G(\nu_{s},y){\rm d}s\big)\big]\quad \mbox{ as }x\to+\infty,\nonumber
\end{align}
which yields the first assertion of this proposition. The second and third assertion follow immediately from the above equation by setting $f\equiv1$ and $g\equiv 1$ respectively.\qed

\smallskip

In the remaining of this section we consider the nondecreasing MAP $(\bar{L}^{-1},\Theta^{+})$. The ordinate $\bar{L}^{-1}$ can be represented by
\begin{equation}\label{eq:L^{-1}_{t}}
\bar{L}^{-1}_{t}=\int_{0}^{t}\ell^+(\Theta^{+}_{s}){\rm d}s+\sum_{s\le t}\Delta \bar{L}^{-1}_{s}\quad \forall t\ge 0
\end{equation}
where $\Delta \bar{L}^{-1}_{s}=\bar{L}^{-1}_{s}-\bar{L}^{-1}_{s-}$.
Note that for any $t\ge 0$, assuming \eqref{prea2} and \eqref{prea8},
\begin{eqnarray}
\bP_{0,\pi^{+}}\left[\bar{L}^{-1}_{t}\right]
&=&\bP_{0,\pi^{+}}\left[\int_{0}^{t}\ell^+(\Theta^{+}_{s}){\rm d}s+\sum_{s\le t}\Delta \bar{L}^{-1}_{s}\right]\nonumber\\
&=&\bP_{0,\pi^{+}}\left[\int_{0}^{t}\left(\ell^+(\Theta^{+}_{s})+{\rm n}^{+}_{\Theta^{+}_{s}}(\zeta)\right){\rm d}s\right]\nonumber\\
&=&\int_{0}^{t}\bP_{0,\pi^{+}}\left[\left(\ell^+(\Theta^{+}_{s})+{\rm n}^{+}_{\Theta^{+}_{s}}(\zeta)\right)\right]{\rm d}s\nonumber\\
&=&t\int_{\s}\left(\ell^+(\theta)+{\rm n}^{+}_{\theta}(\zeta)\right)\pi^{+}({\rm d}\theta)=tc_{\pi^{+}}.\nonumber
\end{eqnarray}
In the last equality we use the fact that $\pi^{+}$ is an invariant distribution for $(\Theta^{+}_{t})_{t\ge 0}$. If we consider the Markov renewal process
$(M^{(q),+}_{n},N^{(q),+}_{n})_{n\ge 0}$,
then we have
\begin{equation}\label{eq:mean N1}
\bP_{0,\pi^{+}}\left[N^{(q),+}_{1}\right]=\bP_{0,\pi+}\left[\bar{L}^{-1}_{\mathbf{e}_{q}}\right]=\int_{0}^{+\infty}q\rme^{-q t}\bP_{0,\pi^{+}}[\bar{L}^{-1}_{t}]{\rm d}t=\frac{1}{q}c_{\pi^{+}}.
\end{equation}
For every $\theta\in\s$, define
$$W^{+}_{\theta}({\rm d}v,{\rm d}r):=\bP_{0,\theta}\left[\int_{0}^{\bar{L}_{\infty}}\1_{\{\Theta^{+}_{s}\in {\rm d}v,\ \bar{L}^{-1}_{s}\in {\rm d}r\}}{\rm d}s\right]$$
and
$G_{\theta}({\rm d}v,{\rm d}r):=\bP_{0,\theta}\left(M^{(q),+}_{1}\in {\rm d}v, N^{(q),+}_{1}\in {\rm d}r\right).$
Let $G^{*0}_{\theta}({\rm d}v,{\rm d}r):=\delta_{\theta}({\rm d}v)\delta_{0}({\rm d}r)$ and $G^{*n}_{\theta}$ be the $n$th convolution of $G_{\theta}$ for $n\ge 1$.
In view of \eqref{eq:mean N1}, under the assumptions of Proposition \ref{prop:stationary overshoots}, it follows by \cite[Theorem 2.1]{Alsmeyer1994} that
\begin{equation}\label{eq:limit of Gn}
\lim_{t\to+\infty}\int_{\s\times [0,t]}g(v,t-r)\sum_{n=0}^{+\infty}G^{*n}_{\theta}({\rm d}v,{\rm d}r)=\frac{q}{c_{\pi^{+}}}\int_{\s\times \R^{+}}g(v,r)
\pi^{+}({\rm d}v){\rm d}r
\end{equation}
for every $\theta\in\s$ and every measurable function $g\in\mathcal{M}$.
By applying similar calculations to $W^{+}_{\theta}({\rm d}v,{\rm d}r)$ as we did to $U^{+}_{\theta}({\rm d}v,{\rm d}z)$, we can show that $q W^{+}_{\theta}({\rm d}v,{\rm d}r)$ is equal to
$\sum_{n=1}^{+\infty}G^{*n}_{\theta}({\rm d}v,{\rm d}r)$. Hence by \eqref{eq:limit of Gn} we have
\begin{equation}\label{eq:limit of W+}
\lim_{t\to+\infty}\int_{\s\times [0,t]}g(v,t-r)W^{+}_{\theta}({\rm d}v,{\rm d}r)=\frac{1}{c_{\pi^{+}}}\int_{\s\times \R^{+}}g(v,r)\pi^{+}({\rm d}v){\rm d}r.
\end{equation}

\smallskip

\begin{lemma}\label{lem:w+}\mbox{}
\begin{itemize}
\item[(i)] The nondecreasing MAP $(\bar{L}^{-1},\Theta^{+})$ has a L\'{e}vy system
$(H^{+},N^{+})$ where $H^{+}_{t}=t\wedge \zeta^{+}$
and
$N^{+}(\theta,{\rm d}v,{\rm d}r):=\Gamma^{+}(\theta,{\rm d}v, {\rm d}r, [0,\infty))$
is a kernel from $\s$ to $\s\times \R^{+}$.

\item[(ii)] For $r>0$, define
$$\bar{d}_{r}:=\inf\{s>r:\ \bar{\xi}_{s}=\xi_{s}\}.$$
Then for every $\theta\in\s$, $\1_{\{\ell^{+}(v)>0\}}W^{+}_{\theta}({\rm d}v,{\rm d}r)$ has a
{kernel}
$w^{+}_{\theta}({\rm d}v,r)$ with respect to the Lebesgue measure ${\rm d}r$ such that
$$\bP_{0,\theta}\left[f(\Theta_{r});\bar{d}_{r}=r\right]=\int_{\s}f(v)\ell^+(v)w^{+}_{\theta}({\rm d}v,r)$$
for every nonnegative measurable function $f:\s\to \R^{+}$ and almost every $r>0$.
Moreover, for every $\theta\in\s$ and every bounded continuous function $h:\s\to\R$, the function $r\mapsto \bP_{0,\theta}\left[h(\Theta_{r});\bar{d}_{r}=r\right]$
is lower semi-continuous on $(0,+\infty)$.
\end{itemize}
\end{lemma}

\proof The claim in (i) follows by taking marginals in Proposition~\ref{prop:Pi+}.
\smallskip

(ii) Since $t\mapsto \bar{L}_{t}$ is a nondecreasing right continuous process, we have $\bar{L}_{r}=\inf\{s>0:\ \bar{L}^{-1}_{s}>r\}$ for every $r>0$.
We also note that $\bar{L}^{-1}(\bar{L}_{r})=\inf\{s>r:\ \bar{\xi}_{s}=\xi_{s}\}=\bar{d}_{r}$.
In view of this, (i) and \eqref{eq:L^{-1}_{t}}, we can apply Proposition \ref{prop:on creeping} to the process $(\bar{L}^{-1},\Theta^{+})$ and deduce that
$\1_{\{\ell^{+}(v)>0\}}W^{+}_{\theta}({\rm d}v,{\rm d}r)$
has a
{kernel}
$w^{+}_{\theta}({\rm d}v,r)$ with respect to the Lebesgue measure ${\rm d}r$ such that
\begin{equation}\label{lemw+.1}
\bP_{0,\theta}\left[f(\Theta_{r});\bar{d}_{r}=r\right]=\bP_{0,\theta}\left[f(\Theta^{+}_{\bar{L}_{r}});\bar{L}^{-1}(\bar{L}_{r})=r\right]=\int_{\s}f(v)\ell^+(v)w^{+}_{\theta}({\rm d}v,r)
\end{equation}
for almost every $r>0$ and every nonnegative measurable function $f:\s\to \R^{+}$.
Now take an arbitrary bounded continuous function $h:\s\to\R$. We have
$$\bP_{0,\theta}\left[h(\Theta_{r});\bar{d}_{r}=r\right]=\bP_{0,\theta}\left[h(\Theta_{r})\right]-\bP_{0,\theta}\left[h(\Theta_{r});\bar{d}_{r}>r\right].$$
It is easy to see that $r\mapsto\bP_{0,\theta}\left[h(\Theta_{r})\right]$ is right continuous on $[0,+\infty)$ since $\Theta$ is a right continuous process. We only need to show that $r\mapsto \bP_{0,\theta}\left[h(\Theta_{r});\bar{d}_{r}>r\right]$ is upper semi-continuous on $(0,+\infty)$. Take an arbitrary sequence $r_{n}\downarrow r\in (0,+\infty)$. Note that, for any $s>0$, $\bar{d}_{s}>s$ if and only if $s\in \cup_{g_{i}\in\bar{G}}[g_{i},d_{i})$. Hence $\{\bar{d}_{r_{n}}>r_{n}\mbox{ i.o.}\}=\{r_{n}\in \cup_{g_{i}\in\bar{G}}[g_{i},d_{i})\mbox{ i.o.}\}\subset \{r\in \cup_{g_{i}\in\bar{G}}[g_{i},d_{i})\}=\{\bar{d}_{r}>r\}$. It follows that $\limsup_{n\to+\infty}\1_{\{\bar{d}_{r_{n}}>r_{n}\}}=\1_{\{\bar{d}_{r_{n}}>r_{n}\mbox{ i.o.}\}}\le \1_{\{\bar{d}_{r}>r\}}$. Thus by the reverse Fatou's lemma, $\bP_{0,\theta}\left[h(\Theta_{r});\bar{d}_{r}>r\right]\ge \limsup_{n\to+\infty}\bP_{0,\theta}\left[h(\Theta_{r_{n}});\bar{d}_{r_{n}}>r_{n}\right]$. We complete the proof.
\qed

\smallskip

\begin{lemma}\label{lemA.7}
Suppose that $((\xi,\Theta),\bP)$ and $((\xi,\Theta),\hat{\bP})$ are a pair of upwards regular MAPs for which condition \eqref{condi:weak reversability} is satisfied.
Under the assumptions of Proposition \ref{prop:stationary overshoots}, we have
\begin{itemize}
\item[(i)] $\int_{\s}\ell^{+}(\theta)\pi^{+}({\rm d}\theta)=0$, and $\int_{0}^{+\infty}\ell^{+}(\Theta_{s}){\rm d}\bar{L}_{s}=0,$ $\bP_{0,\pi}$-a.s.
\item[(ii)] For every $y<0$,
\begin{equation}\label{lemA.7.1}
\hat{H}^{+}_{\theta}(y)\pi({\rm d}\theta)=\frac{1}{c_{\pi^{+}}}\int_{\mathcal{S}}\pi^{+}({\rm d}\phi)
{\rm n}^{+}_{\phi}\left(\int_{0}^{\zeta}\1_{\{\epsilon_{r}\le -y,\ \nu_{r}\in {\rm d}\theta\}}{\rm d}r\right)
\end{equation}
where $\hat{H}^{+}_{\theta}(y)=\hat{\bP}_{y,\theta}\left(\tau^{+}_{0}=+\infty\right)$, and
\begin{equation}\label{lemA.7.2}
\frac{\hat{{\rm n}}^{+}_{\theta}(\zeta=+\infty)}{\ell^+(\theta)+{\rm n}^{+}_{\theta}(\zeta)}\hat{U}^{+}_{\pi}({\rm d}\theta,\R^{+})=\frac{1}{c_{\pi^{+}}}\pi^{+}({\rm d}\theta).
\end{equation}
\end{itemize}
\end{lemma}

\proof
(i)
By \eqref{def:pi+},
we have
\begin{equation}\label{lem6.3.1}
\int_{\s}\ell^{+}(\theta)\pi^{+}({\rm d}\theta)=\frac{1}{\bP_{0,\pi}\left[\bar{L}_{1}\right]}\bP_{0,\pi}\left[\int_{0}^{1}\ell^{+}(\Theta_{s}){\rm d}\bar{L}_{s}\right].
\end{equation}
We note that by \eqref{M2} and Fubini's theorem,
$$\bP_{0,\pi}\left[\int_{0}^{+\infty}\ell^{+}(\Theta_{s}){\rm d}\bar{L}_{s}\right]=\bP_{0,\pi}\left[\int_{0}^{+\infty}\1_{\{s\in\bar{M}\}}{\rm d}s\right]=\int_{0}^{+\infty}\bP_{0,\pi}\left(s\in\bar{M}\right){\rm d}s.$$
By Proposition \ref{prop:equal dist}, we have for any $s>0$,
$$\bP_{0,\pi}\left(s\in\bar{M}\right)=\bP_{0,\pi}\left(\bar{\xi}_{s}-\xi_{s}=0\right)=\hat{\bP}_{0,\pi}\left(\bar{\xi}_{s}=0\right)\le \hat{\bP}_{0,\pi}\left(\tau^{+}_{0}\ge s\right)=0.$$
The last equality is because $((\xi,\Theta),\hat{\bP})$ is upwards regular. It follows that
\begin{equation}\label{lem6.3.3}
\bP_{0,\pi}\left[\int_{0}^{+\infty}\ell^{+}(\Theta_{s}){\rm d}\bar{L}_{s}\right]=0,
 \end{equation}
 and hence by \eqref{lem6.3.1} $\int_{\s}\ell^{+}(\theta)\pi^{+}({\rm d}\theta)=0$.

\smallskip

(ii) First we claim that
\begin{equation}\label{lem6.3.2}
\bP_{0,\pi}\left(\bar{d}_{r}=r\right)=0\quad\forall r>0.
\end{equation}
In fact, by Lemma \ref{lem:w+}(ii) and \eqref{lem6.3.3}, we have
\begin{eqnarray}
\int_{0}^{+\infty}\bP_{0,\pi}\left(\bar{d}_{r}=r\right){\rm d}r
&=&\int_{0}^{+\infty}{\rm d}r\int_{\s}\ell^{+}(v)w^{+}_{\pi}({\rm d}v,r)\nonumber\\
&=&\int_{0}^{+\infty}\int_{\s}\ell^{+}(v)W^{+}_{\pi}({\rm d}v,{\rm d}r)\nonumber\\
&=&\bP_{0,\pi}\left[\int_{0}^{+\infty}\ell^{+}(\Theta_{s}){\rm d}\bar{L}_{s}\right]=0.
\end{eqnarray}
Thus $\bP_{0,\pi}\left(\bar{d}_{r}=r\right)=0$ for almost every $r>0$, and hence for every $r>0$ since $r\mapsto \bP_{0,\pi}\left(\bar{d}_{r}=r\right)$ is lower semi-continuous on $(0,+\infty)$.
By Proposition \ref{prop:equal dist} we have
\begin{equation}\label{lemA.7.3}
\hat{\bP}_{0,\pi}\left[g(\Theta_{0});\bar{\xi}_{t}\le -y\right]=\bP_{0,\pi}\left[g(\Theta_{t});\bar{\xi}_{t}-\xi_{t}\le -y\right]
\end{equation}
for every $y<0$, $t\ge 0$ and every bounded measurable function $g:\s\to\R$. It follows by the bouned convergence theorem that
\begin{eqnarray}
\hat{\bP}_{0,\pi}\left[g(\Theta_{0});\bar{\xi}_{t}\le -y\right]&=&\hat{\bP}_{0,\pi}\left[g(\Theta_{0});\tau^{+}_{-y}>t\right]=\int_{\mathcal{S}}\pi({\rm d}\theta)g(\theta)\hat{\bP}_{0,\theta}\left(\tau^{+}_{-y}>t\right)\nonumber\\
&\to&\int_{\s}\pi({\rm d}\theta)g(\theta)\hat{\bP}_{0,\theta}\left(\tau^{+}_{-y}=+\infty\right)=\int_{\mathcal{S}}\pi({\rm d}\theta)g(\theta)\hat{H}^{+}_{\theta}(y),\label{lemA.7.4}
\end{eqnarray}
as $t\to+\infty$.
On the other hand, we have by \eqref{lem6.3.2}
\begin{equation}
\bP_{0,\pi}\left[g(\Theta_{t});\bar{\xi}_{t}-\xi_{t}\le -y\right]
=\bP_{0,\pi}\left[g(\Theta_{t});\bar{\xi}_{t}-\xi_{t}\le -y,\bar{d}_{t}>t\right]\quad\forall t>0.
\end{equation}
We note that $\bar{d}_{t}>t$ if and only if $t\in \cup_{g_{i}\in \bar{G}}[g_{i},d_{i})$. Hence by \eqref{eq:last exit} the above expectation equals
\begin{align}
\bP_{0,\pi}&\left[g(\Theta_{t});\bar{\xi}_{t}-\xi_{t}\le -y,t\in \cup_{g_{i}\in \bar{G}}[g_{i},d_{i})\right]\nonumber\\
&=\bP_{0,\pi}\left[\int_{0}^{t}{\rm n}^{+}_{\Theta_{s}}\left(g(\nu_{t-s})\1_{\{\epsilon_{t-s}\le -y,t-s<\zeta\}}\right){\rm d}\bar{L}_{s}\right]\nonumber\\
&=\bP_{0,\pi}\left[\int_{0}^{+\infty}\1_{\{\bar{L}^{-1}_{u}\le t\}}{\rm n}^{+}_{\Theta^{+}_{u}}\left(g(\nu_{t-\bar{L}^{-1}_{u}})\1_{\{\epsilon_{t-\bar{L}^{-1}_{u}}\le -y,t-\bar{L}^{-1}_{u}<\zeta\}}\right){\rm d}u\right]\nonumber\\
&=\int_{\s\times [0,t]}W^{+}_{\pi}({\rm d}v,{\rm d}r){\rm n}^{+}_{v}\left(g(\nu_{t-r})\1_{\{\epsilon_{t-r}\le -y,t-r<\zeta\}}\right).\label{lemA.7.12}
\end{align}
By \eqref{eq:limit of W+}, the integral in the right converges as $t\to+\infty$ to
\begin{equation}\nonumber
\frac{1}{c_{\pi^{+}}}\int_{\s\times \R^{+}}\pi^{+}({\rm d}v){\rm d}r\,{\rm n}^{+}_{v}\left(g(\nu_{r})\1_{\{\epsilon_{r}\le -y,r<\zeta\}}\right)
=\frac{1}{c_{\pi^{+}}}\int_{\s}\pi^{+}({\rm d}v){\rm n}^{+}_{v}\left(\int_{0}^{\zeta}g(\nu_{r})\1_{\{\epsilon_{r}\le -y\}}{\rm d}r\right).
\end{equation}
Combining this and \eqref{lemA.7.3}-\eqref{lemA.7.12} we get that
$$\int_{\mathcal{S}}\pi({\rm d}\theta)g(\theta)\hat{H}^{+}_{\theta}(y)=\frac{1}{c_{\pi^{+}}}\int_{\s}\pi^{+}({\rm d}\theta) {\rm n}^{+}_{\theta}\left(\int_{0}^{\zeta}g(\nu_{r})\1_{\{\epsilon_{r}\le -y\}}{\rm d}r\right)$$
for any bounded measurable function $g:\s\to\R$,
which in turn yields \eqref{lemA.7.1}.

\smallskip

Next we prove \eqref{lemA.7.2}. It follows by Proposition \ref{prop:equal dist} that
\begin{equation}
\hat{\bP}_{0,\pi}\left[g(\bar{\Theta}_{t})\right]=\bP_{0,\pi}\left[g(\bar{\Theta}_{t})\right]\quad \forall t\ge 0\label{lemA.7.6}
\end{equation}
for any bounded measurable function $g:\s\to \R$.
Similarly by \eqref{lem6.3.2} and \eqref{eq:last exit} we have
\begin{eqnarray}
\bP_{0,\pi}\left[g(\bar{\Theta}_{t})\right]
&=&\bP_{0,\pi}\left[g(\bar{\Theta}_{t});t\in \cup_{g_{i}\in \bar{G}}[g_{i},d_{i})\right]\nonumber\\
&=&\bP_{0,\pi}\left[\sum_{g_{i}\in\bar{G}}g(\Theta_{g_{i}})\1_{\{g_{i}\le t<d_{i}\}}\right]\nonumber\\
&=&\int_{\s\times [0,t]}W^{+}_{\pi}({\rm d}v,{\rm d}r)g(v){\rm n}^{+}_{v}(t-r<\zeta).\nonumber
\end{eqnarray}
By \eqref{eq:limit of W+}, we get
$$\lim_{t\to+\infty}\bP_{0,\pi}\left[g(\bar{\Theta}_{t})\right]=\frac{1}{c_{\pi^{+}}}\int_{\s}g(\theta){\rm n}^{+}_{\theta}(\zeta)\pi^{+}({\rm d}\theta).$$
It follows by this, \eqref{lemA.7.6}, the bounded convergence theorem and Lemma \ref{lemA.7}(i) that
\begin{eqnarray}\label{lemA.7.11}
\hat{\bP}_{0,\pi}\left[g(\bar{\Theta}_{\mathbf{e}_{q}})\right]&=&\bP_{0,\pi}\left[g(\bar{\Theta}_{\mathbf{e}_{q}})\right]\nonumber\\
&=&\int_{0}^{+\infty}\rme^{-s}\bP_{0,\pi}\left[g(\bar{\Theta}_{s/q})\right]{\rm d}s\nonumber\\
&\to& \frac{1}{c_{\pi^{+}}}\int_{\s}g(\theta){\rm n}^{+}_{\theta}(\zeta)\pi^{+}({\rm d}\theta)=\frac{1}{c_{\pi^{+}}}\int_{\s}g(\theta)\left(\ell^{+}(\theta)+{\rm n}^{+}_{\theta}(\zeta)\right)\pi^{+}({\rm d}\theta)
\end{eqnarray}
as $q\to 0+$. Let $\mathcal{C}$ denote the set of nonnegative bounded measurable functions $h:\s\to\R^{+}$ such that $\theta\mapsto h(\theta)a^{+}(\theta)/\left(\ell^+(\theta)+{\rm n}^{+}_{\theta}(\zeta)\right)$ is a bounded function. On the one hand, by \eqref{lemA.7.11} we have
\begin{equation}\label{lemA.7.10}
\hat{\bP}_{0,\pi}\left[\frac{h(\bar{\Theta}_{\mathbf{e}_{q}})a^{+}(\bar{\Theta}_{\mathbf{e}_{q}})}{\ell^+(\bar{\Theta}_{\mathbf{e}_{q}})+{\rm n}^{+}_{\bar{\Theta}_{\mathbf{e}_{q}}}(\zeta)}\right]
\to \frac{1}{c_{\pi^{+}}}\int_{\s}h(\theta)a^{+}(\theta)\pi^{+}({\rm d}\theta)\quad\mbox{ as }q\to 0+
\end{equation}
for any $h\in\mathcal{C}$. On the other hand, by Proposition \ref{prop:wiener_hopf} we have
\begin{equation}\label{lemA.7.7}
\hat{\bP}_{0,\pi}\left[\frac{h(\bar{\Theta}_{\mathbf{e}_{q}})a^{+}(\bar{\Theta}_{\mathbf{e}_{q}})}{\ell^+(\bar{\Theta}_{\mathbf{e}_{q}})+{\rm n}^{+}_{\bar{\Theta}_{\mathbf{e}_{q}}}(\zeta)}\right]
=\int_{\s\times \R^{+}}\hat{W}^{+}_{\pi}({\rm d}v,{\rm d}r)\rme^{-q r}\frac{h(v)a^{+}(v)}{\ell^+(v)+{\rm n}^{+}_{v}(\zeta)}\left(q \hat\ell^+(v)+\hat{{\rm n}}^{+}_{v}(1-\rme^{-q\zeta})\right).
\end{equation}
If we can show that
\begin{equation}
\lim_{q\to0+}\hat{\bP}_{0,\pi}\left[\frac{h(\bar{\Theta}_{\mathbf{e}_{q}})a^{+}(\bar{\Theta}_{\mathbf{e}_{q}})}{\ell^+(\bar{\Theta}_{\mathbf{e}_{q}})+{\rm n}^{+}_{\bar{\Theta}_{\mathbf{e}_{q}}}(\zeta)}\right]=\int_{\s}\hat{W}^{+}_{\pi}({\rm d}v,\R^{+})\frac{h(v)a^{+}(v)}{\ell^+(v)+{\rm n}^{+}_{v}(\zeta)}\hat{{\rm n}}^{+}_{v}(\zeta=+\infty),\quad\forall h\in\mathcal{C},\label{lemA.7.8}
\end{equation}
then by \eqref{lemA.7.10} and the fact that $\hat{W}^{+}_{\pi}({\rm d}v,\R^{+})=\hat{U}^{+}_{\pi}({\rm d}v,\R^{+})=\hat{\bP}_{0,\pi}\left[\int_{0}^{\bar{L}_{\infty}}\1_{\{\hat{\Theta}^{+}_{s}\in {\rm d}v\}}{\rm d}s\right]$  we get
\begin{equation}
\int_{\s}\hat{U}^{+}_{\pi}({\rm d}v,\R^{+})\frac{h(v)a^{+}(v)}{\ell^+(v)+{\rm n}^{+}_{v}(\zeta)}\hat{{\rm n}}^{+}_{v}(\zeta=+\infty)=\frac{1}{c_{\pi^{+}}}\int_{\s}h(\theta)a^{+}(\theta)\pi^{+}({\rm d}\theta)
\quad\forall h\in\mathcal{C}.\label{eq:6.32}
\end{equation}
Note that for any $q\in (0,1]$ the integrand on the right hand side of \eqref{lemA.7.7} is bounded from above by
$$\|h\|_{\infty}\frac{a^{+}(v)}{\ell^+(v)+{\rm n}^{+}_{v}(\zeta)}\left(\hat\ell^+(v)+\hat{{\rm n}}^{+}_{v}(1-\rme^{-\zeta})\right).$$
Hence to prove \eqref{lemA.7.8} it suffices to prove
\begin{equation}\label{lemA.7.9}
\int_{\s}\hat{W}^{+}_{\pi}({\rm d}v,\R^{+})\frac{a^{+}(v)}{\ell^+(v)+{\rm n}^{+}_{v}(\zeta)}\left(\hat\ell^+(v)+\hat{{\rm n}}^{+}_{v}(1-\rme^{-\zeta})\right)<+\infty.
\end{equation}
By Proposition \ref{prop:equal dist} and Proposition \ref{prop:wiener_hopf} the above integral is equal to
\begin{eqnarray}
\hat{\bP}_{0,\pi}\left[\rme^{  \bar{g}_{e_{1}}}\frac{a^{+}(\bar{\Theta}_{e_{1}})}{\ell^+(\bar{\Theta}_{e_{1}})+{\rm n}^{+}_{\bar{\Theta}_{e_{1}}}(\zeta)}\right]
&=&\bP_{0,\pi}\left[\rme^{(e_{1}-\bar{g}_{e_{1}})}\frac{a^{+}(\bar{\Theta}_{e_{1}})}{\ell^+(\bar{\Theta}_{e_{1}})+{\rm n}^{+}_{\bar{\Theta}_{e_{1}}}(\zeta)}\right]\nonumber\\
&=& \int_{\R^{+}\times \s\times \R^{+}}\rme^{-r}a^{+}(v)V^{+}_{\pi}({\rm d}r,{\rm d}v,{\rm d}z)\nonumber\\
&=& \bP_{0,\pi}\left[\int_{0}^{+\infty}\rme^{-\bar{L}^{-1}_{s}}a^{+}(\Theta^{+}_{s}){\rm d}s\right]\nonumber\\
&=& \bP_{0,\pi}\left[\int_{0}^{\bar{L}_{e_{1}}}a^{+}(\Theta^{+}_{s}){\rm d}s\right].\nonumber
\end{eqnarray}
The finiteness of the final expectation is implied by the finiteness of $\bP_{0,\pi^{+}}\left[\xi^{+}_{1}\right]$. Indeed, by \eqref{def:pi+} and Markov property
\begin{eqnarray}
\bP_{0,\pi^{+}}\left[\xi^{+}_{1}\right]
&=&\frac{1}{\bP_{0,\pi}\left[\bar{L}_{1}\right]}\bP_{0,\pi}\left[\int_{0}^{1}\bP_{0,\Theta_{s}}\left[\xi^{+}_{1}\right]{\rm d} \bar{L}_{s}\right]\nonumber\\
&=&\frac{1}{\bP_{0,\pi}\left[\bar{L}_{1}\right]}\bP_{0,\pi}\left[\int_{0}^{\bar{L}_{1}}\bP_{0,\Theta^{+}_{s}}\left[\xi^{+}_{1}\right]{\rm d}s\right]\nonumber\\
&=&\frac{1}{\bP_{0,\pi}\left[\bar{L}_{1}\right]}\bP_{0,\pi}\left[\int_{0}^{\bar{L}_{1}}\left(\xi^{+}_{s+1}-\xi^{+}_{s}\right){\rm d}s\right].\nonumber
\end{eqnarray}
Since the continuous part of $\xi^{+}_{s+1}-\xi^{+}_{s}$ is $\int_{s}^{s+1}a^{+}(\Theta^{+}_{r}){\rm d}r$, we get by Fubini's theorem that
$$+\infty>\bP_{0,\pi^{+}}\left[\xi^{+}_{1}\right]\bP_{0,\pi}\left[\bar{L}_{1}\right]\ge\bP_{0,\pi}\left[\int_{0}^{\bar{L}_{1}}{\rm d}s
\int_{s}^{s+1}a^{+}(\Theta^{+}_{r}){\rm d}r\right]=\bP_{0,\pi}\left[\int_{0}^{\bar{L}_{1}+1}a^{+}(\Theta^{+}_{r}){\rm d}r\right].$$
By writing $\bP_{0,\pi}\left[\int_{0}^{\bar{L}_{s}}a^{+}(\Theta^{+}_{r}){\rm d}r\right]=\bP_{0,\pi}\left[\int_{0}^{s}a^{+}(\Theta_{r}){\rm d}\bar{L}_{r}\right]$, one can show that $s\mapsto \bP_{0,\pi}\left[\int_{0}^{\bar{L}_{s}}a^{+}(\Theta^{+}_{r}){\rm d}r\right]$ is a subadditive and locally bounded nonnegative function, which implies the finiteness of $\bP_{0,\pi}\left[\int_{0}^{\bar{L}_{e_{1}}}a^{+}(\Theta^{+}_{s}){\rm d}s\right]$.

We deduce therefrom that (\ref{lemA.7.8}) and hence (\ref{eq:6.32}) hold for every $h\in\mathcal{C}.$ Now, for a general nonnegative measurable function $h$, one can always find an nondecreasing sequence of functions $h_{n}\in\mathcal{C}$ such that $h_{n}\to h$ in the pointwise sense. Using this and the monotone convergence theorem, one can show that (\ref{eq:6.32}) holds for any nonnegative function $h$. The identity \eqref{lemA.7.2} follows immediately.\qed

\smallskip

\begin{proposition}\label{prop:rhooplus}
{Suppose that the assumptions of Lemma \ref{lemA.7} hold.}
Then the stationary distribution $\rho^{\oplus}({\rm d}y,{\rm d}v)$ given in Proposition \ref{prop:stationary overshoots} can be represented by
$$\rho^{\oplus}({\rm d}y,{\rm d}v)=\rho^{\oplus}_{1}({\rm d}y,{\rm d}v)+\rho^{\oplus}_{2}({\rm d}y,{\rm d}v)$$
where
$$\rho^{\oplus}_{1}({\rm d}y,{\rm d}v):=\frac{c_{\pi^{+}}}{\mu^{+}}\1_{\{y<0\}}\bar{\Pi}_{v}(-y)\hat{H}^{+}_{v}(y){\rm d}y\,\pi({\rm d}v),$$
and
$$\rho^{\oplus}_{2}({\rm d}y,{\rm d}v):=\frac{c_{\pi^{+}}}{\mu^{+}}\,\frac{a^{+}(v)\hat{{\rm n}}^{+}_{v}(\zeta=+\infty)}{\ell^+(v)+{\rm n}^{+}_{v}(\zeta)}\delta_{0}({\rm d}y)\hat{U}^{+}_{\pi}({\rm d}v,\R^{+}).$$
\end{proposition}

\smallskip

\part{Main results and their proofs}

\section{Assumptions and main results}\label{sec:main results}

{
Suppose $E$ is a locally compact separable metric space and $\mathcal{B}(E)$ is the minimal Borel field in $E$ containing all the open sets. Let $E_{\partial}=E\cup \{\partial\}$ (where $\partial\not\in E$) be its one-point compactification and $C_{0}(E)$ be the class of all continuous functions on $E_{\partial}$ vanishing at $\partial$. Suppose $(Y,\{\mathbf{Q}_{y}:y\in E\})$ is a Markov process on $E$ with lifetime $\zeta$, whose transition semigroup $(Q_{t})_{t\ge 0}$ is given by
$$Q_{t}f(y)=\mathbf{Q}_{y}\left[f(Y_{t}),t<\zeta\right],\quad Q_{0}f(y)=\mathbf{Q}_{y}\left[f(Y_{0})\right]=f(y)$$
for $t>0$, $y\in E$ and nonnegative $\mathcal{B}(E)$-measurable function $f$. In general, $(Q_{t})_{t\ge 0}$ is a subMarkovian semigroup on $(E,\mathcal{B}(E))$. It can be extended to be strictly Markovian on $(E_{\partial},\mathcal{B}(E_{\partial}))$ by setting additionally that: for $t\ge 0$,
\begin{align*}
Q_{t}\1_{\{\partial\}}(x)&=1-Q_{t}\1_{E}(x)\quad\forall x\in E,\\
Q_{t}\1_{E}(\partial)&=0,\quad Q_{t}\1_{\{\partial\}}(\partial)=1.
\end{align*}
This extended transition semigroup naturally defines a Markov process $(Y, \{\mathbf{Q}_{y}:y\in E_{\partial}\})$ on $E_{\partial}$, where we have $\mathbf{Q}_{\partial}\left(Y_{t}=\partial \mbox{ for all }t\ge 0\right)=1$.

\begin{definition}
We say the Markov process $(Y,\{\mathbf{Q}_{y}:y\in E\})$ is a Feller process if its extended transition semigroup on $(E_{\partial},\mathcal{B}(E_{\partial}))$ has the Feller property:
\begin{itemize}
\item[(i)] $Q_{t}f\in C_{0}(E)$ for all $t>0$ and $f\in C_{0}(E)$;
\item[(ii)] $\lim_{t\to 0}\sup_{y\in E_{\partial}}|Q_{t}f(y)-f(y)|=0$ for all $f\in C_{0}(E)$.
\end{itemize}
It is known by \cite[Chaper 2]{ChungWalsh} that under (i) the condition (ii) is equivalent to the apparently weaker condition below:
\begin{itemize}
\item[(ii')] $\lim_{t\to 0}Q_{t}f(y)=f(y)$ for all $f\in C_{0}(E)$ and $y\in E_{\partial}$.
\end{itemize}
\end{definition}

Let $\mathcal{I}=\R\times \s$ and $\mathcal{I}_{\partial}=\mathcal{I}\cup\{\partial\}$ be the one-point compactification of $\mathcal{I}$. Suppose $((\xi,\Theta),\{\bP_{x,\theta}:(x,\theta)\in\mathcal{I}\})$ is a MAP on $\mathcal{I}$. It can be extended to be a Markov process on $\mathcal{I}_{\partial}$ as shown in the above argument. Recall that $\phi(x,\theta)=\theta\mathrm{e}^{x}$ for $(x,\theta)\in\mathcal{I}$. We denote $\phi(\partial)$ by $\triangle$. Let $\h_{\triangle}=\phi(\mathcal{I}_{\partial})=\h\cup\{\triangle\}$.
For every $(x,\theta)\in \mathcal{I}_{\partial}$, let $\P_{\phi(x,\theta)}$ be the law of $X$ given by the Lamperti-Kiu transform \eqref{eq:lamperti_kiu} under $\bP_{x,\theta}$. Here we assume conventionally that $\P_{\triangle}\left(X_{t}=\triangle\mbox{ for all }t\ge 0\right)=1$. Then $(X,\{\P_{z}:z\in \h_{\triangle}\})$ is a Markov process on $\h_{\triangle}$. First we give a lemma which complements the result given in Lemma \ref{lem2.11}.

\begin{lemma}\label{lemn2}
Suppose the following condition holds.
\begin{itemize}
\item[(a1)] $((\xi,\Theta),\{\bP_{x,\theta}:(x,\theta)\in \mathcal{I}\})$ is a Feller process.
\end{itemize}
Then for any $(r_{n},\theta_{n}),\ (r_{0},\theta_{0})\in \mathcal{I}_{\partial}$ with $(r_{n},\theta_{n})\to (r_{0},\theta_{0})$ as $n\to +\infty$, $(X,\P_{\phi(r_{n},\theta_{n})})$ converges to $(X,\P_{\phi(r_{0},\theta_{0})})$ in distribution under the Skorokhod topology.
\end{lemma}
\proof Fix an arbitrary sequence $\{(r_{n},\theta_{n}):n\ge 1\}\subset \mathcal{I}_{\partial}$ such that $(r_{n},\theta_{n})\to (r_{0},\theta_{0})\in \mathcal{I}_{\partial}$. In view of (a1), it follows by \cite[Theorem 4.2.5]{EthKur} that $((\xi,\Theta),\bP_{r_{n},\theta_{n}})$ converges to $((\xi,\Theta),\bP_{r_{0},\theta_{0}})$ in distribution under the Skorokhod topology. (Here we take the convention that $\bP_{\partial}\left((\xi_{t},\Theta_{t})=\partial\mbox{ for all }t\ge 0\right)=1$.) Thus by the Skorokhod representation theorem, there exist a probability space
$(\Omega^{\circ},\mathcal{F}^{\circ},\bP^{\circ})$ and couplings $(\xi^{(n)},\Theta^{(n)})$, $(\xi^{(0)},\Theta^{(0)})$ of the processes $((\xi,\Theta),\bP_{r_{n},\theta_{n}})$ and $((\xi,\Theta),\bP_{r_{0},\theta_{0}})$, respectively, such that $(\xi^{(n)},\Theta^{(n)})$ converges to $(\xi^{(0)},\Theta^{(0)})$ $\bP^{\circ}$-almost surely under the Skorokhod topology. Thus there is a subset $\Omega'\subset \Omega^{\circ}$ with $\bP^{\circ}(\Omega')=1$ such that for all $\omega\in\Omega'$, $(\xi^{(n)}(\omega),\Theta^{(n)}(\omega))$ converges to $(\xi^{(0)}(\omega),\Theta^{(0)}(\omega))$ in $\mathbb{D}_{\R\times\s}$. We fix an arbitrary $\omega\in\Omega'$. For $k\ge 0$ and $t\ge 0$, let $z^{(k)}_{t}(\omega):=\mathrm{e}^{\alpha \xi^{(k)}_{t}(\omega)}\Theta^{(k)}_{t}(\omega)$, $x^{(k)}_{t}(\omega):=\mathrm{e}^{\alpha \xi^{(k)}_{t}(\omega)}$,
$y^{(k)}_{t}(\omega):=\int_{0}^{t}x^{(k)}_{s}(\omega)ds$, and $y^{(k),-1}_{t}(\omega):=\inf\{s>0:\ y^{(k)}_{s}(\omega)>t\}$.  It follows by \cite[Theorem 3.1]{Whitt2} that $z^{(n)}(\omega)$ converges to $z^{(0)}(\omega)$ in $\mathbb{D}_{\R^{d}}$ and $x^{(n)}(\omega)$ convergence to $x^{(0)}(\omega)$ in $\mathbb{D}_{\R}$. Hence by \cite[Proposition 3.5.3(c)]{EthKur}, for every $T\in (0,+\infty)$, there is a sequence of strictly increasing continuous functions $\{\lambda_{n}:[0,T]\to \R^{+}\}$ with $\lambda_{n}(0)=0$ such that
\begin{equation}\label{n2.1}
\lim_{n\to+\infty}\sup_{t\in [0,T]}\left(|x^{(n)}_{t}(\omega)-x^{(0)}_{\lambda_{n}(t)}(\omega)|\vee|\lambda_{n}(t)-t|\right)=0.
\end{equation}
We observe that
\begin{align*}
|y^{(n)}_{t}(\omega)-y^{(0)}_{\lambda_{n}(t)}(\omega)|&=\left|\int_{0}^{t}x^{(n)}_{s}(\omega){\rm d}s-\int_{0}^{\lambda_{n}(t)}x^{(0)}_{s}(\omega){\rm d}s\right|\\
&\le \left|\int_{0}^{t}x^{(n)}_{s}(\omega)-x^{(0)}_{\lambda_{n}(s)}(\omega){\rm d}s\right|
+\left|\int_{0}^{t}x^{(0)}_{\lambda_{n}(s)}(\omega)-x^{(0)}_{s}(\omega){\rm d}s\right|\\
&\quad +\left|\int_{0}^{t}x^{(0)}_{s}(\omega){\rm d}s-\int_{0}^{\lambda_{n}(t)}x^{(0)}_{s}(\omega){\rm d}s\right|\\
&\le \int_{0}^{t}\left|x^{(n)}_{s}(\omega)-x^{(0)}_{\lambda_{n}(s)}(\omega)\right|{\rm d}s
+\int_{0}^{t}\left|x^{(0)}_{\lambda_{n}(s)}(\omega)-x^{(0)}_{s}(\omega)\right|{\rm d}s+\int_{t\wedge \lambda_{n}(t)}^{t\vee \lambda_{n}(t)}|x^{(0)}_{s}(\omega)|{\rm d}s.
\end{align*}
Hence for $T\in (0,+\infty)$,
\begin{align}
\sup_{t\in [0,T]}|y^{(n)}_{t}(\omega)-y^{(0)}_{\lambda_{n}(t)}(\omega)|&\le T\sup_{t\in [0,T]}|x^{(n)}_{t}(\omega)-x^{(0)}_{\lambda_{n}(t)}(\omega)|+
\int_{0}^{T}\left|x^{(0)}_{\lambda_{n}(s)}(\omega)-x^{(0)}_{s}(\omega)\right|{\rm d}s\notag\\
&\quad+\sup_{t\in [0,T]}|\lambda_{n}(t)-t|\cdot\sup_{s\in [0,T+|\lambda_{n}(T)-T|]}|x^{(0)}_{s}(\omega)|.\label{lem6.2.4}
\end{align}
Immediately by \eqref{n2.1} the first and third terms on the right hand side converge to $0$ as $n\to +\infty$. Since $\lambda_{n}(s)\to s$ for every $s\in [0,T]$, one has
$x^{(0)}_{\lambda_{n}(s)}(\omega)-x^{(0)}_{s}(\omega)\to 0$ at every continuous point $s\in [0,T]$ of the function $t\mapsto x^{(0)}_{t}(\omega)$. Thus by the right continuity of $t\mapsto x^{(0)}_{t}(\omega)$ and the bounded convergence theorem the second term on the right hand side of \eqref{lem6.2.4} converges to $0$. Therefore we get $\sup_{t\in [0,T]}|y^{(n)}_{t}(\omega)-y^{(0)}_{\lambda_{n}(t)}(\omega)|\to 0$, and again by \cite[Proposition 3.5.3(c)]{EthKur} $y^{(n)}(\omega)$ converges to $y^{(0)}(\omega)$ in $\mathbb{D}_{\R}$. It then follows by Theorem 7.2 and Theorem 3.1 in \cite{Whitt2} that $y^{(n),-1}(\omega)$ converges to $y^{(0),-1}(\omega)$ in $\mathbb{D}_{R}$ and  $z^{(n)}\circ y^{(n),-1}(\omega)$ converges to $z^{(0)}\circ y^{(0),-1}(\omega)$ in $\mathbb{D}_{\R^{d}}$. The above argument shows that $z^{(n)}\circ y^{(n),-1}$ converges to $z^{(0)}\circ y^{(0),-1}$ $\bP^{\circ}$-almost surely under the Skorokhod topology. We note that for $k\ge 0$ $z^{(k)}\circ y^{(k),-1}$ corresponds to the process $(\xi^{(k)},\Theta^{(k)})$ via the Lamperti-Kiu transform, and thus $\left(z^{(k)}\circ y^{(k),-1},\bP^{\circ}\right)$ is equal in law to $(X,\P_{\phi(r_{k},\theta_{k})})$. Hence we prove that $(X,\P_{\phi(r_{n},\theta_{n})})$ converges to $(X,\P_{\phi(r_{0},\theta_{0})})$ in distribution under the Skorokhod topology.\qed

\smallskip
}

{In what follows, we assume that}
$(X,\{\P_{z}, z\in \h\})$ is an $\h$-valued ssMp and $((\xi,\Theta),\bP)$ is the corresponding MAP via the Lamperti-Kiu transform, for which its L\'evy system $(H,\Pi)$  satisfies $H_t = t$ until killing.
{We assume condition (a1) and the following additional conditions hold.}

\begin{itemize}

\item[(a2)] The modulator of $((\xi,\Theta),\bP)$ is a positive recurrent process having an invariant distribution $\pi$ which is fully supported on $\s$. The continuous part of $\xi^{+}$ of $((\xi^{+},\Theta^{+}),\bP)$ can be represented by $\int_{0}^{t}a^{+}(\Theta^{+}_{s}){\rm d}s$, {either for some strictly positive measurable function $a^{+}$ on $\s$, or such that $a^+(v) =0$ for all $v\in\s$}.

\item[(a3)] $((\xi,\Theta),\bP)$ and $((\xi,\Theta),\hat{\bP})$ are a pair of upwards regular MAPs for which \eqref{condi:weak reversability} is satisfied.

\item[(a4)] $((\xi,\Theta),\hat{\bP})$ satisfies condition \eqref{condi:Harris type}.

\item[(a5)] $\bP_{0,\pi}\left[\sup_{s\in [0,1]}
{|\xi_{s}|}
\right]<+\infty$.

\item[(a6)] The modulator  of the ascending ladder height process $((\xi^{+},\Theta^{+}),\bP)$  is a nonarithmetic aperiodic Harris recurrent process having an invariant distribution $\pi^{+}$ on $\s$ with full support such that
$\bP_{0,\pi^{+}}\left[\xi^{+}_{1}\right]<+\infty$.

\item[(a7)] $\hat{{\rm n}}^{+}_{v}(\zeta=+\infty)>0$  for every $v\in\s$.

\item[(a8)] $\inf_{v\in\s}\left[\ell^{+}(v)+{\rm n}^{+}_{v}\left(1-\mathrm{e}^{-\zeta}\right)\right]>0$ and ${\rm n}^{+}_{v}(\zeta)<+\infty$ for every $v\in\s$.
\end{itemize}

As noted in Section \ref{sec:stationary overshoots}, given conditions (a2) and (a8),
it follows by Corollary \ref{cor:existence of invariant distribution} that
$$\pi^{+}(\cdot)=\frac{1}{\bP_{0,\pi}\left[\bar{L}_{1}\right]}\bP_{0,\pi}\left[\int_{0}^{1}\1_{\{\Theta_{s}\in\cdot\}}{\rm d}\bar{L}_{s}\right]$$
is an invariant distribution for $\Theta^{+}$. Moreover, the Harris recurrence of $\Theta^{+}$ given in (a6) implies that $\pi^{+}$ is the unique invariant distribution for $\Theta^{+}$.

\begin{theorem}\label{thm:main}
    Under assumptions (a1)-(a8), the conclusions (C1)-(C5) in the Introduction are true.
     \end{theorem}

We conclude this section by considering a slight adjustment of the sufficient conditions (a1)-(a8), such that (a5) and (a7) can be replaced by the stronger sufficient conditions (i.e. ones that imply (a5) and (a7)). Our principal aim here is to produce conditions that can be identified in terms of the components of the ascending ladder process of $((\xi, \Theta),{\bP})$ and the ascending ladder process of the dual process $((\xi, \Theta),\hat{\bP})$.
More precisely, we have the following alternative conditions to Theorem \ref{thm:main}.

\begin{theorem}\label{a5a7}
Suppose conditions (a5) and (a7) are replaced by:
\begin{itemize}
\item[(a5)'] The modulator $(\Theta^{+}_{t})_{t\ge 0}$ of the ascending ladder height process $((\xi^{+},\Theta^{+}),\hat{\bP})$,  is an aperiodic Harris recurrent  process having an invariant distribution $\hat{\pi}^{+}$ on $\s$ with full support such that $$\int_{\mathcal{S}} \hat{\pi}^+(\d v) \left[\hat{a}^+(v)  +\hat{\emph{{\rm n}}}^+_{v}(| \epsilon_\zeta | ; \zeta <\infty)\right]<+\infty.$$

\item[(a7)']\quad $\inf_{v\in\mathcal{S}}\hat{\emph{\rm n}}^{+}_{v}(\zeta=+\infty)>0$.

\end{itemize}

Then the conclusion of Theorem \ref{thm:main} is still valid.
\end{theorem}

\begin{remark}
\rm
Before continuing to the proof, let us note that the condition in (a5)' is the natural analogue of (a6). Indeed, note that $\bP_{0,\pi^{+}}\left[\xi^{+}_{1}\right]= \int_{\mathcal{S}} \pi^+(\d v) \left[ a^+(v)  +\emph{{\rm n}}^+_{v}(| \epsilon_\zeta | ; \zeta <\infty)\right]$.
\end{remark}

\noindent\textit{Proof of Theorem \ref{a5a7}}\quad
Condition (a7)' obviously implies (a7).
The proof is based around showing that the new conditions
together with (a1)-(a4) and (a8) imply (a5).
Suppose that $\mathbf{e}_q$ is an independent exponentially distributed random variable with rate $q>0$. On account of the fact that
$t\mapsto\bP_{0,\pi}\left[\sup_{s\in[0,t]}|\xi_{s}|\right]$ is increasing, to show (a5) it suffices to show that
$$\bP_{0,\pi}\left[\sup_{s\in[0,\mathbf{e}_q]}|\xi_{s}|\right] =\int_0^\infty q {\rme}^{-q t}\bP_{0,\pi}\left[\sup_{s\in[0,t]}|\xi_{s}|\right] \d t<\infty.$$
 For the latter, we note that
\begin{align*}
\bP_{0,\pi}\left[\sup_{s\in[0,\mathbf{e}_q]}|\xi_{s}|\right]& \leq \bP_{0,\pi}\left[\bar{\xi}_{\mathbf{e}_q}\right]  -  \bP_{0,\pi}\left[\underline{\xi}_{\mathbf{e}_q}\right].
\end{align*}
Next define
\begin{equation}
\Lambda^+_v(q) :=
\ell^+(v)q +
{\rm n}^+_v(1-{\rm e}^{-q\zeta}),\qquad  q\geq 0.
\label{Lambda}
\end{equation}
Note from Proposition \ref{prop:wiener_hopf} that
\begin{equation}
 \bP_{0,\pi}\left[\bar{\xi}_{\mathbf{e}_q}\right] =  \bP_{0,\pi}\left[\int_{0}^{\infty}\1_{\{\xi^+_s<\infty\}}
     \mathrm{e}^{-q \bar{L}^{-1}_s } \xi^+_s  \Lambda^+_{\Theta^+_s}(q) \d s \right] .
     \label{applyCOM}
     \end{equation}
     Next define the change of measure
     \begin{equation}
\left.\frac{ \d\bP_{0,\theta}^{(q)}  }{\d \bP_{0,\theta}}\right|_{\mathcal{G}_t} =    \mathrm{e}^{-q \bar{L}^{-1}_t + \int_0^t\Lambda^+_{\Theta^+_s}(q)\d s}
     \label{expmg}
     \end{equation}
     for $\theta\in\mathcal{S}$, where $\mathcal{G}_t = \sigma((\bar{L}^{-1}_s, \xi^+_s,\Theta^+_s), s\leq t)$.
     To see why the right-hand side of \eqref{expmg} is a martingale, it suffices to note that $(\bar{L}^{-1}_t ,\Theta^+_t )_{t\ge 0}$ is a MAP and that, for $\theta\in\mathcal{S}$,
      {\[
     \bP_{0,\theta}[  \mathrm{e}^{-q \bar{L}^{-1}_t}| \Theta^+_s: s\leq t] =
     \mathrm{e}^{-\int_0^t\Lambda^+_{\Theta^+_s}(q)\d s}, \qquad t\geq 0,
     \]}
which follows from the the definition \eqref{Lambda} and the fact that the constituent parts of $\Lambda^+_v$, namely $\ell^+(v)$ and ${\rm n}^+_v(1-{\rm e}^{-q\zeta})$ describe the rate at which
{$\bar{L}^{-1}_{s}$}
moves continuously and with jumps given
$\Theta^+_{s}=v$, for $v\in\mathcal{S}$.

     \smallskip

 Using \eqref{expmg}
      in \eqref{applyCOM}, we have
     \[
      \bP_{0,\pi}\left[\bar{\xi}_{\mathbf{e}_q}\right] =  \bP_{0,\pi}^{(q)}\left[\int_{0}^{\infty}
     \mathrm{e}^{-\int_0^s\Lambda^+_{\Theta^+_u}(q)\d u} \xi^+_s  \Lambda_{\Theta^+_s}(q) \d s \right] =   \bP_{0,\pi}^{(q)}\left[\int_{0}^{\infty}
     \mathrm{e}^{-\int_0^s\Lambda^+_{\Theta^+_u}(q) \d u} \d \xi^+_s \right],
     \]
     where the final equality follows by a straightforward integration by parts (recall that the process $\xi^+$ is non-decreasing and therefore has bounded variation paths).
     From (a8),
     we now have that there exists a constant $c>0$ such that for any $q\ge 1$
     \begin{equation}
      \bP_{0,\pi}\left[\bar{\xi}_{\mathbf{e}_q}\right] \leq  \bP_{0,\pi}^{(q)}\left[\int_{0}^{\infty}
     \mathrm{e}^{-c s} \d \xi^+_s \right] =c\bP^{(q)}_{0,\pi} \left[\int_{0}^{\infty}
\mathrm{e}^{-c s} \xi^+_s \d s \right] = c\int_{0}^{\infty}
\mathrm{e}^{-c s} \bP^{(q)}_{0,\pi} \left[\xi^+_s  \right]\d s,
\label{lingrow}
      \end{equation}
      where, again, we have performed an integration by parts.
      Next note that, given $\Theta^+$, the exponent associated to $(\bar{L}^{-1}_t, \xi^+_t)_{t\geq 0}$, is given by
      \[
      \bP_{0,\pi}^{(q)} [{\rm e}^{-\alpha \bar{L}^{-1}_t - \beta\xi^+_t}|\Theta^+] = \exp\left\{- \int_0^t\d s\left[\alpha\ell^+(\Theta^+_s) + \beta a^+(\Theta^+_s) + {\rm n}^+_{\Theta^+_s}((1- {\rm e}^{-\alpha \zeta -\beta \epsilon_\zeta }){\rm e}^{- q \zeta  }; \zeta <\infty)\right]\right\},
      \]
      for $\alpha,\beta,t\geq 0$.
      From this it is easily deduced by differentiation that
      \begin{align*}
       \bP^{(q)}[\xi^+_t|\Theta^+] &= \int_0^t \d s  \left[a^+(\Theta^+_s)  +{\rm n}^+_{\Theta^+_s}(| \epsilon_\zeta | {\rm e}^{- q \zeta  }; \zeta <\infty)\right]\\
       &\leq \int_0^t a^+(\Theta^+_s)  +{\rm n}^+_{\Theta^+_s}(| \epsilon_\zeta | ; \zeta <\infty)\d s \\
       &= \bP[\xi^+_t|\Theta^+] .
       \end{align*}
      Using the ergodic properties of $\Theta^+$ under $\bP$, we can invoke  Theorem 1.1. of \cite{F} and conclude that
       \begin{align*}
       \limsup_{t\to\infty} \frac{1}{t}\bP^{(q)}_{0,\pi}[\xi^+_t] &\leq \lim_{t\to\infty} \frac{1}{t}\bP_{0,\pi}[\xi^+_t] \\
       &= \lim_{t\to\infty}\frac{1}{t}\bP_{0,\pi}\left[\int_0^t a^+(\Theta^+_s)  +{\rm n}^+_{\Theta^+_s}(| \epsilon_\zeta | ; \zeta <\infty)\d s\right]\\
       &
       = \int_{\mathcal{S}} \pi^+(\d v) \left[ a^+(v)  +{\rm n}^+_{v}(| \epsilon_\zeta | ; \zeta <\infty)\right]\\
&       = \bP_{0,\pi^+}[\xi^+_1]
       \end{align*}
Using the above linear growth, it follows from \eqref{lingrow} that $   \bP_{0,\pi}\left[\bar{\xi}_{\mathbf{e}_q}\right]<\infty$.

         \smallskip

Using obvious notation, the analogous object to $\Lambda_v^+(q)$ for the descending ladder height MAP takes the form
\[
\Lambda^-_v(q) =  {\rm n}^-_v(\zeta =+\infty)+\ell^-(v)q + {\rm n}^-_v(1-{\rm e}^{-q\zeta}; \zeta<\infty),\qquad  q\geq 0
\]
(Specifically, we cannot rule out the possibility of killing.)
Let us momentarily assume that the modulator of the descending ladder height process $((\xi^{-},\Theta^{-}),\bP)$  is an aperiodic Harris recurrent  process with  an invariant distribution $\pi^{-}$ on $\s$ with full support such that
$\int_{\mathcal{S}} \pi^-(\d v) \left[ a^-(v)  +{\rm n}^-_{v}(| \epsilon_\zeta | ; \zeta <\infty)\right]<+\infty$ and $\inf_{v\in\mathcal{S}}{\rm n}^{-}_{v}(\zeta=+\infty)>0$.
Following the above computations, albeit using the last lower bound to
justify the lower bounding constant $c$ in \eqref{lingrow},
we can show that $\bP_{0,\pi}\left[\underline{\xi}_{\mathbf{e}_q}\right] <\infty$.

 \smallskip

 To complete the proof, we need to show that the assumptions in the last paragraph match those in the statement of the theorem by verifying that  $\bP_{0,\pi^-}\left[\xi^-_1\right]= \hat{\bP}_{0,\hat{\pi}^+}\left[\xi^+_1\right]$.
 Thanks to
 the weak reversal relation
 between  $\bP$ and $\tilde\bP$ (see the discussion below Lemma \ref{lem:walsh_killing_time}), we have that  $\bP_{0,\pi^-}[\xi^-_1] = \tilde\bP_{0,\tilde\pi^-}[\xi^-_1]$, where $\tilde{\pi}^-$ plays the role of $\pi^-$ but for $((\xi,\Theta),\widetilde{\bP})$.  The relation between $\tilde\bP$ and $\hat\bP$ then implies that $\tilde\pi^- = \hat\pi^+$ and $ \tilde\bP_{0,\tilde\pi^-}[\xi^-_1] = \hat\bP_{0,\pi^+}[\xi^+_1] $ as required.\qed
\smallskip

The remainder of the paper is devoted to the proof of Theorem \ref{thm:main}; as such,  we always assume conditions
(a1)-(a8)
hold unless otherwise stated. Before moving to the proof of Theorem \ref{thm:main}, we first engage in a little discussion concerning its applicability.

\section{Applicability of results}\label{sec:application}
There are two immediate cases of interest: The case of a Brownian motion in a cone and the case of a stable processes in a cone.  The general philosophy of the proofs of Theorems \ref{thm:main} and \ref{a5a7} is  to use a judiciously chosen
harmonic function to construct a process whose dual can serve as the desired ssMp  entering at the origin. In the two examples below, we verify the general criteria of the aforesaid theorem(s). The reader will note that, in both cases, it is first necessary to construct the notion of the process conditioned to remain in a cone (appealing to an appropriate change of measure, which ultimately comes from  a harmonic function constructed on the Martin boundary of the base process killed on exiting a cone) and transfer that notion into the general setting of required criteria on MAPs. One may argue that it may prove to be easier to construct the candidate process entering at the origin in a direct way rather than via Theorem  \ref{thm:main} or \ref{a5a7}. Indeed this was the approach in \cite{KRS}. It is also worthy of noting that the majority of the conditions in Theorem  \ref{thm:main}  (and hence Theorem  \ref{a5a7}) boils down to controlling the stability of the ssMp in order to obtain the Skorokhod limit in (C4), $\lim_{x\to0}\mathbb{P}_x$, rather than the existence of the limit $\mathbb{P}_0$.

In the remainder of this section we will also discuss further open problems that could in principle be analysed appealing to the fluctuation theory of the Lamperti--Kiu decomposition. Moreover, we also discuss   the reason why general fluctuation theory of MAPs is deserving of further investigation given what has been laid out in this article.

\subsection{Brownian motion in a cone}\label{bcone}
We are interested in cones of the form
\begin{equation}\label{def:cone}
\Gamma = \{x=r\theta :r> 0,\theta\in \Omega\},
\end{equation}
where $\Omega$ is a non-empty open subset on $\mathbb{S}^{d-1}$. We assume further that there is a complete set of eigenfunctions
$\{m_j:\mathbb{S}^{d-1}\to\mathbb{R},\ j\ge 1\}$, which are orthonormal with respect to the surface measure on $\Omega$ with eigenvalues $0<\lambda_1<\lambda_2\leq \lambda_3\leq \ldots$, which satisfy
\begin{eqnarray*}
\Delta^{d-1}m_j=-\lambda_jm_j &on& \Omega\\
m_j=0 &on& \partial\Omega,
\qquad  j\geq 1
\end{eqnarray*}
 such that every boundary point of $\Omega$ is regular for the above Dirichlet problem,
where $\Delta^{d-1}$ denotes the Laplace-Beltrami operator on $\mathbb{S}^{d-1}$.

Suppose {$(B,\P)$ is a $d$-dimensional $(d\ge 2)$ standard Brownian motion and }that $\Gamma$ is a regular cone in $\mathbb{R}^d$. Let 
$\tau^\Gamma := \inf\{t>0: B_t \not \in \Gamma\}$.
From Theorem B of \cite{garbit2014} (see also \cite{banuelos1997}), it is known that there exists a constant $\kappa>0$ such that
\begin{equation*}
\mathbb{P}_x(\tau^{\Gamma}>t)=\kappa \|x\|^{p}m_1(\arg(x)) t^{-p/2}(1+o(1)),\;\;\;t\rightarrow \infty,
\end{equation*}
where $p=\sqrt{\lambda_1+
(d/2-1)^2}-(d/2-1)$
.
It thus relatively straightforward to show that,
for $A\in \F_t=\sigma (B_u:u\leq t)$ and $x\in \Gamma$,
\begin{equation}\mathbb{P}_x^{\Gamma}(A):=        \lim_{s\rightarrow\infty}\mathbb{P}_x(A \,|\, \tau^{\Gamma}>t+s)
\label{conecondition}
\end{equation}
defines a family of conservative probabilities on the space of continuous paths such that
\begin{equation}\label{eq-changemeasure}
\left.\frac{d\mathbb{P}_x^{\Gamma}}{d\mathbb{P}_x}\right\vert_{\F_t}:=\frac{M(B_t)}{M(x)}\1_{\{t<\tau^{\Gamma}\}},\;\; t\geq 0,
\end{equation}
where \[M(x) = \|x\|^{p}m_1(\arg(x))\] is a harmonic function in the cone. In particular, the right hand side of \eqref{eq-changemeasure} is a martingale. Furthermore, if $\mathbb{P}^{\Gamma}=\{\mathbb{P}_x^{\Gamma},x\in \Gamma\}$, then the process $(B,\mathbb{P}^{\Gamma})$ is a ssMp.
 \smallskip

 Brownian motion conditioned to stay in a cone has previously been considered in \cite{garbit2009, garbit2014, banuelos1997, Shimura}. Only in \cite{garbit2009} was the notion of a Brownian motion conditioned to stay in the cone issued from the apex considered. In that case, the authors built the law of the Brownian motion conditioned to remain in the cone and survive for at least one unit of time from a point $x$ away from the apex, and then showed the weak limit on path space as $x$ converges to the apex of the cone. The authors described their construction as analogous to the construction of the Brownian meander in the upper half-line (i.e., a Brownian meander in the cone).  Independently  to the aforementioned work, Theorem \ref{thm:main} provides an easy route to the construction of Brownian motion conditioned to stay in a cone, issued from the apex, as the weak limit on path space of the conditioned process issued from any other point in the interior of the cone.

To understand how to see $(B,\mathbb{P}^{\Gamma})$ as a ssMp, we first consider standard Brownian motion as such. Its Lamperti--Kiu decomposition has  MAP $(\xi, \Theta)$, with probabilities $\mathbf{P} =\{\mathbf{P}_{y, \theta}, y\in\mathbb{R},\theta\in \mathbb{S}^{d-1}\}$, has the property that $\xi$ is independent of $\Theta$. Moreover,  $\xi$ is a Brownian motion  with  drift on $\mathbb{R}$, whose Laplace exponent is given by   $\psi(\theta) = \theta^2+(d-2)\theta$.
{Since $d\ge 2$, the drift is strictly positive.}
Note that
\begin{equation}
\frac{M(B_t) }{M(x)}\mathbf{1}_{\{t<\tau^\Gamma\}}= {\rm e}^{p(\xi_{\varphi(t)} -\log\|x\|)} \frac{m_1(\Theta_{\varphi(t)})}{m_1(\arg(x))}\mathbf{1}_{\{t<\tau^\Gamma\}}, \qquad t\geq 0.
\label{mapMG}
\end{equation}
Hence, recalling that $\varphi(t)$ is a stopping time, we can think of the change of measure \eqref{eq-changemeasure} as the product of an Esscher transform on $\xi$ and the natural change of measure on $\Theta$ corresponding to a Doob $h$-transform with $h = m_1$. Indeed, treating $s = \varphi(t)$ as the natural timescale of $(\xi,\Theta)$, we have the martingale on the right-hand side of \eqref{mapMG} equal to
\begin{align}
{\rm e}^{p(\xi_s -\log\|x\|)} \frac{m_1(\Theta_s)}{m_1(\arg(x))}\mathbf{1}_{\{s<\texttt{k}^\Omega\}} &= {\rm e}^{p(\xi_s-\log\|x\|) -\psi(p)s}  \times {\rm e}^{\psi(p)s} \frac{m_1(\Theta_s)}{m_1(\arg(x))}\mathbf{1}_{\{s<\texttt{k}^\Omega\}} \notag\\
& = {\rm e}^{p(\xi_s-\log\|x\|) -\psi(p)s} \times {\rm e}^{\lambda_1 s}\frac{m_1(\Theta_s)}{m_1(\arg(x))}\mathbf{1}_{\{s<\texttt{k}^\Omega\}} ,\qquad s\ge 0,
\label{xiThCOM}
\end{align}
where  $\psi(p) = p^2+(d-2)p =\lambda_1$ and $\texttt{k}^\Omega = \inf\{t>0: \Theta_t \not\in \Omega\}$.
\smallskip

As such, we see that the Brownian motion conditioned to stay in a cone has underlying MAP $(\xi, \Theta)$ with probabilities $\mathbf{P}^\Gamma = \{\mathbf{P}^\Gamma_{y,\theta}, y\in \mathbb{R}, \theta\in\Omega\}$, such that $\mathbf{P}^\Gamma$ is absolutely continuous with respect to $\mathbf{P}$, with the change of measure on $\sigma((\xi_u,\Theta_u), u<s )$ given by
the right-hand side of \eqref{xiThCOM}. Because of the factorisation of \eqref{xiThCOM} into a martingale acting on the law of $\xi$ and another acting on the law of $\Theta$, it is clear that the components of the pair $(\xi,\Theta)$ under $\mathbf{P}^\Gamma$ are still independent, moreover, $\xi$ still belongs to the family of Brownian motion with strictly positive drift; the latter being characterised by the Laplace exponent
\[
\psi^\Gamma(\theta) = \psi(\theta + p) -\psi(p)
=\theta^2 +   (d + 2p -2)\theta
\]
Now referring back to Theorem \ref{thm:main}, we may consider the question of whether the process conditioned to remain in the cone may be issued from its apex. To the best of our knowledge this would offer a new result.

\smallskip

To this end, let us consider each of the assumptions (a1)-(a8). 
\begin{itemize}[(a4)-(a5):]
\item[(a1):] In $d$-dimensions, as alluded to previously, $\xi$ and $\Theta$ are independent. Whilst we have identified above $\xi$ as a Brownian motion with drift, we should mention that $\Theta$ is known from Section 7.15 of \cite{IM}. The Feller property is an easy consequence of known results given there.

\item[(a2):] As $\xi$ is a Brownian motion with drift, independent of $\Theta$, the continuous part of $\xi^+$ is nothing more than a pure drift. Hence $a^+(\theta) = a^+>0$.
 The existence of a discrete spectrum of the Laplace--Beltrami operator on $\Omega$, the independence of $\Theta$ from $\xi$ and thus developing the semigroup of $\Theta$ as a spectral expansion allows us to easily deduce that the latter converges weakly to a stationary distribution, $\pi$. In addition, for bounded measurable function $g:\Omega\to\mathbb{R}$, we can identify $\pi$ via
\[
\int_\Omega g(\theta) \pi(\d \theta)=\int_\Omega g(\theta) \tilde{m}_1(\theta)m_1(\theta)\Sigma(\d\theta),
\]
{where $\Sigma(\d \theta)$ is the uniform measure on the sphere $\mathbb{S}^{d-1}$ and $\tilde{m}_{1}$ is the left eigenfunction associated to $\lambda_{1}$, which is normalized to satisfy that $\int_{\Omega}\tilde{m}_{1}(\theta)m_{1}(\theta)\Sigma(\d \theta)=1$.}
\smallskip

\item[(a3):] This condition is easily satisfied thanks to the duality properties of Brownian with drift and the reversibility of $\Theta$ with respect to $\pi$.
\item[(a4)-(a5):] Note that the dual process $((\xi,\Theta),\hat{\bP}^\Gamma)$ to $((\xi,\Theta), \mathbf{P}^\Gamma)$ is equal in law to $((-\xi,\Theta), \mathbf{P}^\Gamma)$ (thanks to the independence of $\xi$ and $\Theta$ and the fact that $\xi$ is a Brownian motion with drift).  Rather than verifying the criteria (a4) and (a5) directly, we can refer back to their use in the  proof of Theorem  \ref{thm:slln} and Proposition \ref{prop:tail}. In the spirit of
Remark \ref{HTremark}, we note that it suffices to show that $\lim_{t\to\infty}\xi_t/t$ exists $\mathbf{P}^\Gamma$-almost surely.
 This is,
 of course, a trivial consequence of the independence of $\Theta$ and $\xi$ as well as the fact that $\xi$ under
 $\bP^{\Gamma}$
 is a Brownian motion with drift.
%

\item[(a6):] This requirement is fundamentally needed for Section \ref{sec:stationary overshoots}, in order to analyse overshoot distributions of the MAP using
    the Markov renewal theory.
    In the current setting, on account of the fact that there are no overshoots, only a creeping term, and that $\xi$ and $\Theta$ are independent, the only requirement needed is that a
    stationary distribution
    for $\Theta^+$ exists. However, this was dealt with in (a2).


\item[(a7):] The excursion measure $\hat{{\rm n}}^{+}_{v}$ of $((\xi,\Theta),\hat{\bP}^{\Gamma})$ from its running maximum does not depend on $v\in\Omega$, and it agrees with the excursion measure of $((\xi,\Theta),\bP^{\Gamma})$ from its running minimum. As $\xi_t\to\infty$, $\mathbf{P}^\Gamma$-a.s., we find easily that $\hat{{\rm n}}^{+}(\zeta=+\infty)>0$.

\item[(a8):] Once again, independence of $\xi$ and $\Theta$ under $\mathbf{P}^\Gamma$ ensures that none of the items in this assumption depend on  $v\in\Omega$. As such, we note that  ${\rm n}^{+}$ is played by the role of the excursion measure of $((\xi,\Theta),\bP^{\Gamma})$ from its maximum. Moreover, $\ell^+$  is the drift coefficient of the inverse local time of $\xi$ at its maximum, which is zero. The first part of the assumption follows immediately as ${\rm n}^{+}\left(1-\mathrm{e}^{-\zeta}\right)>0$. In order to verify the second part of the assumption, i.e. that   ${\rm n}^{+}(\zeta)<+\infty$, it suffices to recall the Wiener-Hopf factorisation for a Brownian motion with strictly positive drift. It is classically known (see, e.g., Section 6.5 of \cite{Kbook})
    that the inverse local time at the maximum is a tempered stable process, and thus has finite mean, which necessarily implies ${\rm n}^{+}(\zeta)<+\infty$.
\end{itemize}

 \subsection{Stable process in a cone}\label{scone}

Recently \cite{KRS} resolved the matter of conditioning a $d$-dimensional ($d\ge 2$) isotropic $\alpha$-stable ($\alpha\in (0,2)$) process to remain in a cone $\Gamma$ of the form defined in \eqref{def:cone} where $\Omega$ is an open set on $\mathbb{S}^{d-1}$. As here, their approach relied on the Markov renewal theory, albeit the application was not undertaken explicitly in the context of excursion theory.

As in the Brownian setting,
the conditioning \eqref{conecondition} can be made sense of, resulting in a change of measure as in \eqref{eq-changemeasure}. Also similarly to the Brownian case, the harmonic function $M$ {is a locally bounded function on $\R^{d}$ which vanishes on $\R^{d}\setminus\Gamma$ and} has the property that it can be written in the form
\begin{equation}
M(x) = \|x\|^p M(\arg(x))
\label{M7}
\end{equation}
for some
$p\in(0,\alpha)$.
Unlike the Brownian case, the underlying MAP $((\xi,\Theta),\bP)$ of the isotropic stable process does not have the property that the ordinate is independent of the modulator. The coupling between the two is complicated to describe (cf. \cite{KALEA}), moreover, it forces the process $((\xi, \Theta), \mathbf{P}^\Gamma)$ to similarly display coupling.

Nonetheless, what is similar to the Brownian setting is that, thanks to isotropy itself, the  ordinate process of the  MAP $((\xi, \Theta), \mathbf{P})$
is a L\'evy process. This  is known in explicit detail via its characteristic exponent
\begin{equation}
\Psi_\xi(\theta) = \frac{\Gamma(\frac{1}{2}(-{\rm i}\theta +\alpha ))}{\Gamma(-\frac{1}{2}{\rm i}\theta)}
\times
\frac{\Gamma(\frac{1}{2}({\rm i}\theta +d))}{\Gamma(\frac{1}{2}({\rm i}\theta +d-\alpha))}, \qquad \theta\in\mathbb{R}.
\label{a}
\end{equation}
We can also identify the associated Laplace exponent of $(\xi, \mathbf{P})$ via the relation $\psi(\lambda) = -\Psi_{\xi}(-{\rm i}\lambda)$, providing ${\rm Re}(\lambda)\in (-d, \alpha)$.
As the ordinate $\xi$ alone is a L\'evy process, it follows that, when seen as a change of measure on $(\xi, \Theta)$,  \eqref{eq-changemeasure} can be better written as
{\begin{equation}\label{eq-changemeasure2}
\left.\frac{{\rm d}\mathbf{P}_{0,\theta}^{\Gamma}}{{\rm d}\mathbf{P}_{0,\theta}}\right\vert_{\sigma((\xi_{s},\Theta_{s}):s\le t)}:={\rm e}^{p\xi_t-\psi(p)t}\frac{M(\Theta_{t})}{M(\theta)}{\rm e}^{\psi(p)t}\1_{\{t<\texttt{k}^{\Gamma}\}},\;\; t\geq 0,
\end{equation}}

From \eqref{a} its Wiener-Hopf factorisation (indicated by its multiplicative sign) is also explicit, see \cite{KALEA, Deep3} for more details. Moreover, the ordinate process $\xi$ under $\mathbf{P}^\Gamma$ can also be seen as the result of a generalised Esscher transform. However, the loss of isotropy in the cone means that the ordinate is no longer a L\'evy process. Nonetheless, the excursion theory of $((\xi, \Theta), \mathbf{P}^\Gamma)$ can be related back to that of
$((\xi,\Theta),\bP)$.

In terms of verifying assumptions (a1)-(a8), it is unsurprising that some of them formed part of the proof in \cite{KRS}, or follow as easy corollaries of those proved there. As such Theorem \ref{thm:main} offers an alternative way of assembling many of the smaller results in \cite{KRS}.
As in Section \ref{bcone}, we run through (a1)-(a8) below. We want to stress that all but (a2) hold for in any dimension $d\geq 2$, however, for technical reasons we are only able to deal with $d= 2$ (see Remark \ref{d>2} below).

\begin{itemize}[(a4)-(a5):]
\item[(a1):] The MAP underlying the isotropic stable process is Feller. Being a generalised Esscher transform (i.e., a Doob $h$-transform) of the former, it is straightforward to verify that $((\xi,\Theta), \mathbf{P}^\Gamma)$ is a Feller process.

\item[(a2):]
{From Example \ref{eg3} we know that $((\xi,\Theta),\bP)$ is dual to $((\xi,\Theta),\hat{\bP})$ with respect to the measure $\mathrm{Leb}\otimes \pi$ where $(\xi,\Theta)$ under $\hat{\bP}_{x,\theta}$ is equal in law to $(-\xi,\Theta)$ under $\bP_{-x,\theta}$ and $\pi$ is the uniform measure on $\mathbb{S}^{d-1}$.  We claim in the following that condition \eqref{condi:weak reversability} is satisfied by $((\xi,\Theta),\bP^{\Gamma})$ and $((\xi,\Theta),\hat{\bP}^{\Gamma})$ with respect to the measure $\pi^{\Gamma}({\rm d}\theta):=\1_{\{\theta\in\Omega\}}M(\theta)^{2}\pi(d\theta)$, where $\hat{\bP}^{\Gamma}=\{\hat{\bP}^{\Gamma}_{x,\theta},\ x\in\R,\ \theta\in\Omega\}$ and $\hat{\bP}^{\Gamma}_{x,\theta}$ denotes the law of $(-\xi,\Theta)$ under $\bP^{\Gamma}_{-x,\theta}$. In fact, for any $t>0$ and any nonnegative measurable functions $f,h:\Omega\to \R^{+}$ and $g:\R\to\R^{+}$ we have
\begin{align*}
\int_{\Omega}f(\theta)\bP^{\Gamma}_{0,\theta}\left[g(\xi_{t})h(\Theta_{t})\right]\pi^{\Gamma}({\rm d}\theta)
&=\int_{\Omega}M(\theta)f(\theta)\bP_{0,\theta}\left[\mathrm{e}^{p\xi_{t}}M(\Theta_{t})g(\xi_{t})h(\Theta_{t})\1_{\{t<\texttt{k}^{\Omega}\}}\right]\pi({\rm d}\theta)\\
&=\bP_{0,\pi}\left[M(\Theta_{0})f(\Theta_{0})\mathrm{e}^{p\xi_{t}}M(\Theta_{t})g(\xi_{t})h(\Theta_{t})\1_{\{t<\texttt{k}^{\Omega}\}}\right]\\
&=\hat{\bP}_{0,\pi}\left[M(\Theta_{t})f(\Theta_{t})\mathrm{e}^{p(\xi_{0}-\xi_{t})}M(\Theta_{0})g(\xi_{0}-\xi_{t})h(\Theta_{0})\1_{\{t<\texttt{k}^{\Omega}\}}\right]\\
&=\bP_{0,\pi}\left[M(\Theta_{0})h(\Theta_{0})\mathrm{e}^{p\xi_{t}}M(\Theta_{t})f(\Theta_{t})g(\xi_{t})\1_{\{t<\texttt{k}^{\Omega}\}}\right]\\
&=\int_{\Omega}M(\upsilon)h(\upsilon)\bP_{0,\upsilon}\left[\mathrm{e}^{p\xi_{t}}M(\Theta_{t})f(\Theta_{t})g(\xi_{t})\1_{\{t<\texttt{k}^{\Omega}\}}\right]\pi({\rm d}\upsilon)\\
&=\int_{\Omega}h(\upsilon)\bP^{\Gamma}_{0,\upsilon}\left[g(\xi_{t})f(\Theta_{t})\right]\pi^{\Gamma}({\rm d}\upsilon).
\end{align*}
The third equality follows from Lemma \ref{lem:time reversal}. In conclusion we've proved that
\begin{equation}\nonumber
\bP^{\Gamma}_{0,\theta}\left(\xi_{t}\in {\rm d} z,\Theta_{t}\in {\rm d}\upsilon\right)\pi^{\Gamma}({\rm d}\theta)=\bP^{\Gamma}_{0,\upsilon}\left(\xi_{t}\in {\rm d}z,\Theta_{t}\in {\rm d}\theta\right)\pi^{\Gamma}({\rm d}\upsilon).
\end{equation}
By integrating the above equation over $z$, it follows that $(\Theta,\bP^{\Gamma})$ is dual to itself with respect to $\pi^{\Gamma}$. Without loss of generality we assume $\pi^{\Gamma}(\Omega)=1$. Thus $\pi^{\Gamma}$ is an invariant distribution of $(\Theta,\bP^{\Gamma})$.
}
\smallskip

{Next, note from Proposition 20.17 in \cite{Kallenberg}, as soon as we know that $(\Theta,\mathbf{P}^\Gamma)$ is a {\it regular Feller} process, on account of the fact that \eqref{eq-changemeasure} describes a conservative process, it follows that  $\Theta$ is Harris recurrent. Here, by `regular Feller', the process has the Feller property as well as its transition semigroup being absolutely continuous with respect to some
{locally finite}
measure, such that the transition density is jointly continuous in time and its spatial variables.
\smallskip

 {We consider $d=2$ only.} We know from \cite{BW} that, under $\mathbf{P}$,  $ \Theta$ can be written in the form of $\exp({\rm i} \vartheta_t)$, where $\vartheta: =\{\vartheta_t, t\geq 0\}$ is a pure jump L\'evy process (the winding number) whose L\'evy measure is known. From
 {\cite[Section 5]{Kal}}
 (see also \cite[Theorem 1]{KS} for more recent results in this direction) and the form of the L\'evy measure given in \cite{BW},  $m({\rm d}u)= c\int_{\mathbb{R}^2} \mathbf{1}_{\{\arg(1+x)\in \d u\}}\|x\|^{-(2+\alpha)}\d x$ (where $c$ is an unimportant constant), it is straightforward to verify that the sufficient condition (up to an unimportant multiplicative constant)
 \[
 \lim_{\varepsilon\to0}\frac{\int_{(-\varepsilon, \varepsilon)}u^2 m(\d u)}{\varepsilon^2|\log \varepsilon|} =
  \lim_{\varepsilon\to0}\frac{c\int_{(-\varepsilon, \varepsilon)}u^2  \int_{0}^\infty  |r{\rm e}^{{\rm i}u}-1|^{-(2+\alpha)}   r\d r \d u }{\varepsilon^2|\log \varepsilon|}
 =+\infty
 \]
holds, in which case $\vartheta$ is a L\'evy process
{having transition density function with respect to the Lebesgue measure, which is bounded continuous and vanishes at infinity in its spatial variable.}
Note that the Feller property also ensures that it is continuous in its temporal variable.
Suppose we write the latter as $\{\texttt{p}_t(x), t\geq 0, x\in\mathbb{R}\}$.  This automatically provides us with a density for the process $\vartheta$ killed on exiting $\Omega$. (Note, we are slightly abusing notation here and we are now interpreting $\Omega$ as a subset of $(-\pi,\pi]$.) Indeed, we can write the latter as
 \[
 \texttt{p}^\dagger_t(x,y) =  \texttt{p}_t(y-x) - \mathbf{P}_{0,x}[\mathbf{1}_{\{\texttt{k}^\Omega<t\}} \texttt{p}_{t-\texttt{k}^\Omega} (y-\vartheta_{\texttt{k}^\Omega})], \qquad x,y\in\Omega.
 \]
 Note that the Feller property again gives us continuity in $t$.
Moreover, from e.g. the discussion between (2.18) and (2.19) in \cite{BPW}, it
is straightforward to show that $M$ is continuous. This tells us that $(\vartheta, \mathbf{P}^\Gamma)$ also has
{a transition density function which is continuous in both of its spatial variables,}
i.e., in light of \eqref{eq-changemeasure2}, $\texttt{p}_t^\Gamma(x, y) =\exp(\psi(p)t)M(y)\texttt{p}_t^\dagger(x,y)/M(x)$, $x,y\in\Omega$; and similarly continuity in $t$ follows from the Feller property. The joint continuity of $(t,x,y)\mapsto \texttt{p}_t^\Gamma(x, y)$ now follows by a classical  epsilon-delta chasing argument with multiple use of the triangle inequaltiy.
Returning to the previous paragraph, we thus conclude that $(\Theta, \mathbf{P}^\Gamma)$ is Harris recurrent. This is slightly weaker than the required positive recurrence, however Theorem 20.20 in \cite{Kal} ensures that we must be in the positive recurrent setting.}
\smallskip

%

Finally for the third requirement,  the ordinate of the ascending ladder process $((\xi^{+},\Theta^{+}),\bP)$ underlying the  isotropic stable process has no continuous part (as $\xi$ alone is a L\'evy process with no drift to its ascending ladder process, cf. \eqref{a}), hence after the change of measure corresponding to conditioning the stable process to remain in a cone, the same is true for $((\xi^{+},\Theta^{+}),\bP^{\Gamma})$. That is to say, for  $((\xi^+,\Theta^+), \mathbf{P}^\Gamma)$, $a^{+}(v)=0$ for all $v\in\s$.

\item[(a3):] {
It was verified in (a2) that condition \eqref{condi:weak reversability} holds for $((\xi,\Theta),\bP^{\Gamma})$ and $((\xi,\Theta),\hat{\bP}^{\Gamma})$.} Moreover, since the former can be described as a Doob $h$-transform with respect to
    $((\xi,\Theta),\bP)$, and since this process is both upwards and downwards regular, then the same is true of $((\xi,\Theta),\bP^\Gamma)$ and its dual.

\item[(a4):] 
{In the spirit of Remark \ref{HTremark} it suffices to show that $(\Theta,\bP^{\Gamma})$ has a skeleton process of the form $\{\Theta_{n\delta}:n\ge 1\}$ for some $\delta>0$ that is Harris recurrent. {This follows from our computations in (a2)
{and Example 3.1 in \cite[Chapter VII, section 3]{AsmussenQueue}}.} 
}

\item[(a5):] Define $\xi^*_t = \sup_{s\leq t}|\xi_s|$. Appealing to the change of measure between $((\xi,\Theta), \mathbf{P})$ and $((\xi,\Theta), \mathbf{P}^\Gamma)$, we have, for some constant $C>0$ (whose value may change from line to line) and {{$0<\varepsilon\ll 1$}},
\begin{align*}
\bP^\Gamma_{0,\pi^\Gamma}\left[\xi^*_1
\right]
&=
\int_{\Omega}{M(\theta)^{2}}\pi(\d \theta) \bP_{0,\theta}\left[\xi^*_1{\rm e}^{p\xi_1}\frac{M(\Theta_1)}{M(\theta)}\mathbf{1}_{\{t<\texttt{k}^\Omega\}}
\right]\\
&\leq C \int_{\Omega}\pi(\d \theta) \bP_{0,\theta}\left[\xi^*_1{\rm e}^{p\xi^*_1}
\right]\\
&\leq C \bP_{0}\left[{\rm e}^{(p+\varepsilon)\xi^*_1}
\right],
\end{align*}
where $\varepsilon$ is chosen so that $p+\varepsilon<\alpha$ and we have abused notation and written
{$(\xi,\bP_{0})$ to denote the L\'{e}vy process with characteristic exponent \eqref{a} issued from the origin.}
Appealing to Theorem of 24.18 of \cite{Sato}, it follows that
$
\bP_{0}\left[{\rm e}^{(p+\varepsilon)\xi^*_1}\right]<+\infty
$
if and only if
$
\bP_{0}\left[{\rm e}^{(p+\varepsilon)|\xi_1|}\right]<+\infty.
$
The latter condition occurs if and only if
$\bP_{0}\left[{\rm e}^{(p+\varepsilon){\xi}_1}\vee{\rm e}^{-(p+\varepsilon){\xi}_1}\right]<+\infty.
$ However,
\[
\bP_{0}\left[{\rm e}^{(p+\varepsilon){\xi}_1}\vee{\rm e}^{-(p+\varepsilon){\xi}_1}\right]\leq
\bP_{0}\left[{\rm e}^{(p+\varepsilon){\xi}_1}+{\rm e}^{-(p+\varepsilon){\xi}_1}\right] =
{\mathrm{e}^{\psi(p+\varepsilon)}+\mathrm{e}^{\psi(-(p+\varepsilon))},}
\]
and the latter is finite on account of the fact $\psi(\lambda)$ analytically extends to the interval $(-d, \alpha)$ and $p+\varepsilon<\alpha<d$.

\item[(a6):] 
The modulator  of the ascending ladder height process $((\xi^{+},\Theta^{+}),\bP^{\Gamma})$ was shown to be an nonarithmetic aperiodic Harris recurrent process having an invariant distribution $\pi^{\Gamma,+}$ on $\Omega$ with full support in
{Theorem 4.2 and Theorem 4.3 of \cite{KRS}.}
 Moreover, from Remark 4.1
 of \cite{KRS}, we recall that $M(v)^{-1}\pi^{\Gamma, +}(\d v)$ is an invariant distribution for
 {the process $((\Theta^{+}_{t}\1_{\{t<\texttt{k}^{+,\Omega}\}})_{t\ge 0},\bP)$, where $\texttt{k}^{+,\Omega} = \inf\{t>0: \Theta^+_t\not\in \Omega\}$.}
 We can thus write
 \begin{align}
\bP^\Gamma_{0,\pi^{\Gamma, +}}\left[\xi^{+}_{1}\right]
&= \int_\Omega \frac{\pi^{\Gamma, +}(\d \theta)}{M(\theta)}\bP_{0,\theta}
\left[\xi^{+}_{1}{\rm e}^{p\xi^+_1} M(\Theta^+_1)\1_{\{
1<\texttt{k}^{+,\Omega}\}}\right]
\leq \int_\Omega \frac{\pi^{\Gamma, +}(\d \theta)}{M(\theta)}
\bP_{0,\theta}
\left[{\rm e}^{(p+\varepsilon)\xi^+_1} \right],
\label{rhsfinite1}
\end{align}
where
{$0<\varepsilon\ll1$.}
The right-hand side of \eqref{rhsfinite1} is finite on account of the fact that $(\xi^+,\mathbf{P})$ is a subordinator, whose Laplace exponent is given by
\begin{equation}
\kappa(\lambda): = \frac{\Gamma((\lambda+\alpha )/2)}{\Gamma(\lambda/2)}, \qquad \lambda>-\alpha.
\label{asclap}
\end{equation}
It is clear from the Laplace exponent above that $\xi^+$ has finite exponential moments, in particular, the right hand side of \eqref{rhsfinite1} is finite, providing $p+\varepsilon<\alpha$
.

\item[(a7):] As mentioned above, the dual of $((\xi,\Theta),\bP^\Gamma)$ is equal in law to $((-\xi,\Theta),{\bP}^\Gamma)$.
It follows that the ascending ladder process of the dual is equal in law to the descending ladder process of $((\xi,\Theta),\bP^\Gamma)$.
It follows that the associated excursion measure at the maximum of the dual MAP, $\hat{\rm n}^{\Gamma, +}_v$, agrees with the excursion measure at the minimum, ${\rm n}^{\Gamma, -}_v$. Hence our objective is to understand ${\rm n}^{\Gamma,-}_v(\zeta = +\infty)$.
Suppose now that $((\xi^-,\Theta^-),\bP)$ is the descending laddder height process of the isotropic stable process. From \eqref{a}, we know that $\xi^-$ alone is a killed pure jump subordinator with Laplace exponent
\[
\hat\kappa(\lambda):  =   \frac{\Gamma(\frac{1}{2}(\lambda +d))}{\Gamma(\frac{1}{2}(\lambda +d-\alpha))}, \qquad \lambda \geq 0,
\]
and hence has killing rate $q^- = {\Gamma(d/2)}/{\Gamma((d-\alpha)/2)}$. When taking account of the killing rate for the coupled system $((\xi^-,\Theta^-), \bP)$, isotropy ensures that the killing rate of the pair is also equal to $q^-$. In other words, the lifetime of the pair $((\xi^-,\Theta^-), \bP)$, written as $\texttt{k}^-$,
is an independent and exponentially distributed random variable with parameter $q^-$. 
{From \cite[Remark 4.1]{KRS}, we know that $((\xi^{-},\Theta^{-}),\bP^{\Gamma})$ can be described as a Doob $h$-transform with respect to $((\xi^{-},\Theta^{-}),\bP)$. Hence the lifetime of $((\xi^{-},\Theta^{-},\bP^{\Gamma}))$ is an independent and exponentially distributed random variable with parameter $q^{-}$.} In conclusion we have that  $\hat{{\rm n}}^{\Gamma,+}_{v}(\zeta=+\infty) = {\rm n}_{v}^{\Gamma,-}(\zeta = +\infty) = q^->0$  for every $v\in\Omega$.

\item[(a8):] For the MAP $((\xi,\Theta), \bP)$ underlying the isotropic stable process, the local time of the ordinate at its maximum is simply that of the L\'evy process at its maximum. As such the drift component of the inverse local time, say $\ell^+$, does not depend on $v\in\Omega$. In fact, as a consequence of the fact that $(\xi,\bP)$ is a L\'evy process with unbounded variation, we can show that $\ell^+=0$. However, the excursion measure ${\rm n}^+_\upsilon$ does depend on $\upsilon$. 
\smallskip

Denote by $\ell^{\Gamma,+}$ the drift component of the additive functional describing local time of $((\xi,\Theta),\bP^\Gamma)$ at the running maximum in the spirit of \eqref{M1}.  Applying a similar analysis to the verification of (a7),
we can deduce in a straightforward manner that
$\ell^{\Gamma, +}  = \ell^+$ and
\begin{align*}
{\rm n}_v^{\Gamma, +}(1-{\rm e}^{-\zeta})&= \frac{1}{M(v)}{\rm n}_v^{ +}\left({\rm e}^{p|\epsilon_\zeta|}M(\nu_\zeta)(1-{\rm e}^{-\zeta}); \zeta< \texttt{k}^\Omega\right)\\
&\geq \frac{1}{M(v)}{\rm n}_v^{ +}\left({\rm e}^{p|\epsilon_\zeta|}M(\nu_\zeta); \zeta< \texttt{k}^\Omega\right),
\end{align*}
where we have abused slightly our notation for $\texttt{k}^\Omega$ in the obvious way.  A consequence of the fact that ${\rm e}^{p\xi_t}M(\Theta_t)\1_{(t<\texttt{k}^\Omega)}$
is a martingale is that
\[
   \ell^{+}M(v) + {\rm n}_v^{ +}\left({\rm e}^{p|\epsilon_\zeta|}M(\nu_\zeta); \zeta< \texttt{k}^\Omega\right) = M(v).
\]
It follows that
\[
\inf_{v\in\Omega}\left[\ell^{\Gamma,+}+{\rm n}^{\Gamma,+}_{v}\left(1-\mathrm{e}^{-\zeta}\right)\right]\geq 1
\]
as required.
\smallskip

To show that
 ${\rm n}^{\Gamma,+}_{v}(\zeta)<+\infty$ for every $v\in\Omega$, we note that
 \begin{equation}
 {\rm n}^{\Gamma,+}_{v}(\zeta) = \frac{1}{M(v)}{\rm n}_v^{ +}\left(\zeta{\rm e}^{p|\epsilon_\zeta|}M(\nu_\zeta); \zeta< \texttt{k}^\Omega\right)\leq C\frac{1}{M(v)}{\rm n}_v^{ +}\left(\zeta{\rm e}^{p|\epsilon_\zeta|}\right)
 \label{meanexlengthfinite}
 \end{equation}
 for some constant $C>0$. Note that ${\rm n}_v^{p, +}(\cdot): = {\rm n}_v^{ +}\left({\rm e}^{p|\epsilon_\zeta|};\,\cdot\right)$ is the excursion measure at the maximum of the MAP $((\xi,\Theta),\bP^{p})$ that results by changing measure with respect to the law of $((\xi,\Theta),\bP)$ via the martingale ${\rm e}^{p\xi_t +\Psi_\xi(-{\rm i}p)t}$. Note the latter martingale is well defined thanks to the analytic extension of \eqref{a} to a moment generating cumulant on $(-d,\alpha)$ and the fact that $p\in(0,\alpha)$. On account of the fact that $\xi$ drifts to $+\infty$ under $\bP$, i.e., its mean is strictly positive, the same is true for $\xi$ under $\bP^p$. It is a general fact
 {that for a L\'{e}vy process which drifts to $+\infty$,}
 the inverse local time at the maximum has finite mean.
 This implies that ${\rm n}^{p,+}_v(\zeta)<\infty$, which, together with \eqref{meanexlengthfinite}, implies that ${\rm n}^{\Gamma,+}_{v}(\zeta)<\infty$ for all $v\in\Omega$.
\end{itemize}

{\begin{remark}\rm\label{d>2}
The reader will note that, in light of Theorem \ref{a5a7}, the criteria (a5)'
 and (a7)' would have been equally easy to check.
 We also note that as soon as the regular Feller property can be verified in (a2) for the angular process $\Theta$ in dimension $d\geq3$, then the verification above also works in that setting. It is worth emphasising that the reason why the case $d = 2$ is more tractable here is that $\Theta = \exp({\rm i}\vartheta)$ is such that $\vartheta$ is a L\'evy process. This fact is not true in higher dimensions. Nonetheless, on account of the fact that we know the jump rate of $\Theta$, see e.g. \cite{KALEA}, it is likely that, with a little effort, one can develop an argument along the lines of the $d=2$ case given in (a2) above.
 \end{remark}}

\subsection{Open problems}

The examples given in Sections \ref{bcone} and \ref{scone} have one thing in common, namely that they originate from isotropic processes (Brownian motion and isotropic stable processes respectively).
{One major open problem that we can register here for future inspiration pertains to the obvious removal of the assumption that the underlying process is isotropic.}
\smallskip

In the setting of Brownian motion, some foundations already exist in that direction in \cite{garbit2014}. However, in the stable setting, very little is known of these anisotropic L\'evy processes  beyond some formalities, e.g. that, up to a multiplicative constant, their L\'evy measures can be written in generalised polar coordinate form
$
r^{-(\alpha +1)}\varpi(\d \theta)$,  $r>0$, $\theta\in\mathbb{S}^{d-1},
$
for some (anisotropic) measure  $\varpi$ on $\mathbb{S}^{d-1}$. A straightforward example thereof could be seen as e.g. $\varpi(\d\theta) = \Sigma(\d\theta)\mathbf{1}_{\{\theta\in \Omega\}}$, where $\Sigma$ is the surface measure on $\mathbb{S}^{d-1}$ and
$\Omega$ is a non-empty open convex domain of $\mathbb{S}^{d-1}$. Another example takes the form $\varpi(\d\theta) = h(\theta)\Sigma(\d \theta)$ for some anisotropic function $h: \mathbb{S}^{d-1}\to (0,\infty)$.
\smallskip

On account of the fact that all of the aforementioned processes are L\'evy processes, being able to issue them from the origin, their associated weak continuity in the point of issue and their Feller property is automatic.
However,
an outstanding challenge, in light of the results of this article, is to understand how to condition them to remain in a cone and show that the apex can be included in the state space, in the sense discussed in this article.
\smallskip

Whereas understanding conditioning boils down to a matter of considering the existence of harmonic functions of the cone, as we have seen in Sections \ref{bcone} and \ref{scone}, the matter of including the apex as an entrance point boils down to an understanding of how such harmonic functions interact with the underlying MAP of the anisotropic Brownian motion or L\'evy process.

\subsection{Remarks on the fluctuation theory of MAPs}
The calculations we have made that concern general MAPs have, on the one hand, been guided by the particular application to constructing the entrance law at the origin of a general ssMp. However, it is also notable that most of what has been developed here aligns with  fluctuation theory for L\'evy processes that underpins a large body of applied probability literature. There is an existing body of literature which considers applications of MAPs in the more basic setting (relative to the context in this article) that the modulator is an ergodic Markov chain with a finite number of states.  Among the many applications, this includes (with an example item of literature from the many): Aspects of queuing theory \cite{AsmussenQueue}; Fluid queues and dams \cite{Rogers, Prabhu};  Ruin problems for surplus risk processes \cite{AA, FS};  Optimal stopping problems \cite{Hongler}; Stochastic differential equations and stochastic control \cite{LY}; Multi-type branching and fragmentation processes \cite{Robin, Bertfrag}. \smallskip
The more exotic setting of a general Markov  modulator for the MAP  opens the door to much richer categories of  models with far more subtle questions. Whilst we have provided some core results in this paper for general MAPs, it is remarkable that there is a significant amount of material that is still missing from the literature in the general setting. The papers \cite{Cinlar1, Cinlar2, Kaspi83, KaspiWH} seem to be some of the very few general treatments of MAPs. As such, the variety of results proved here for MAPs lends weight to the perspective that it is a relatively tractable theory which should now be better understood since the theory of one-dimensional L\'evy processes has largely been resolved.

\section{Construction of entrance law}\label{sec:entrance}

We define the killed process $(\xi^\dagger,\Theta^\dagger)$ by setting
  \[
    (\xi^\dagger_t,\Theta^\dagger_t) :=\begin{cases}
      (\xi_t,\Theta_t) & \text{ if }t < \tau^{+}_{0}\\
      \partial & \text{ if }t \geq \tau^{+}_{0}.
    \end{cases}
  \]
  The next lemma is the analogue of Hunt's switching identity (see \cite[Theorem II.5]{Bertoin} for the case of L\'evy processes). It follows from the proof of \cite[Theorem(11.3)]{GS}, we include it here for completeness.

  \begin{lemma}\label{lem:hunt switch}
    $((\xi^\dagger,\Theta^\dagger),\hat{\bP})$ and $((\xi^\dagger,\Theta^\dagger),\bP)$ are dual with respect to ${\rm Leb} \otimes \pi$.
  \end{lemma}
  \proof
    Let $\mu:={\rm Leb}\otimes \pi$ and fix an arbitrary $t >0$. Then from Proposition \ref{lem:easy switch} and Lemma \ref{lem:walsh_fixed_time} we see that the process $((\xi_{(t-s)-},\Theta_{(t-s)-})_{s \leq t }, \bP_\mu)$ has the same law as $((\xi_s,\Theta_s)_{s \leq t },\hat{\bP}_\mu)$. It follows that the triple process $((\xi_{(t-s)-},\Theta_{(t-s)-},\bar\xi_{t})_{s \leq t }, \bP_\mu)$ has the same law as $((\xi_s,\Theta_s,\bar\xi_t)_{s\leq t },\hat{\bP}_\mu)$. Thus for any nonnegative measurable functions $f,g: \R\times \s \rightarrow \R^{+}$,
    \begin{align*}
      \int_{\R\times\s }\mu({\rm d}y,{\rm d}\theta)g(y,\theta)\hat{\bP}_{y,\theta}[f(\xi^\dagger_t,\Theta^\dagger_t)] &= \int_{\R\times\s }\mu({\rm d}y,{\rm d}\theta)\hat{\bP}_{y,\theta}[g(\xi_0,\Theta_0)f(\xi_t,\Theta_t)\1_{\{\bar\xi_t\le 0\}}]\\
      &=\int_{\R\times\s }\mu({\rm d}y,{\rm d}\theta)\bP_{y,\theta}[g(\xi_{t-},\Theta_{t-})f(\xi_0,\Theta_0)\1_{\{\bar\xi_t\le 0\}}]\\
      &=\int_{\R\times\s }\mu({\rm d}y,{\rm d}\theta)\bP_{y,\theta}[g(\xi_{t},\Theta_{t})f(\xi_0,\Theta_0)\1_{\{\bar\xi_t\le 0\}}]\\
      &= \int_{\R\times \s }\mu({\rm d}y,{\rm d}\theta)f(y,\theta)\bP_{y,\theta}[g(\xi^\dagger_t,\Theta^\dagger_t)],
    \end{align*}
    where in the third equality we have used the quasi-left continuity of $((\xi,\Theta),\bP)$.\qed

\smallskip

Recall the definition of $\varphi$ from \eqref{eq:time_change}. Let us define the time-changed process $(\xi^\varphi,\Theta^\varphi)$ by setting
$$(\xi^\varphi_t,\Theta^\varphi_t) := (\xi_{\varphi(t)},\Theta_{\varphi(t)}) \qquad \forall\ 0\le t <\bar \zeta,$$
  where $\bar\zeta := \int_0^\infty \exp\{\alpha\xi_u\}\,{\rm d}u$
  is the lifetime of  $(\xi^\varphi,\Theta^\varphi)$.
  We denote by $(\xi^{\varphi,\dagger},\Theta^{\varphi,\dagger})$ the process of $(\xi^\varphi,\Theta^\varphi)$ killed after the time $\tau^{\varphi,+}_0:=\inf\{t \geq 0:\xi^\varphi_t>0\}$.

  \begin{lemma}\label{lem:time change dual}
    The processes $((\xi^{\varphi,\dagger},\Theta^{\varphi,\dagger}),\bP)$ and $((\xi^\varphi,\Theta^\varphi),\hat\bP^{\downarrow})$ are dual with respect to the measure
    \[
      \nu_{0}({\rm d}y, {\rm d}\theta):= \1_{\{y<0\}}\frac{c_{\pi^{+}}}{\mu^{+}}\,{\rm e}^{\alpha y}\hat{H}^{+}_\theta(y){\rm d}y \pi({\rm d}\theta).
    \]
  \end{lemma}
  \proof
    Let $f, g: \R\times \s \rightarrow\R^{+}$ be two nonnegative measurable functions.
    By the definition of $\hat\bP^\downarrow$ given in Section \ref{sec:condi-negative} we have
    \begin{align}\label{lem:time change dual.0}
    &\int_{(-\infty,0)\times\s } {\rm d}y \pi({\rm d}\theta)\, \hat{H}^{+}_\theta(y) g(y,\theta) \hat\bP^{\downarrow}_{y,\theta}[f(\xi_t,\Theta_t)]\nonumber\\
      &=\int_{(-\infty,0)\times\s }{\rm d}y \pi({\rm d}\theta) \hat{H}^{+}_\theta(y) g(y,\theta)\hat\bP_{y,\theta}\left[f(\xi_t,\Theta_t)\frac{\hat{H}^{+}_{\Theta_t}(\xi_t)}{\hat{H}^{+}_\theta (y)}; t<\tau^{+}_{0}\right]\nonumber\\
      &=\int_{(-\infty,0)\times\s } {\rm d}y \pi({\rm d}\theta) g(y,\theta)\hat\bP_{y,\theta}\left[f(\xi_t,\Theta_t)\hat{H}^{+}_{\Theta_t}(\xi_t); t<\tau^{+}_{0}\right]\nonumber\\
      &=\int_{(-\infty,0)\times\s } {\rm d}y \pi({\rm d}\theta) \hat{H}^{+}_\theta(y) f(y,\theta) \bP_{y,\theta}\left[g(\xi_t,\Theta_t); t<\tau^{+}_{0}\right].
    \end{align}
  In the final equality we have applied Lemma \ref{lem:hunt switch}. The above equations show that $((\xi^\dagger,\Theta^\dagger),\bP)$ and $((\xi,\Theta),\hat\bP^\downarrow)$ are dual with respect to the measure
  \[
    \mu({\rm d}y, {\rm d}\theta) := \1_{\{y < 0\}}\frac{c_{\pi^{+}}}{\mu^{+}}\hat{H}^{+}_\theta(y){\rm d}y \pi({\rm d}\theta).
  \]
Next for $t \geq 0$, define
  \[
    A_t:=\int_0^t \exp\{\alpha \xi_u\}\,{\rm d}u.
  \]
  Then $A_t$ is an additive functional in the sense that
  \[
    A_{t+s}=A_t + A'_t \circ \theta_s\qquad t,s \geq 0
  \]
  where $\theta$ is the shift operator and $A'$ is an independent copy of $A$. Since $\varphi$ is the right inverse of $A$, \cite[Theorem 4.5]{Walsh} states that the time-changed processes $((\xi^{\varphi,\dagger},\Theta^{\varphi,\dagger}),\bP)$ and $((\xi^\varphi,\Theta^\varphi),\hat\bP^{\downarrow})$ are dual with respect to the Revuz measure $\nu$ associated with $A_{t}$, which is determined by the following formula:
  \begin{equation}
  \int_{\R\times\s }f(y,\theta)\nu({\rm d}y,{\rm d}\theta)=\lim_{t\to 0+}\frac{1}{t}\int_{\R\times\s }\mu({\rm d}z,{\rm d}v)\bP_{z,v}\left[\int_{0}^{t}f(\xi^{\dagger}_{s},\Theta^{\dagger}_{s})dA_{s}\right]\label{lem:time change dual.1}
  \end{equation}
for every nonnegative measurable functions $f:\R\times\s \to \R^{+}$.
By Fubini's theorem and the duality relation obtained in \eqref{lem:time change dual.0} we have
\begin{eqnarray}
\mbox{RHS of }\eqref{lem:time change dual.1}&=&\lim_{t\to 0+}\frac{1}{t}\int_{\R\times\s }\mu({\rm d}z,{\rm d}v)\bP_{z,v}\left[\int_{0}^{t}f(\xi^{\dagger}_{s},\Theta^{\dagger}_{s})\mathrm{e}^{\alpha \xi^{\dagger}_{s}}{\rm d}s\right]\nonumber\\
&=&\lim_{t\to 0+}\int_{\R\times \s }\mu({\rm d}z,{\rm d}v)\mathrm{e}^{\alpha z}f(z,v)\frac{1}{t}\int_{0}^{t}\hat{\bP}^{\downarrow}_{z,v}(s<\zeta){\rm d}s\nonumber\\
&=&\int_{\R\times\s }\mu({\rm d}z,{\rm d}v)\mathrm{e}^{\alpha z} f(z,v).\nonumber
\end{eqnarray}
In the final equality we use the dominated convergence theorem.
Hence the processes $((\xi^{\varphi,\dagger},\Theta^{\varphi,\dagger}),\bP)$ and $((\xi^{\varphi},\xi^{\varphi}),\hat{\bP}^{\downarrow})$ are dual with respect to ${\rm e}^{\alpha y}\mu({\rm d}y, {\rm d}\theta)=\1_{\{y<0\}}\frac{c_{\pi^{+}}}{\mu^{+}}\,{\rm e}^{\alpha y}\hat{H}^{+}_\theta(y){\rm d}y \pi({\rm d}\theta)$.\qed

\smallskip

Now we wish to apply Lemma \ref{lem:walsh_killing_time} to the dual processes $((\xi^{\varphi,\dagger},\Theta^{\varphi,\dagger}),\bP)$ and $((\xi^{\varphi},\Theta^{\varphi}),\hat{\bP}^{\downarrow})$.
In order to do so, we need to check the integral condition given in Lemma \ref{lem:walsh_killing_time}.
We will show the integral condition in Lemma \ref{lem:walsh_killing_time} by breaking it up into two lemmas as follows.

\begin{lemma}\label{lem:eta1}
For every nonnegative measurable function $f:\R\times\s \rightarrow \R^{+}$,
  \[
    \int_{\R\times \s }\rho^{\oplus}_{1}({\rm d}y,{\rm d}\theta) \hat\bP^{\downarrow}_{y,\theta}\left[\int_0^{\bar\zeta}f(\xi^{\varphi}_t,\Theta^{\varphi}_t)\, {\rm d}t\right]=\int_{\R\times \s } \nu_{0}({\rm d}y,{\rm d}\theta)f(y,\theta)\bP_{y,\theta}(\xi_{\tau^{+}_{0}}>0).
  \]
\end{lemma}
\proof
  Let $f:\R\times \s  \rightarrow \R$ be an arbitrary nonnegative measurable function. We have
  \begin{align}
  &\int_{\R\times \s }\rho^{\oplus}_{1}({\rm d}y,{\rm d}\theta) \hat\bP^{\downarrow}_{y,\theta}\left[\int_0^{\bar\zeta}f(\xi^{\varphi}_t,\Theta^{\varphi}_t)\, {\rm d}t\right]\nonumber\\
  &=\frac{c_{\pi^{+}}}{\mu^{+}}\int_{(-\infty,0)\times \s }{\rm d}y {\rm e}^{\alpha y}\hat{H}^{+}_\theta(y)\pi({\rm d}\theta) {\rm e}^{-\alpha y}\bar\Pi_\theta(-y) \hat\bP^{\downarrow}_{y,\theta}\left[\int_0^{\bar\zeta}f(\xi^{\varphi}_t,\Theta^{\varphi}_t)\, {\rm d}t\right]\nonumber\\
  &=\int_{\R\times \s }\nu_{0}({\rm d}y,{\rm d}\theta)\mathrm{e}^{-\alpha y}\bar{\Pi}_{\theta}(-y)\hat\bP^{\downarrow}_{y,\theta}\left[\int_0^{\bar\zeta}f(\xi^{\varphi}_t,\Theta^{\varphi}_t)\, {\rm d}t\right]\nonumber\\
  &=\int_{\R\times \s }\nu_{0}({\rm d}y,{\rm d}\theta)f(y,\theta)\bP_{y,\theta}\left[\int_{0}^{\tau^{\varphi,+}_{0}}\mathrm{e}^{-\alpha \xi^{\varphi}_{t}}
  \bar{\Pi}_{\Theta^{\varphi}_{t}}(-\xi^{\varphi}_{t}){\rm d}t\right],\nonumber
  \end{align}
where $\bar{\Pi}_{v}(z)=\Pi(v,\s,(z,+\infty))$.
The last equality follows from Lemma \ref{lem:time change dual}.
We undo the time-change and write
  \[
    \bP_{y,\theta}\left[\int_0^{\tau_0^{\varphi,+}}{\rm e}^{-\alpha \xi^{\varphi}_t}\bar\Pi_{\Theta^{\varphi}_t}(-\xi^{\varphi}_t)\, {\rm d}t\right] = \bP_{y,\theta}\left[\int_0^{\tau_0^+}\bar\Pi_{\Theta_t}(-\xi_t)\, {\rm d}t\right].
  \]
Hence we get
\begin{equation}\label{lem4.4.0}
\int_{\R\times \s }\rho^{\oplus}_{1}({\rm d}y,{\rm d}\theta) \hat\bP^{\downarrow}_{y,\theta}\left[\int_0^{\bar\zeta}f(\xi^{\varphi}_t,\Theta^{\varphi}_t)\, {\rm d}t\right]=\int_{\R\times \s }\nu_{0}({\rm d}y,{\rm d}\theta)f(y,\theta)\bP_{y,\theta}\left[\int_0^{\tau_0^+}\bar\Pi_{\Theta_t}(-\xi_t)\, {\rm d}t\right].
\end{equation}
On the other hand, by the L\'{e}vy system representation given in \eqref{levysys}, we have
 \begin{eqnarray}
 \bP_{y,\theta}\left(\xi_{\tau^{+}_{0}}>0\right)
 &=&\bP_{y,\theta}\left(\sum_{s\le \tau^{+}_{0}}\1_{\{\xi_{s}>0\}}\right)\nonumber\\
 &=&\bP_{y,\theta}\left[\int_{0}^{\tau^{+}_{0}}{\rm d}s\int_{\s \times \R}\1_{\{\xi_{s}+z>0\}}\Pi(\Theta_{s},{\rm d}v, {\rm d}z)\right]\nonumber\\
 &=&\bP_{y,\theta}\left[\int_{0}^{\tau^{+}_{0}}\Pi(\Theta_{s},\ \s,\ (-\xi_{s},+\infty)){\rm d}s\right]\nonumber\\
 &=&\bP_{y,\theta}\left[\int_{0}^{\tau^{+}_{0}}\bar{\Pi}_{\Theta_{s}}(-\xi_{s}){\rm d}s\right].\label{lem4.4.1}
 \end{eqnarray}
The lemma now follows by plugging \eqref{lem4.4.1} into the right-hand side of \eqref{lem4.4.0}.\qed

\smallskip

\begin{lemma}\label{lem:eta2}
 For every nonnegative measurable function $f:\R\times\s  \rightarrow \R^{+}$,
  \begin{equation}\label{lem:eta2.1}
    \int_{\R\times \s }\rho^{\oplus}_{2}({\rm d}r,{\rm d}\theta) \hat\bP^{\downarrow}_{r,\theta}\left[\int_0^{\bar\zeta}f(\xi^{\varphi}_t,\Theta^{\varphi}_t)\, {\rm d}t\right]=\int_{\R\times \s } f(r,\theta)\nu_{0}({\rm d}r,{\rm d}\theta)\bP_{r,\theta}(\xi_{\tau^{+}_{0}}=0).
  \end{equation}
\end{lemma}

\proof Without loss of generality we assume that $f$ is a nonnegative compactly supported function for which the integral on the right-hand side of \eqref{lem:eta2.1} is finite. First we undo the time-change and write
\begin{equation}\nonumber
\hat{\bP}^{\downarrow}_{r,\theta}\left[\int_{0}^{\bar{\zeta}}f(\xi^{\varphi}_{t},\Theta^{\varphi}_{t}){\rm d}t\right]
=\hat{\bP}^{\downarrow}_{r,\theta}\left[\int_{0}^{\zeta}\mathrm{e}^{\alpha \xi_{t}}f(\xi_{t},\Theta_{t}){\rm d}t\right].
\end{equation}
Let $F(x,\theta):=\mathrm{e}^{\alpha x}\hat{H}^{+}_{\theta}(x)f(x,\theta)$ for $(x,\theta)\in \R\times \s$. By \eqref{def:condi at 0} and Fubini's theorem we have
\begin{equation}\nonumber
\hat{\bP}^{\downarrow}_{0,\theta}\left[\int_{0}^{\zeta}\mathrm{e}^{\alpha \xi_{t}}f(\xi_{t},\Theta_{t}){\rm d}t\right]
=\hat{\bP}^{\downarrow}_{0,\theta}\left[\int_{0}^{\zeta}\hat{H}^{+}_{\Theta_{t}}(\xi_{t})^{-1}F(\xi_{t},\Theta_{t}){\rm d}t\right]
=\frac{\hat{{\rm n}}^{+}_{\theta}\left[\int_{0}^{\zeta}F(-\epsilon_{s},\nu_{s}){\rm d}s\right]}{\hat{{\rm n}}^{+}_{\theta}\left(\zeta=+\infty\right)}.
\end{equation}
Hence by the definition of $\rho^{\oplus}_{2}$ we get
\begin{align}
&\int_{\R\times \s }\rho^{\oplus}_{2}({\rm d}r,{\rm d}\theta) \hat\bP^{\downarrow}_{r,\theta}\left[\int_0^{\bar\zeta}f(\xi^{\varphi}_t,\Theta^{\varphi}_t)\, {\rm d}t\right]\nonumber\\
&=\frac{c_{\pi^{+}}}{\mu^{+}}\int_{\s}\hat{U}^{+}_{\pi}({\rm d}\theta,\R^{+})\frac{a^{+}(\theta)\hat{{\rm n}}^{+}_{\theta}\left(\zeta=+\infty\right)}{\ell^+(\theta)
+{\rm n}^{+}_{\theta}(\zeta)}\hat{\bP}^{\downarrow}_{0,\theta}\left[\int_{0}^{\zeta}\mathrm{e}^{\alpha \xi_{t}}f(\xi_{t},\Theta_{t}){\rm d}t\right]\nonumber\\
&=\frac{c_{\pi^{+}}}{\mu^{+}}\int_{\s}\hat{U}^{+}_{\pi}({\rm d}\theta,\R^{+})\frac{a^{+}(\theta)\hat{{\rm n}}^{+}_{\theta}\left[\int_{0}^{\zeta}
F(-\epsilon_{s},\nu_{s}){\rm d}s\right]}{\ell^+(\theta)+{\rm n}^{+}_{\theta}(\zeta)}.\label{lem7.4.1}
\end{align}
On the other hand by Proposition \ref{prop:on creeping} and Fubini's theorem we have
\begin{align}
&\int_{\R\times \s }\nu_{0}({\rm d}y,{\rm d}\theta)f(y,\theta)\bP_{y,\theta}\left(\xi_{\tau^{+}_{0}}=0\right)\nonumber\\
&=\frac{c_{\pi^{+}}}{\mu^{+}}\int_{\R^{-}\times \s }{\rm d}y \pi({\rm d}\theta)\mathrm{e}^{\alpha y}\hat{H}^{+}_{\theta}(y)f(y,\theta)\bP_{0,\theta}\left(\xi_{\tau^{+}_{-y}}=-y\right)\nonumber\\
&=\frac{c_{\pi^{+}}}{\mu^{+}}\int_{\R^{+}\times \s }{\rm d}z \pi({\rm d}\theta)\mathrm{e}^{-\alpha z}\hat{H}^{+}_{\theta}(-z)f(-z,\theta)\int_{\s }a^{+}(v)u^{+}_{\theta}({\rm d}v,z)\nonumber\\
&=\frac{c_{\pi^{+}}}{\mu^{+}}\int_{\s }\pi({\rm d}\theta)\int_{\s \times \R^{+}}U^{+}_{\theta}({\rm d}v,{\rm d}z)F(-z,\theta)a^{+}(v).\nonumber
\end{align}
From this and \eqref{lem7.4.1} we can see that to show \eqref{lem:eta2.1}, it suffices to show
\begin{equation}\label{lem:eta2.2}
\int_{\s}\hat{U}^{+}_{\pi}({\rm d}\theta,\R^{+})\frac{a^{+}(\theta)\hat{{\rm n}}^{+}_{\theta}\left[\int_{0}^{\zeta}
F(-\epsilon_{s},\nu_{s}){\rm d}s\right]}{\ell^+(\theta)+{\rm n}^{+}_{\theta}(\zeta)}=\int_{\s }\pi({\rm d}\theta)\int_{\s \times \R^{+}}U^{+}_{\theta}({\rm d}v,{\rm d}z)F(-z,\theta)a^{+}(v).
\end{equation}
By Proposition \ref{prop:equal dist} the following equation holds for all $q>0$:
\begin{equation}\label{lem:eta2.3}
\bP_{0,\pi}\left[\rme^{q\bar{g}_{\mathbf{e}_{q}}}\frac{F(-\bar{\xi}_{\mathbf{e}_{q}},\Theta_{0})a^{+}(\bar{\Theta}_{\mathbf{e}_{q}})}{q\left(\ell^+(\bar{\Theta}_{\mathbf{e}_{q}})+{\rm n}^{+}_{\bar{\Theta}_{\mathbf{e}_{q}}}(\zeta)\right)}\right]
=\hat{\bP}_{0,\pi}\left[\rme^{q\left(\mathbf{e}_{q}-\bar{g}_{\mathbf{e}_{q}}\right)}\frac{F(-(\bar{\xi}_{\mathbf{e}_{q}}-\xi_{\mathbf{e}_{q}}),\Theta_{\mathbf{e}_{q}})a^{+}(\bar{\Theta}_{\mathbf{e}_{q}})}{q\left(\ell^+(\bar{\Theta}_{\mathbf{e}_{q}})+{\rm n}^{+}_{\bar{\Theta}_{\mathbf{e}_{q}}}(\zeta)\right)}\right].
\end{equation}
By Proposition \ref{prop:wiener_hopf}, the expectation on the left-hand side equals
\begin{align}
\int_{\s}&\pi({\rm d}\theta)\int_{\s\times \R^{+}}U^{+}_{\theta}({\rm d}v,{\rm d}z)F(-z,\theta)a^{+}(v)\frac{\ell^+(v)+{\rm n}^{+}_{v}\left(\int_{0}^{\zeta}\rme^{-q s}{\rm d}s\right)}{\ell^+(v)+{\rm n}^{+}_{v}(\zeta)}\nonumber\\
&\to\int_{\s}\pi({\rm d}\theta)\int_{\s\times \R^{+}}U^{+}_{\theta}({\rm d}v,{\rm d}z)F(-z,\theta)a^{+}(v)\label{lem:eta2.4}
\end{align}
as $q\to 0+$ by the monotone convergence theorem and condition (a8) that ${\rm n}^{+}_{v}(\zeta)<+\infty$ for all $v\in\s$. Similarly by Proposition \ref{prop:wiener_hopf} and the monotone convergence theorem, the expectation on the right-hand side of \eqref{lem:eta2.3} equals
\begin{align}
\int_{\R^{+}\times \s\times\R^{+}}&\hat{V}^{+}_{\pi}({\rm d}r,{\rm d}v,{\rm d}z)\rme^{-q r}\frac{a^{+}(v)\hat{{\rm n}}^{+}_{v}\left(\int_{0}^{\zeta}F(-\epsilon_{s},\nu_{s}){\rm d}s\right)}{\ell^+(v)+{\rm n}^{+}_{v}(\zeta)}\nonumber\\
&\to\int_{\s}\hat{U}^{+}_{\pi}({\rm d}v,\R^{+})\frac{a^{+}(v)\hat{{\rm n}}^{+}_{v}\left(\int_{0}^{\zeta}F(-\epsilon_{s},\nu_{s}){\rm d}s\right)}{\ell^+(v)+{\rm n}^{+}_{v}(\zeta)}\label{lem:eta2.5}
\end{align}
as $q\to 0+$. Hence \eqref{lem:eta2.2} follows immediately by combining \eqref{lem:eta2.3}-\eqref{lem:eta2.5}.\qed

\smallskip

Finally, we show that the process $((\xi^\varphi,\Theta^\varphi),\hat{\bP}^{\downarrow})$ has a finite lifetime.

\begin{lemma}\label{lemma:finite_killing}
For every $x\le 0$, $\theta\in \s $,
$$\hat{\bP}^{\downarrow}_{x,\theta}\left(\int_{0}^{+\infty}\mathrm{e}^{\alpha \xi_{t}}{\rm d}t<+\infty\right)=1.$$
  In particular, the lifetime $\bar{\zeta}$ of the process $((\xi^\varphi,\Theta^\varphi),\hat\bP^\downarrow_{x,\theta})$ is finite almost surely and
  $\xi^\varphi_{\bar{\zeta}-}=-\infty$ $\hat\bP^\downarrow_{x,\theta}$-a.s.
\end{lemma}
\proof Since the lifetime of the time-changed process $(\xi^{\varphi},\Theta^{\varphi})$ equals $\int_{0}^{+\infty}\mathrm{e}^{\alpha \xi_{t}}{\rm d}t$, we only need to prove the first assertion. We first consider the case where $x<0$ and $\theta\in\s $.
Recall that $\hat{\bP}^{\downarrow}_{x,\theta}$ is defined from $\hat{\bP}_{x,\theta}$ through a martingale change of measure with $W_{t}:=\hat{H}^{+}_{\Theta_{t}}(\xi_{t})\1_{\{t<\tau^{+}_{0}\}}/\hat{H}^{+}_{\theta}(x)$ being the martingale.
Since $\hat{H}^{+}_{v}(y)=\hat{\bP}_{y,v}\left(\tau^{+}_{0}=+\infty\right)\in [0,1]$, $W_{t}$ is a bounded $\hat{\bP}_{x,\theta}$-martingale and hence has an almost sure limit $W_{\infty}$ such that $W_{t}\to W_{\infty}$ in $L^{1}(\hat{\bP}_{x,\theta})$. This implies that
$\hat{\bP}^{\downarrow}_{x,\theta}(A)=\hat{\bP}_{x,\theta}\left[W_{\infty}\1_{A}\right]$ for all $A\in \mathcal{F}_{\infty}$.
Hence we get
\begin{equation}\label{lem4.11.1}
\hat{\bP}^{\downarrow}_{x,\theta}\left(\int_{0}^{+\infty}\mathrm{e}^{\alpha \xi_{t}}{\rm d}t<+\infty\right)=\hat{\bP}_{x,\theta}\left[W_{\infty}\1_{\{\int_{0}^{+\infty}\mathrm{e}^{\alpha \xi_{t}}{\rm d}t<+\infty\}}\right].
\end{equation}
It follows by Lemma \ref{lem:time reversal} that
$$\hat{\bP}_{0,\pi}\left[\sup_{s\in [0,1]}|\xi_{s}|\right]=\bP_{0,\pi}\left[\sup_{s\in [0,1]}|\xi_{s}-\xi_{1}|\right]\le 2\bP_{0,\pi}\left[\sup_{s\in [0,1]}|\xi_{s}|\right]<+\infty.$$
Hence the MAP $((\xi,\Theta),\hat{\bP})$ exhibits exactly one of the tail behaviors described in Proposition \ref{prop:tail}. We have proved in Proposition \ref{prop:martingale}(i) that $\hat{\bP}_{x,\theta}\left(\tau^{+}_{0}=+\infty\right)>0$. This together with Proposition \ref{thm:slln} implies that under $\hat{\bP}_{x,\theta}$ the ordinate $\xi_{t}$ drifts to $-\infty$ at a linear rate.
Hence we have
$$\hat{\bP}_{x,\theta}\left(\int_{0}^{+\infty}\mathrm{e}^{\alpha \xi_{t}}{\rm d}t<+\infty\right)=1.$$
By this and \eqref{lem4.11.1} we get
$$\hat{\bP}^{\downarrow}_{x,\theta}\left(\int_{0}^{+\infty}\mathrm{e}^{\alpha \xi_{t}}{\rm d}t<+\infty\right)=\hat{\bP}_{x,\theta}[W_{\infty}]=1.$$
Now we consider the case where $x=0$. We have proved in Proposition \ref{prop:cond at 0} that under $\hat{\bP}^{\downarrow}_{0,\theta}$, $\xi_{t}$ leaves $0$ instantaneously and that the process $(\xi_{t},\Theta_{t})_{t>0}$ has the same transition rates as $((\xi_{t},\Theta_{t})_{t>0},\hat{\bP}^{\downarrow}_{y,\theta})$ where $(y,\theta)\in (-\infty,0)\times \s$. By the Markov property we have
\begin{equation}\nonumber
\hat{\bP}^{\downarrow}_{0,\theta}\left(\int_{s}^{+\infty}\mathrm{e}^{\alpha \xi_{t}}{\rm d}t<+\infty\right)=
\hat{\bP}^{\downarrow}_{0,\theta}\left[\hat{\bP}^{\downarrow}_{\xi_{s},\Theta_{s}}\left(\int_{0}^{+\infty}\mathrm{e}^{\alpha \xi_{t}}{\rm d}t<+\infty\right)\right]
\end{equation}
for any $s>0$. Hence we get $\hat{\bP}^{\downarrow}_{0,\theta}\left(\int_{0}^{+\infty}\mathrm{e}^{\alpha \xi_{t}}{\rm d}t<+\infty\right)=1$.\qed

\smallskip

By Lemma \ref{lem:time change dual} the processes $((\xi^{\varphi,\dagger},\Theta^{\varphi,\dagger}),\bP)$ and $((\xi^\varphi,\Theta^\varphi),\hat\bP^{\downarrow})$ are dual with respect to $\nu_{0}$. By Proposition \ref{prop:rhooplus}, Lemma \ref{lem:eta1} and Lemma \ref{lem:eta2} one has
\begin{equation}\label{eq:rhoodot}
\int_{\R\times \s }\rho^{\oplus}({\rm d}r,{\rm d}\theta)\hat{\bP}^{\downarrow}_{r,\theta}\left[\int_{0}^{\bar{\zeta}}f(\xi^{\varphi}_{t},\Theta^{\varphi}_{t}){\rm d}t\right]
=\int_{\R\times \s }f(r,\theta)\nu_{0}({\rm d}r,{\rm d}\theta)
\end{equation}
for every nonnegative measurable function $f:\R\times \s \to \R^{+}$.
We define the time-changed reversed process $(\tilde\xi,\tilde \Theta)$ by setting
$$(\tilde\xi_t,\tilde \Theta_t):= \left(\xi^\varphi_{(\bar\zeta-t)-},\Theta^\varphi_{(\bar\zeta-t)-}\right) \quad \mbox{ for }0\le t<\bar\zeta.$$
In view of \eqref{eq:rhoodot} and Lemma \ref{lemma:finite_killing} we can apply Lemma \ref{lem:walsh_killing_time} to deduce that $((\tilde {\xi}_{t},\tilde{\Theta}_{t})_{0<t<\bar{\zeta}},\hat\bP^{\downarrow}_{\rho^{\oplus}})$ is a right continuous strong Markov process having the same transition rates as $((\xi^{\varphi,\dagger},\Theta^{\varphi,\dagger}),\bP )$.
In conclusion we have just shown the following proposition.

\begin{proposition}\label{prop:P_{0} exist}
  Let $\varrho$ be the image of the probability measure $\rho^{\oplus}$ under the map
  $\phi:\ (y,\theta)\mapsto \theta\mathrm{e}^{y}$.
  Let $\P^{\searrow}_{\varrho}$ be the law of the process $(\tilde{X}_{t}:=\mathrm{e}^{\tilde{\xi}_{t}}\tilde{\Theta}_{t})_{t<\bar{\zeta}}$ under $\hat{\bP}^{\downarrow}_{\rho^{\oplus}}$.
  Then the process $((\tilde{X}_{t})_{t<\bar{\zeta}},\P^{\searrow}_{\varrho})$ is a right continuous Markov process such that $\tilde{X}_{0}=0$ and $\tilde{X}_{t}\not=0$ for all $t>0$ $\P^{\searrow}_{\varrho}$-a.s. Moreover, $((\tilde{X}_{t})_{0<t<\bar{\zeta}},\P^{\searrow}_{\varrho})$ is a strong Markov process having the same transition rates as the self-similar Markov process $(X,\{\P_{z},z\in \mathcal{H}\})$ killed upon exiting the unit ball.
\end{proposition}

By applying the scaling property of ssMp, we can describe the law of the process killed when exiting the ball of radius $r$, for any $r>0$.
Thus we see that there exists a process $(X,\P_0)$ started at the origin such that for any $r >0$,
$((X_{t})_{t<\tau^{\ominus}_{r}},\P_0)$ is equal in law to $((r\tilde{X}_{r^{-\alpha}t})_{t<r^{\alpha}\bar{\zeta}},\P^{\searrow}_{\varrho}).$

\section{Convergence of entrance law}\label{sec:convergence}

In the following we give a convergence lemma, which gives sufficient conditions for the candidate law $\P_{0}$ defined in Section \ref{sec:entrance} to be the weak limit of $\lim_{\mathcal{H}\ni z\to 0}\P_{z}$. The idea of its proof is from \cite[Proposition 7]{DDK}. For completeness we also give details here.

\begin{lemma}\label{lemma:convergence_conditions}
Suppose $\{\mu_{n}:n\ge 0\}$ is a sequence of probability measures on $\h$ which converges weakly to $\delta_{0}$.  Then $\P_{0}=\mbox{w-}\lim_{n\to +\infty}\P_{\mu_{n}}$ in the Skorokhod space if the following two conditions are satisfied:
  \begin{itemize}
  \item[(i)]
    $\lim_{\delta\to 0}\limsup_{n\to +\infty}\P_{\mu_{n}}\left[\tau^{\ominus}_{\delta}\wedge 1\right]=0,$

  \item[(ii)]  There exists a $\Delta>0$ such that for every $\delta\in (0,\Delta)$, $(X_{\tau^{\ominus}_{\delta}},\P_{\mu_{n}}) \rightarrow (X_{\tau^{\ominus}_{\delta}},\P_0)$ in distribution as $n\to +\infty$.
  \end{itemize}

\end{lemma}

\proof Let $\mathbb D_{\R^d}$ be the space of (possibly killed) c\`{a}dl\`{a}g functions $\omega:[0,\infty)\rightarrow \R^d$, equipped with the Skorokhod topology. We work with the Prokhorov's metric $\mathrm{d}(\cdot,\cdot)$ which is compatible with the Skorokhod convergence: for $m\in \mathbb{N}$ and two paths $x,y$ in $\mathbb{D}_{\R^{d}}$, define
  $$\mathrm{d}_{m}(x,y):=\inf_{\lambda\in\Lambda_{m}}\{\sup_{t\in [0,m]}|\lambda(t)-t|\vee\sup_{t\in [0,m]}|x(t)-y\circ \lambda(t)|\},$$
  where $\Lambda_{m}$ denotes the set of strictly increasing continuous functions $\lambda:[0,m]\to\R^{+}$ with $\lambda(0)=0$, and define
  $$\mathrm{d}(x,y):=\sum_{m=1}^{+\infty}2^{-m}\left(\mathrm{d}_{m}(x,y)+\mathrm{d}_{m}(y,x)\right)\wedge 1.$$
  To prove $\P_{0}=\mbox{w-}\lim_{n\to+\infty}\P_{\mu_{n}}$ in the Skorokhod space, it suffices to prove that for an arbitrary Lipschitz continuous function $f:\mathbb{D}_{\R^{d}}\to \R$ with Lipschitz constant $\kappa>0$,
  \begin{equation}
  \lim_{n\to +\infty}\P_{\mu_{n}}\left[f(X)\right]=\P_{0}\left[f(X)\right].\label{lem6.1.1}
  \end{equation}
    We note that by Proposition \ref{prop:P_{0} exist} $\left((X_{t+\tau^{\ominus}_\delta})_{t \ge 0},\P_0\right)$ is a Markov process having the same transition rates as $(X,\{\P_{z}, z\in \h\})$.
    {In view of (a1) and condition (ii), Lemma \ref{lemn2} yields that}
    for every $\delta\in (0,\Delta)$
    $$\left((X_{t+\tau^{\ominus}_{\delta}})_{t \ge 0},\P_{\mu_{n}}\right)
    \rightarrow \left((X_{t+\tau^{\ominus}_\delta})_{t \ge 0},\P_0\right)$$
    in distribution under the Skorokhod topology as $n\to +\infty$.
    Thus by the representation theorem, there exist an appropriate probability space $(\Omega^{*},\mathcal{F}^{*},\P^{*})$ and couplings $Y^{(n)}$, $Y^{(0)}$ of the processes $(X,\P_{\mu_{n}})$ and $(X,\P_0)$, respectively, such that
    $$(Y^{(n)}_{t+\varsigma_{n}})_{t \ge 0}\rightarrow (Y^{(0)}_{t+\varsigma_{0}})_{t\ge 0}\quad\mbox{ as }n\to +\infty$$
    $\P^{*}$-almost surely in the Skorokhod space, where for $k\ge 1$, $\varsigma_{k}:=\inf\{t \ge0: \|Y^{(k)}_t\|>\delta\}$ and $\varsigma_{0}:=\inf\{t \ge0: \|Y^{(0)}_t\|>\delta\}$.
    We observe that for $n\ge 1$,
    \begin{equation}
    \mathrm{d}(Y^{(n)},Y^{(0)})\le 4\delta+2\left|\varsigma_{n}-\varsigma_{0}\right|\wedge 1 +\mathrm{d}(Y^{(n)}_{\cdot+\varsigma_{n}},Y^{(0)}_{\cdot+\varsigma_{0}}).\nonumber
    \end{equation}
    Thus by the Lipschitz continuity of $f$,
    \begin{equation}\label{lem6.1.2}
    \left|\P^{*}\left[f(Y^{(n)})\right]-\P^{*}\left[f(Y^{(0)})\right]\right|\le 4\kappa\delta+2\kappa\P^{*}\left[\left|\varsigma_{n}-\varsigma_{0}\right|\wedge 1\right] +\kappa\P^{*}\left[\mathrm{d}(Y^{(n)}_{\varsigma_{n}+\cdot},Y^{(0)}_{\varsigma_{n}+\cdot})\right].
    \end{equation}
    Obviously the third term converges to $0$ as $n\to +\infty$ by the dominated convergence theorem. Note that
    $$\P^{*}\left[\left|\varsigma_{n}-\varsigma_{0}\right|\wedge 1\right]\le \P^{*}\left[\varsigma_{n}\wedge 1\right]+\P^{*}\left[\varsigma_{0}\wedge 1\right].$$
    Condition (i) implies that
    $$\lim_{\delta\to 0}\limsup_{n\to +\infty}\P^{*}\left[\varsigma_{n}\wedge 1\right]=\lim_{\delta\to 0}\limsup_{n\to+\infty}\P_{\mu_{n}}\left[\tau^{\ominus}_{\delta}\wedge 1\right]=0,$$
    and the right continuity of $(Y^{(0)},\P^{*})$ implies that $\lim_{\delta\to 0}\P^{*}\left[\varsigma_{0}\wedge 1\right]=0$.
    Hence we get by \eqref{lem6.1.2} that $\limsup_{n\to +\infty}\left|\P^{*}\left[f(Y^{(n)})\right]-\P^{*}\left[f(Y^{(0)})\right]\right|\le 4\kappa\delta$. Hence \eqref{lem6.1.1} follows immediately by letting $\delta\to 0$.\qed

\smallskip

\begin{lemma}\label{lem6.2}
For any $\delta>0$ and any bounded continuous function $f:\h\to \R$,
$z\mapsto \P_{z}\left[\tau^{\ominus}_{\delta}\wedge 1\right]$ and $z\mapsto\P_{z}\left[f(X_{\tau^{\ominus}_{\delta}})\right]$ are continuous on $\h$.
\end{lemma}

\proof  Fix an arbitrary $\delta>0$. Suppose $z_{n},z_{\infty}\in \h$ satisfies that $\lim_{n\to+\infty}z_{n}=z_{\infty}$.
{It follows by Lemma \ref{lemn2} that $(X,\P_{z_{n}})\to (X,\P_{z_{\infty}})$ in distribution under the Skorokhod topology.}
For $n\ge 0$, let $(Y^{(n)},\P^{*})$ and $(Y ,\P^{*})$ be couplings of $(X,\P_{z_{n}})$ and $(X,\P_{z_{\infty}})$ respectively, such that $Y^{(n)}\rightarrow Y $ $\P^{*}$-a.s. in the Skorokhod topology.  Let $S:=\inf\{t\ge 0:\ \|Y_{t}\|>\delta\}$ and $\varsigma_{n}:=\inf\{t\ge 0:\ \|Y^{(n)}_{t}\|>\delta\}$ for $n\ge 0$.
Since $X$ is a sphere-exterior regular process, so is $Y$, which implies that $\|Y_{t}\|\not=\delta$ for any $t<S$ $\P^{*}$-a.s. In view of this, it follows by \cite[Theorem 13.6.4]{Whitt} that
$$(\varsigma_{n},Y^{(n)}_{\varsigma_{n}})\to (S,Y_{S})\quad \P^{*}\mbox{-a.s.}$$
as $n\to+\infty$. Consequently $\left((\tau^{\ominus}_{\delta},X_{\tau^{\ominus}_{\delta}}),\P_{z_{n}}\right)$ converges in distribution to $\left((\tau^{\ominus}_{\delta},X_{\tau^{\ominus}_{\delta}}),\P_{z_{\infty}}\right)$, and hence this lemma follows.\qed

\smallskip

   \begin{lemma}\label{lem6.3}
   For any sequence $\{z_{n}:n\ge 0\}\subset \h$ with $\lim_{n\to+\infty}z_{n}=0$, we have
   \begin{equation}\label{eq:limit tau-delta}
   \lim_{\delta\to 0}\limsup_{n\to +\infty}\P_{z_{n}}\left[\tau^{\ominus}_{\delta}\wedge 1\right]=0.
   \end{equation}
   \end{lemma}
\proof 
Without loss of generality we assume $\s$ is a compact subset of $\mathbb{S}^{d-1}$.
It suffices to prove \eqref{eq:limit tau-delta} for a sequence $\{z_{n}:n\ge 0\}$ with $\lim_{n\to +\infty}\|z_{n}\|=0$ and $\lim_{n\to +\infty}\arg(z_{n})=\theta$ for some $\theta\in\s$. We first consider the case where $\arg(z_{n})=\theta$ for $n$ sufficiently large. By Lamperti-Kiu transform one has
$$\left(\tau^{\ominus}_{\delta},\P_{x}\right)\stackrel{d}{=}\left(\int_{0}^{\tau^{+}_{\log\delta}}\mathrm{e}^{\alpha \xi_{u}}{\rm d}u,\bP_{\log\|x\|,\arg(x)}\right)\quad\forall \delta>0,x\in\h.$$
Taking expectations of both sides and using the translation invariance of $\xi$ and Fubini's theorem, we have for every $x\in \h$ with $\|x\|<\delta$,
\begin{eqnarray}\label{eq:exp_time1}
    \P_{x}[\tau^{\ominus}_\delta]
    &=&\bP_{\log \|x\|, \arg(x)}\left[\int_{0}^{\tau^{+}_{\log \delta}}\mathrm{e}^{\alpha\xi_u} {\rm d}u\right]\nonumber\\
    &=&\delta^\alpha\bP_{\log ( \|x\|/\delta), \arg(x)}\left[\int_0^{\tau^{+}_{0}}   \mathrm{e}^{\alpha\xi_u} {\rm d}u\right]\nonumber\\
    &=&\delta^\alpha \int_0^\infty{\rm d} u\, \bP_{\log ( \|x\|/\delta), \arg(x)}\left[\mathrm{e}^{-\alpha(\bar\xi_u- \xi_u)}\mathrm{e}^{\alpha\bar\xi_u} \1_{\{\bar\xi_u\le 0\}}\right]\nonumber \\
    &=&\delta^\alpha \lim_{q\downarrow 0} \frac{1}{q}\bP_{\log ( \|x\|/\delta), \arg(x)}\left[\mathrm{e}^{-\alpha(\bar\xi_{\mathbf{e}_q}- \xi_{\mathbf{e}_q})}\mathrm{e}^{\alpha\bar\xi_{\mathbf{e}_q}} \1_{\{\bar\xi_{\mathbf{e}_q}\le 0\}}\right].\label{lem6.2.1}
  \end{eqnarray}
Set $y=\log(\|x\|/\delta)<0$ and $u =\arg(x)$. By Proposition \ref{prop:wiener_hopf} and the monotone convergence theorem we have
  \begin{align}
  \frac{1}{q}&\bP_{y,u}\left[\mathrm{e}^{-\alpha(\bar\xi_{\mathbf{e}_q}- \xi_{\mathbf{e}_q})}\mathrm{e}^{\alpha\bar\xi_{\mathbf{e}_q}} \1_{\{\bar\xi_{\mathbf{e}_q}\le 0\}}\right]\nonumber\\
&=\frac{1}{q}\bP_{0,u}\left[\mathrm{e}^{-\alpha(\bar\xi_{\mathbf{e}_q}- \xi_{\mathbf{e}_q})}\mathrm{e}^{\alpha\left(\bar\xi_{\mathbf{e}_q}-|y|\right)} \1_{\{\bar\xi_{\mathbf{e}_q}\le |y|\}}\right]\nonumber\\
&=\int_{\R^{+}\times\s\times [0,|y|] }\mathrm{e}^{-q r}\mathrm{e}^{\alpha(z-|y|)}\left[\ell^+(v)+{\rm n}^{+}_{v}\left(\int_{0}^{\zeta}\mathrm{e}^{-q s-\alpha s}{\rm d}s\right)\right]V^{+}_{u}({\rm d}r,{\rm d}v,{\rm d}z)\nonumber\\
&\rightarrow\int_{\s \times [0,|y|]}\mathrm{e}^{-\alpha(|y|-z)}\left[\ell^+(v)+{\rm n}^{+}_{v}\left(\int_{0}^{\zeta}\mathrm{e}^{-\alpha s}{\rm d}s\right)\right]U^{+}_{u}({\rm d}v,{\rm d}z)\label{lem6.2.2}
\end{align}
as $q\downarrow 0$. It follows from \eqref{lem6.2.1} and \eqref{lem6.2.2} that
\begin{equation}\label{lem6.2.3}
\P_{z_{n}}\left[\tau^{\ominus}_{\delta}\right]=\delta^{\alpha}\int_{\s \times [0,|y_{n}|]}\mathrm{e}^{-\alpha(|y_{n}|-z)}\left[\ell^+(v)+{\rm n}^{+}_{v}\left(\int_{0}^{\zeta}\mathrm{e}^{-\alpha s}{\rm d}s\right)\right]U^{+}_{\theta}({\rm d}v,{\rm d}z)
\end{equation}
where $y_{n}=\log(\|z_{n}\|/\delta)$. Since $|y_{n}|\to +\infty$ as $n\to +\infty$, by \eqref{eq:limit of U+} the integral on the right-hand side converges to
$$\frac{1}{\alpha}\int_{\s }\left[\ell^+(v)+{\rm n}^{+}_{v}\left(\int_{0}^{\zeta}\mathrm{e}^{-\alpha s}{\rm d}s\right)\right]\pi^{+}({\rm d}v),$$
which is bounded from above
by $c_{\pi^{+}}/\alpha$.
Hence \eqref{eq:limit tau-delta} follows by letting $\delta\to 0$ in \eqref{lem6.2.3}.
For a more general sequence $\{z_{n}:n\ge 0\}$ which satisfies the conditions stated in the beginning of this proof, we set $z^{*}_{n}:=\|z_{n}\|\theta$. The above argument shows that $\lim_{\delta\to 0}\limsup_{n\to+\infty}\P_{z^{*}_{n}}\left[\tau^{\ominus}_{\delta}\wedge 1\right]=0$. Since $\lim_{n\to +\infty}\|z^{*}_{n}-z_{n}\|=0$ and by Lemma \ref{lem6.2} the function $z\mapsto \P_{z}\left[\tau^{\ominus}_{\delta}\wedge 1\right]$ is uniformly continuous on any compact subset of $\h$, we have $\lim_{n\to+\infty}\left|\P_{z^{*}_{n}}\left[\tau^{\ominus}_{\delta}\wedge 1\right]-\P_{z_{n}}\left[\tau^{\ominus}_{\delta}\wedge 1\right]\right|=0$ and hence \eqref{eq:limit tau-delta} follows.\qed

\smallskip

\begin{lemma}\label{lem:stationary overshoots for X}
Suppose $\{z_{n}:n\ge 0\}\subset \h$ satisfies $\lim_{n\to +\infty}z_{n}=0$. Then for any $\delta>0$, the probability measures $\P_{z_{n}}\left(X_{\tau^{\ominus}_{\delta}}\in\cdot\right)$
converges weakly to a proper distribution $\mu_{\delta}(\cdot)$ on $\h$.
\end{lemma}

\proof We need to show that there exists a distribution $\mu_{\delta}$ on $\h$ such that
\begin{equation}\label{lem6.4.1}
\lim_{n\to+\infty}\P_{z_{n}}\left[f(X_{\tau^{\ominus}_{\delta}})\right]=\int_{\h}f\mathrm{d}\mu_{\delta}
\end{equation}
for every bounded continuous function $f:\h\to\R$. In view of Lemma \ref{lem6.2} and the argument in the end of the above proof, we only need to prove that \eqref{lem6.4.1} holds for a sequence $\{z_{n}:n\ge 0\}$ where
$\lim_{n\to+\infty}\|z_{n}\|=0$ and $\arg(z_{n})=\theta$ for $n$ sufficiently large. By Lamperti-Kiu transform we have
\begin{eqnarray}
\P_{z_{n}}\left[f(X_{\tau^{\ominus}_{\delta}})\right]
&=&\bP_{\log\|z_{n}\|,\theta}\left[f\left(\exp\{\xi_{\tau^{+}_{\log \delta}}\}\Theta_{\tau^{+}_{\log\delta}}\right)\right]\nonumber\\
&=&\bP_{0,\theta}\left[f\left(\mathrm{e}^{\log \delta}\exp\left\{\xi_{\tau^{+}_{\log \frac{\delta}{\|z_{n}\|}}}-\log\frac{\delta}{\|z_{n}\|}\right\}\Theta_{\tau^{+}_{\log\frac{\delta}{\|z_{n}\|}}}\right)\right].\nonumber
\end{eqnarray}
Since $\|z_{n}\|\to 0$ and $\log(\delta/\|z_{n}\|)\to +\infty$, Proposition \ref{prop:stationary overshoots} yields that the distribution of $(\xi_{\tau^{+}_{\log(\delta/\|z_{n}\|)}}-\log(\delta/\|z_{n}\|),\Theta_{\tau^{+}_{\log(\delta/\|z_{n}\|)}})$ converges weakly to $\rho^{\ominus}$. Thus the 
{expectation}
on the right-hand side of the above equation converges to $\int_{\R^{+} \times \s}f\left(\mathrm{e}^{\log\delta}\mathrm{e}^{z}v\right)\rho^{\ominus}({\rm d}z,{\rm d}v)$.
Hence, by setting 
$$\mu_{\delta}(\cdot)=\int_{\R^{+} \times \s}\1_{\{\mathrm{e}^{\log \delta}\mathrm{e}^{z}v\in \cdot\}}\rho^{\ominus}({\rm d}z,{\rm d}v),$$ we get \eqref{lem6.4.1}.\qed

\smallskip

\begin{lemma}\label{lem:6.5}
For any $\delta>0$, we have $\P_{0}\left(X_{\tau^{\ominus}_{\delta}}\in\cdot\right)=\mu_{\delta}(\cdot)$.
\end{lemma}
\proof Suppose $f:\h\to \R$ is an arbitrary bounded continuous function and $\sigma_{n}:=1/n$ for $n\ge 1$. By
{the strong Markov property,}
we have for any $0<\sigma_{n}<\delta$,
\begin{equation}
\P_{0}\left[f(X_{\tau^{\ominus}_{\delta}})\right]=\P_{0}\left[\P_{X_{\tau^{\ominus}_{\sigma_{n}}}}\left[f(X_{\tau^{\ominus}_{\delta}})\right]\right]
=\P_{0}\left[g(X_{\tau^{\ominus}_{\sigma_{n}}})\right]\label{lem6.5.1}
\end{equation}
where $g(x):=\P_{x}\left[f(X_{\tau^{\ominus}_{\delta}})\right]$. Since under $\P_{0}$ the process $X_{t}$ leaves $0$ instantaneously and continuously,
we have $X_{\tau^{\ominus}_{\sigma_{n}}}\to 0$ $\P_{0}$-a.s. as $n\to +\infty$. Hence by Lemma \ref{lem:stationary overshoots for X}, $g(X_{\tau^{\ominus}_{\sigma_{n}}})=
{\P}
_{X_{\tau^{\ominus}_{\sigma_{n}}}}\left[f(X_{\tau^{\ominus}_{\delta}})\right]\to \mu_{\delta}(f)$ $\P_{0}$-a.s. By letting $n\to+\infty$ in \eqref{lem6.5.1} we get that $\P_{0}\left[f(X_{\tau^{\ominus}_{\delta}})\right]= \mu_{\delta}(f)$, which yields this lemma.\qed

\bigskip

\noindent\textit{Proof of Theorem \ref{thm:main}:}
 The statements of
 (C1), (C2) and (C3)
 are from Propositions \ref{prop:cond <0}-\ref{prop:cond at 0}, Proposition \ref{prop:stationary overshoots} and Proposition \ref{prop:P_{0} exist}, respectively. Hence we only need to show
  (C4) and (C5).

 \smallskip

\noindent(C4): We get $\P_{0}=\mbox{w-}\lim_{\h\ni z\to 0}\P_{z}$ by a combination of Lemmas \ref{lemma:convergence_conditions}-\ref{lem:6.5}. Properties \eqref{p:distribution} and \eqref{p:instant} are direct consequences of the construction of $(X,\P_{0})$ given in Section \ref{sec:entrance}.
Next we show that $(X,\{\P_{z},z\in \h_{0}\})$ is a Feller process.
{Let $\h^{0}_{\partial}=\h_{0}\cup\{\partial\}$ (resp. $\h_{\triangle}=\h\cup\{\triangle\}$) be the one-point compactification of $\h_{0}$ (resp. $\h$). Both $\h^{0}_{\partial}$ and $\h_{\triangle}$ are compact separable metric spaces. Let $C_{0}(\h_{0})$ (resp. $C_{0}(\h)$) be the class of continuous functions on $\h^{0}_{\partial}$ (resp. $\h_{\triangle}$) vanishing at $\partial$ (resp. $\triangle$).}
 Fix an arbitrary $f\in C_{0}(\h_{0})$, and let $P_{t}f(z):=\P_{z}\left[f(X_{t})\right]$ for $z\in \h_{0}$ and $t\ge 0$. To show the Feller property, it suffices to show that $P_{t}f\in C_{0}(\h_{0})$ for all $t>0$ and $\lim_{t\to 0+}P_{t}f(z)=f(z)$ for all $z\in \h_{0}$ 
. The latter holds naturally since $(X,\{\P_{z},z\in\h_{0}\})$ is a right continuous process.
We only need to show $P_{t}f\in C_{0}(\h_{0})$ for $t>0$. Suppose $x_{n},x\in\h_{0}$ and $x_{n}\to x$.
{It is proved by the above argument and Lemma \ref{lemn2} that w-$\lim_{n\to+\infty}\P_{x_{n}}=\P_{x}$ in the weak sense of measures on the Skorokhod space.}
If
\begin{equation}\label{thm4.1.1}
\P_{x}\left(X_{t-}\not=X_{t}\right)=0
\end{equation}
for $t>0$, then it follows by \cite[Proposition VI.2.1]{Jacod} that $(X_{t},\P_{x_{n}})$ converges in distribution to $(X_{t},\P_{x})$ and hence $\lim_{n\to+\infty}P_{t}f(x_{n})=\lim_{n\to+\infty}\P_{x_{n}}\left[f(X_{t})\right]=\P_{x}\left[f(X_{t})\right]=P_{t}f(x)$.
{Note that for $x\in \h$ and $t>0$,
$$\P_{x}\left(X_{t-}\not=X_{t}\right)=\bP_{\log\|x\|,\arg(x)}\left((\xi_{\varphi(t)-},\Theta_{\varphi(t)-})\not=(\xi_{\varphi(t)},\Theta_{\varphi(t)})\right),$$
where $\varphi(t)$ defined in \eqref{eq:time_change} is a stopping time with respect to the process $((\xi,\Theta),\bP)$. Hence \eqref{thm4.1.1} holds by the quasi-left continuity of $((\xi,\Theta),\bP)$.
}
 For $x=0$, we have by the Markov property that
$$\P_{0}\left(X_{t-}\not=X_{t}\right)=\P_{0}\left(\P_{X_{t/2}}\left(X_{\frac{t}{2}-}\not=X_{\frac{t}{2}}\right)\right)=0\quad \forall t>0.$$
Thus we have proved \eqref{thm4.1.1} holds for all $x\in\h_{0}$ and $t>0$. Hence $z\mapsto P_{t}f(z)$ is continuous on $\h_{0}$.
Next we show $P_{t}f$ vanishes at
{$\partial$}.
{Let $C^{*}_{0}(\h_{0})$ be the subclass of $C_{0}(\h_{0})$ vanishing at $0$. We observe that if a sequence $\{x_{n}:n\ge 1\}\subset \h$ converges to either $\partial$ or $0$ in the space $\h^{0}_{\partial}$, it also converges to $\triangle$ in the space $\h_{\triangle}$. So
a function $g$ in $C^{*}_{0}(\h_{0})$ can be viewed as a function in $C_{0}(\h)$ by setting $g(\triangle)=0$.
Thus by Lemma \ref{lemn2}, for any $\h\ni x_{n}\to \partial$ we have
\begin{equation}\label{n1}
\lim_{n\to+\infty}P_{t}g(x_{n})=\lim_{n\to+\infty}\P_{x_{n}}\left[g(X_{t})\right]=0
\end{equation}
for all $g\in C^{*}_{0}(\h_{0})$ and $t>0$.
If, in particular, we take $x_{n}=\mathrm{e}^{r_{n}}\theta$ where $(r_{n},\theta)\in\R\times \s$ and $r_{n}\to +\infty$, then by the scaling property and the bounded convergence theorem, we have
$$\lim_{n\to+\infty}P_{t}h(x_{n})=\lim_{n\to+\infty}\P_{\theta}\left[h(\mathrm{e}^{r_{n}}X_{\mathrm{e}^{-\alpha r_{n}}t})\right]=0$$
for all $h\in C_{0}(\h_{0})$ and $t>0$. This combined with \eqref{n1} implies that as $x_{n}\to\partial$ the distribution of $(X_{t},\P_{x_{n}})$ converges weakly to the Dirac measure at $\partial$.}
We use $B(0,\delta)$ to denote the $\delta$-neighborhood of $0$. It follows that
\begin{equation}\label{thm4.1.2}
\lim_{\h\ni x\to \partial}\P_{x}\left(X_{t}\in B(0,\delta)\right)=0\quad\forall \delta>0.
\end{equation}
Note that for every $x\in\h$ and $\delta>0$,
\begin{eqnarray}
|P_{t}f(x)|&\le&\left|\P_{x}\left[f(X_{t});X_{t}\in B(0,\delta)\right]\right|+\left|\P_{x}\left[f(X_{t});X_{t}\not\in B(0,\delta)\right]\right|\nonumber\\
&\le&\|f\|_{\infty}\P_{x}\left(X_{t}\in B(0,\delta)\right)+\sup_{y\in \h\setminus B(0,\delta)}|f(y)|.\nonumber
\end{eqnarray}
In view of \eqref{thm4.1.2} and the fact that $f$ vanishes at $\partial$, by letting $x\to \partial$ and then $\delta\to +\infty$ in the above inequality, we get that $\lim_{\h\ni x\to\partial}\left|P_{t}f(x)\right|=0$. Hence $P_{t}f\in C_{0}(\h_{0})$.
Therefore $(X,\{\P_{z},z\in \h_{0}\})$ is a Feller process.

\smallskip

Recall that $((X_{t})_{t>0},\P_{0})$ has the same transition rates as the ssMp $(X,\{\P_{z}, z\in \h\})$. Thus by Markov property, to show $(X,\P_{0})$ is self-similar, we only need to show that $(X_{t},\P_{0})\stackrel{d}{=}(cX_{c^{-\alpha}t},\P_{0})$ for every $t>0$ and $c>0$ , and this is true since
$$
(X_{t},\P_{0})=\mbox{w-}\lim_{\h\ni z\to 0}(X_{t},\P_{cz})\stackrel{d}{=}\mbox{w-}\lim_{\h\ni z\to 0}(cX_{c^{-\alpha}t},\P_{z})=(cX_{c^{-\alpha}t},\P_{0}).
$$

Finally we show the uniqueness of $\P_{0}$.
Suppose there exists another probability measure $\P^{*}_{0}$ for which the property \eqref{p:Feller} is satisfied. Using Feller property twice we get
$$\P^{*}_{0}\left(X_{t}\in\cdot\right)=\mbox{w-}\lim_{\h\ni z\to 0}\P_{z}\left(X_{t}\in\cdot\right)=\P_{0}\left(X_{t}\in \cdot\right)\mbox{ for every }t>0.$$
Hence by Markov property $\P^{*}_{0}$ is equal to $\P_{0}$. Suppose now that, instead, $\P^{*}_{0}$ satisfies the property \eqref{p:instant}. Then for any $t>0$
and any bounded continuous function $h:\s\to\R$,
\begin{eqnarray}
\P^{*}_{0}\left[h(X_{t})\right]
&=&\lim_{\epsilon\to 0+}\P^{*}_{0}\left[h(X_{t+\epsilon})\right]\nonumber\\
&=&\lim_{\epsilon\to 0+}\P^{*}_{0}\left[\P_{X_{\epsilon}}\left[h(X_{t})\right]\right]\nonumber\\
&=&\P_{0}\left[h(X_{t})\right].\nonumber
\end{eqnarray}
We used in the first equality the fact that $(X,\P^{*}_{0})$ is a right continuous process and in the second equality the Markov property. The fact that $\lim_{\epsilon\to 0+}X_{\epsilon}=0$ $\P^{*}_{0}$-a.s. and the Feller property of $(X,\{\P_{z},z\in\h_{0}\})$ imply that $\P_{X_{\epsilon}}\left(X_{t}\in\cdot\right)$ converges weakly to $\P_{0}\left(X_{t}\in\cdot\right)$ $\P^{*}_{0}$-a.s. This is used in third equality.
The above equation implies that $\P^{*}_{0}(X_{t}\in\cdot)=\P_{0}(X_{t}\in\cdot)$ for all $t>0$, and therefore $\P^{*}_{0}$ is equal to $\P_{0}$ again by the Markov property.

\smallskip

\noindent(C5): By the strong Markov property and the sphere-exterior regularity of $(X,\{\P_{z},z\in\h\})$, we have
\begin{eqnarray}
\P_{0}\left(\|X_{t}\|=\delta\mbox{ for some } t\in (0,\tau^{\ominus}_{\delta})\right)
&=&\P_{0}\left(\|X_{t}\|=\delta\mbox{ for some } t\in [\tau^{\ominus}_{\delta/2},\tau^{\ominus}_{\delta}),\tau^{\ominus}_{\delta/2}<\tau^{\ominus}_{\delta}\right)\nonumber\\
&=&\P_{0}\left[\P_{X_{\tau^{\ominus}_{\delta/2}}}\left(\|X_{t}\|=\delta\mbox{ for some } t<\tau^{\ominus}_{\delta})\right);\tau^{\ominus}_{\delta/2}<\tau^{\ominus}_{\delta}\right]\nonumber\\
&=&0.\nonumber
\end{eqnarray}
In view of this and  the fact that w-$\lim_{\h\ni z\to 0}\P_{z}=\P_{0}$ in the Skorokhod space, it follows by the Skorokhod representation theorem and \cite[Theorem 13.6.4]{Whitt} that
$\left((X_{\tau^{\ominus}_{\delta}-},X_{\tau^{\ominus}_{\delta}}), \P_{z}\right)$ converges in distribution to $\left((X_{\tau^{\ominus}_{\delta}-},X_{\tau^{\ominus}_{\delta}}), \P_{0}\right)$ as $z\to 0$. We note that for any $x>0$ and $\theta\in\s$
\begin{align}
&\P_{\theta\mathrm{e}^{-x}}\left(\arg(X_{\tau^{\ominus}_{1}-})\in {\rm d}v,\ \log\|X_{\tau^{\ominus}_{1}-}\|\in {\rm d}y,\ \arg(X_{\tau^{\ominus}_{1}})\in {\rm d}\phi,\ \log\|X_{\tau^{\ominus}_{1}}\|\in {\rm d}z\right)\nonumber\\
&=\bP_{-x,\theta}\left(\Theta_{\tau^{+}_{0}-}\in {\rm d}v,\ \xi_{\tau^{+}_{0}-}\in {\rm d}y,\ \Theta_{\tau^{+}_{0}}\in {\rm d}\phi,\ \xi_{\tau^{+}_{0}}\in {\rm d}z\right)\nonumber\\
&=\bP_{0,\theta}\left(\Theta_{\tau^{+}_{x}-}\in {\rm d}v,\ \xi_{\tau^{+}_{x}-}-x\in {\rm d}y,\ \Theta_{\tau^{+}_{x}}\in {\rm d}\phi,\ \xi_{\tau^{+}_{x}}-x\in {\rm d}z\right).\nonumber
\end{align}
By Proposition \ref{prop:stationary overshoots} the last distribution converges weakly to $\rho({\rm d}v,{\rm d}y,{\rm d}\phi,{\rm d}z)$ as $x\to+\infty$. Hence by the above argument we get
\begin{align}
&\mbox{w-}\lim_{\h\ni z\to 0}\P_{z}\left(\arg(X_{\tau^{\ominus}_{1}-})\in {\rm d}v,\ \log\|X_{\tau^{\ominus}_{1}-}\|\in {\rm d}y,\ \arg(X_{\tau^{\ominus}_{1}})\in {\rm d}\phi,\ \log\|X_{\tau^{\ominus}_{1}}\|\in {\rm d}z\right)\nonumber\\
&=\P_{0}\left(\arg(X_{\tau^{\ominus}_{1}-})\in {\rm d}v,\ \log\|X_{\tau^{\ominus}_{1}-}\|\in {\rm d}y,\ \arg(X_{\tau^{\ominus}_{1}})\in {\rm d}\phi,\ \log\|X_{\tau^{\ominus}_{1}}\|\in {\rm d}z\right)\nonumber\\
&=\rho({\rm d}v,{\rm d}y,{\rm d}\phi,{\rm d}z).\nonumber
\end{align}
This completes the proof.
\qed

\bibliographystyle{amsplain}

\end{document}